\documentclass[12pt,letterpaper]{article}
	\pdfoutput=1
	
	\usepackage{graphics,epsfig,graphicx,color}
	\graphicspath{{/EPSF/}{../Figures/}{Figures/}}
	
	\usepackage{amssymb,latexsym}
	
	\usepackage{setspace}
	\usepackage{amsmath}
	
	\usepackage{mathrsfs,amsfonts}
	
	\usepackage{amsthm,amsxtra}
	
	\usepackage[font={small,it}]{caption}

	\usepackage[utf8]{inputenc}

	\usepackage{graphics,epsfig,graphicx,color}
	\usepackage{subcaption}
	
	\newtheorem{remark}{Remark}[section]

	\newcommand{\beqa}{\begin{eqnarray}}
	\newcommand{\eeqa}[1]{\label{#1}\end{eqnarray}}
	\newcommand{\beq}{\begin{equation}}
	\newcommand{\eeq}[1]{\label{#1}\end{equation}}

	\newcommand{\GD}{\Delta}

	\def\Bn{{\bf n}}

	\def\Bv{{\bf v}}
	
	\def\Bx{{\bf x}}
	\def\By{{\bf y}}

	\def\B0{{\bf 0}}

	\def \RR {{\mathbb R}}
	
	\def \ba {\begin{array}}
	\def \ea {\end{array}}
	\newtheorem {Thm} {Theorem} [section]

	\newtheorem {Adef} [Thm] {Definition}
	
	\newtheorem {Arem} [Thm] {Remark}
	
	\newtheorem {Aexa} [Thm] {Example}
	
	\newtheorem {Anot} [Thm] {Notation}

	\def \refe #1.{(\ref{#1})}
	\def \reff #1.{figure~\ref{#1}}
	\def \refs #1.{section~\ref{#1}}
	\def \refss #1.{subsection~\ref{#1}}
	\def \refD #1.{Definition~\ref{#1}}
	\def \refT #1.{Theorem~\ref{#1}}
	\def \refL #1.{Lemma~\ref{#1}}
	\def \refC #1.{Corollary~\ref{#1}}
	\def \refP #1.{Proposition~\ref{#1}}
	\def \refR #1.{Remark~\ref{#1}}
	\def \refE #1.{Example~\ref{#1}}
	\def \refN #1.{Notation~\ref{#1}}
	\newif\ifPDF
	\ifx\pdfoutput\undefined
	    \PDFfalse
	\else
	    \ifnum\pdfoutput > 0
	        \PDFtrue
	    \else
	        \PDFfalse
	    \fi
	\fi
	
	\ifPDF
	    \usepackage{pdftricks}
	    \begin{psinputs}
	        \usepackage{pstricks}
	        \usepackage{pstcol}
	        \usepackage{pst-plot}
	        \usepackage{pst-tree}
	        \usepackage{pst-eps}
	        \usepackage{multido}
	        \usepackage{pst-node}
	        \usepackage{pst-eps}
	    \end{psinputs}
	\else
	    \usepackage{pstricks}
	\fi
	
	\ifPDF
	    \usepackage[debug,pdftex,colorlinks=true,
	    linkcolor=blue, bookmarksopen=false,
	    plainpages=false,pdfpagelabels]{hyperref}
	\else
	    \usepackage[dvips]{hyperref}
	\fi
	
	\usepackage{fancyhdr}
	
\usepackage{calc}
\usepackage{authblk}

\title{On the synthesis of acoustic sources with controllable near fields}
\author[1]{D. Onofrei}
\author[2]{E. Platt}
\affil[1]{University of Houston, Department of Mathematics: onofrei@math.uh.edu}
\affil[2]{University of Houston, Department of Mathematics: eplatt@math.uh.edu}

\newcommand{\figsizeC}{0.3278\textwidth}
\newcommand{\figsizeD}{0.3025\textwidth} 

\begin{document}

\maketitle

\begin{abstract}
 
 In this paper we present a strategy for the the synthesis of acoustic sources with controllable near fields in free space and finite depth homogeneous ocean environments. We first present the theoretical results at the basis of our discussion and then, to illustrate our findings we focus on the following three particular examples: 
 \begin{enumerate}
 \item{} acoustic source approximating a prescribed field pattern in a given bounded sub-region of its near field.\vspace{-0.6cm}\\
 \item{} acoustic source approximating different prescribed field patterns in given disjoint bounded near field sub-regions.\vspace{-0.6cm}\\
 \item{} acoustic source approximating a prescribed back-propagating field in a given bounded near field sub-region while maintaining a very low far field signature.\vspace{0cm}
 \end{enumerate}
  For each of these three examples, we discuss the optimization scheme used to approximate their solutions and support our claims through relevant numerical simulations.

 \end{abstract}

\section{Introduction and main results}

The problem of active control of acoustic fields is well studied in the literature with a multitude of ideas and techniques presented (see monographs \cite{1,2}). The main strategies for active sound control are based on the use of boundary controls or secondary sources.
		 	
Applications of sound field control ideas are very important and they include: active noise control \cite{5} (see also the pioneer works \cite{3,4}), acoustic field reproduction \cite{6,7,8,9,10} and active control of scattered sound fields \cite{11,12,13,14,15,16,17,18,18'}. A rigorous comparative analysis of the theoretical similarities and respective challenges for these three areas of applications is done in \cite{19}. 
        
 In a recent development in \cite{20} (see also \cite{21} for the low frequency approximation), a general analytical approach based on the theory of boundary layer potentials was proposed for the active acoustic control problem in homogeneous environments. Then, in \cite{22}, building up on \cite{20}, the authors presented a thorough two dimensional sensitivity analysis for the synthesis of time-harmonic weak radiators with controllable patterns in some exterior region and proved that such acoustic sources will be feasible only if the region of control is in the reactive near-field of the source. 
 
 The work presented in this paper uses ideas from, and is relevant to, a wide array of important research areas: acoustic wave field synthesis, inverse source problems, optimization, personal audio techniques, acoustic near field control. 
 We are making use of the theoretical results developed in \cite{20} and, through a Tikhonov regularization procedure (with Morozov discrepancy principle for the choice of the regularization parameter), we synthesize acoustic sources in one of the following scenarios:
 
 \medskip
 
 $1.$ Sources approximating a given pattern in a prescribed exterior near field sub-region.
 
 \medskip
 
 $2.$ Sources approximating a given pattern in a prescribed sub-region of their near field while having a null in a different given sub-region of their near field.
 
 \medskip
 
 $3.$ Sources which have a very weak field in a given (sufficiently far) exterior annuli while approximating a given pattern in a prescribed sub-region of their near field.
 
 \medskip
 
The first type of sources are relevant for the problem of acoustic rendering \cite{5',6}. The second type of sources above present an interest for the problem of personal audio studied in \cite{8,9,10,11} where we assume that by superposition our strategy will imply the possibility to approximate, with one source, different given sound patterns in disjoint regions of space. For the third type of sources above, although our strategy works well for the general question of synthesis of weak acoustic radiators approximating any given pattern in the near field control region, we focused on the problem of characterizing the necessary inputs (normal velocity or pressure) on the boundary of the source so that it approximates a backward propagating plane wave in the region of control while maintaining a very weak field in the given exterior annuli. This problem is relevant for the question of acoustic shielding or cloaking since by using a similar strategy we believe we can synthesize a planar array with similar properties: having a very weak field in an exterior annuli (where enemy detection measurements are taken) while approximating a given backward propagating plane wave in a near field region in front of it. Thus, by superposition, such an array could, when paired with a time control loop for the detection of interrogating signals, annihilate  through destructive interference any incoming signal in its near field region without a large signature in its far field (i.e., shielding an object located behind the array). Then, a compact volume surrounded by a similar conformal array would lead to an active cloaking device for any object located inside.

The results presented in the literature regarding pattern synthesis use arrays of secondary sources (usually approximated by point sources) to control the field in interior regions (i.e., located in the interior of the geometric convex hull of the point sources), or focus on planar rendering (i.e. control in a horizontal plane) or assume that the field to be approximated propagate away from the source to be synthesized.  

In the present paper we propose a theoretical optimization strategy for the synthesis of acoustic sources which approximate different prescribed field patterns in given disjoint {\textit{ exterior } regions in free space and finite depth homogeneous ocean environments. To simplify the exposition, for the numerical support we consider only the case of sources in free space and focus on the three particular cases listed above.


 The paper is organized as follows: In Section \ref{theoretical} we present the theoretical results in two parts: first, in Subsection \ref{fs} we first briefly recall the theoretical results of \cite{20} for acoustic control in free space and then, in Subsection \ref{ho} we discuss their extension to the problem of underwater acoustic control in the context of a constant depth homogeneous ocean environments. In Section \ref{numerical}, we build up on our previous results in \cite{22} and discuss the $L^2$ - Tikhonov regularization with Morozov discrepancy numerical approximation for the acoustic control problem in 3D and (assuming the superposition principle) without loosing the generality present numerical simulations in the three important situations listed above: first, in Subsection \ref{O} we present the synthesis of an acoustic source approximating a prescribed plane wave in a give near field sub-region; then in Subsection \ref{I} we present the synthesis of an acoustic source with a null in a given sub-region of its near field and approximating an outgoing plane wave in a disjoint near field sub-region; and finally, in Subsection \ref{II} we synthesize a very weak acoustic radiator (almost non-radiating source (ANR)) approximating, in a sub-region of its near field, a given backward propagating (propagating towards the source) plane wave. 

\section{Theoretical results}
\label{theoretical}

In this section we will present the theoretical results behind our optimization scheme presented below. In Section \ref{fs} we will recall the results of \cite{20} developed for the free space environments (i.e., homogeneous media with no boundaries and radiating condition at infinity) and then in Section \ref{ho} we will present their extension to the case of finite depth homogeneous ocean environments as introduced in \cite{23,24} (i.e., infinite rectangular waveguide with constant depth along z direction, $z\in[h,0]$ for some $h<0$, and pressure release boundary at the water-air interface $z=0$, total reflecting boundary at the ocean bottom interface $z=h$ together with radiation condition at infinity). 

We consider the source support represented by $D_a$, a compact region of space with smooth enough boundary (Lipschitz boundary will be enough) and as in \cite{20}, assume that $D_1\Subset \RR^3$ and $D_2\Subset \RR^3$, with $D_1\cap D_2=\emptyset$ and $\{D_1\cup D_2\}\cap D_a=\emptyset$. Without losing the generality we assume $u_1$ is a solution of the Helmholtz equation in neighbourhoods of $D_1$ and $u_2=0$. With the these general hypotheses, in what follows the following three geometrical situations will be considered:
\begin{eqnarray}
&\;i)\;& D_1 \mbox{ bounded }, D_2=\emptyset.\nonumber\\
&\;ii)\;& D_1 \mbox{ bounded }, D_2 \mbox{ bounded }.\label{geom-setting}\\
&\;iii)& D_1 \mbox{ bounded }, D_2 \mbox{ unbounded with } D_1\Subset\RR^3\!\setminus\!D_2.\nonumber
\end{eqnarray}
At this point we mention that the theoretical results of \cite{20} hold true for any finite number of mutually disjoint regions $D_i$, $i\in\{1,...,n\}$ satisfying $\{\cup_i D_i\}\cap D_a=\emptyset$ and respectively $n$ scalar acoustic fields $u_i$ in $D_i$. 

 In the case of free space environments the main question is to characterize boundary inputs (normal velocity or pressure) on the boundary of the source such that the acoustic field radiated by it has the property that it approximates $u_1$ in  $D_1$ and $u_2=0$ in region $D_2$ respectively (the condition on $D_2$ is not needed in case \eqref{geom-setting} $i)$ above). 

\subsection{Acoustic control in free space}
\label{fs}

In this section we will recall the result obtained in \cite{20} in the geometrical setting described at \eqref{geom-setting}: A source approximating prescribed acoustic patterns $u_1$ and $u_2$ in two given disjoint exterior regions, $D_1$ and respectively $D_2$.  Mathematically this can be written as follows:

\vspace{0.6cm}
{\bf Problem 1. }

 Find normal velocity ${v}_n$ (or pressure $p_b$) on the boundary of the antenna $\partial D_a$ so that,

\begin{equation}
\label{1}
\vspace{0.15cm}\left\{\vspace{0.15cm}\begin{array}{llll}
		 \GD u+k^2u=0 \mbox{ in }\RR^3\!\setminus \!\overline{D}_a\vspace{0.15cm},\\
		 \nabla u\cdot \Bn=v_n, (\mbox{ or } u=p_b)\mbox{ on }\partial D_a\vspace{0.15cm},\\
		 \displaystyle\left(\frac{\Bx}{|\Bx|},\nabla u(\Bx)\right)\! -\!iku(\Bx)\!=\!o\left(\frac{1}{|\Bx|}\!\right)\!,\mbox{ as }|\Bx|\rightarrow\infty
		 \mbox{ uniformly for all $\displaystyle\frac{\Bx}{|\Bx|},$ }\end{array}\right.
\end{equation}
(where $\GD$ denotes the 3D Laplace operator, $\nabla$ denotes the 3D gradient, $q=o(a)$ means $\displaystyle\lim_{a\rightarrow 0}q=0$,  $\partial S$ denotes the boundary of the set $S\Subset\RR^3$ and here $\Bn$ denotes the exterior normal to $\partial D_a$) and the following approximations hold true,
		\begin{equation}\label{1-c}u\approx u_1, \mbox{ in } D_1, \mbox{ and } u\approx 0, \mbox{ in } D_2,\end{equation}
		where the approximation in \eqref{1-c} is in the sense of smooth norms (e.g., twice differentiable functions) 
		
As a consequence of results in \cite{20} we have that {\bf Problem 1} above can be answered in the affirmative in all the geometrical configurations described at \eqref{geom-setting}. Indeed, if $k$ is not a resonance (i.e., for all wave numbers $k$ except a discrete set \cite{20}) we have that there exists an infinite class of smooth functions  $w$ (i.e., infinitely differentiable) so that normal velocity $v_n$ (or pressure $p_b$ ) given by,
\begin{eqnarray}
& v_n(\Bx) =&\displaystyle \frac{-i}{\rho c k}\frac{\partial}{\partial\Bn_{\Bx}}\int_{\partial D_{a'}}w(y)\frac{\partial
		\Phi(\Bx,\By)}{\partial \Bn_{\By}}ds_{\By}, \mbox{ for } \Bx\in\partial D_a,\label{2}\\
& p_b(\Bx)=&\displaystyle \int_{\partial D_{a'}}w(y)\frac{\partial
		\Phi(\Bx,\By)}{\partial \Bn_{\By}}ds_{\By}, \mbox{ for } \Bx\in\partial D_a,\label{2''}
\end{eqnarray}
(where $\rho$ denotes the density of the surrounding medium, $D_a'\Subset D_a$ is a smooth (i.e. with $C^2$ boundary) compact region, $\Bn_\By$ denotes the exterior normal to $\partial D_a'$ computed in $\By\in \partial D_a'$ and $\Phi$ is the free space fundamental solution of the Helmholtz equation),  will generate the required acoustic field $u$ satisfying \eqref{1} and \eqref{1-c}. 

\begin{remark}
\label{rem0}
Note that the fact that $D_a'$ is smooth with $D_a'\Subset D_a$ in \eqref{2} or \eqref{2''} implies that the boundary input $v_n$ or $p_b$ is smooth on $\partial D_a$. Moreover, this permits us to assume minimal smoothness for the boundary of the actual physical source $\partial D_a$ (i.e., just enough to have the exterior problem well posed and thus Lipschitz will suffice) which may be very important for some applications.
\end{remark}

\begin{remark}
\label{rem1}
Observe that the normal velocity $v_n$ (or pressures $p_b$) defined at \eqref{2}, (or \eqref{2''}) generate a solution $u$ of \eqref{1}, \eqref{1-c} represented as a double layer potential defined by $$u(\Bx)=\displaystyle \int_{\partial D_{a'}}w(y)\frac{\partial
		\Phi(\Bx,\By)}{\partial \Bn_{\By}}ds_{\By}, \mbox{ for } \Bx\in \RR^3\!\setminus \!{\overline D}_a,$$
		but, it is elementary to see how the results of \cite{20} can be extended to obtain solutions of \eqref{1}, \eqref{1-c} represented as single layer potentials or a linear combination between double layer and single layer potentials. 

\end{remark}

\subsection{Acoustic control in homogeneous oceans of constant depth}
\label{ho}

For the case of constant depth homogeneous ocean environments we model the surrounding homogeneous media as an infinite rectangular wave-guide, with constant depth $h$, i.e,  $z\in[h,0]$ (where $z$ denotes the vertical coordinate in a rectangular coordinate system and $h<0$), and assume a pressure release condition at the water-air interface, i.e., zero pressure at $z=0$, and total pressure reflection at the bottom ocean interface, zero normal pressure at $z=h$ interface \cite{23,24}.

 Let $\RR_h^3=\{\Bx=(\tilde{\Bx},z)\in\RR^3, \tilde{\Bx}\in\RR^2, h\leq z\leq 0\}$ and consider domains $D_1$, $D_2$ and functions $u_1$ and $u_2=0$ as in Section \ref{fs}. Assuming cylindrical coordinates and using the same notations as in \eqref{1}, \eqref{1-c} the problem can be formulated mathematically as follows:

\medskip

{\bf Problem 2. }

Find normal velocity ${v}_n$ (or the pressures $p_b$) on the boundary of the antenna $\partial D_a$ so that $u$, the solution of,
 
\begin{equation}
\label{3}
\vspace{0.15cm}\left\{\vspace{0.15cm}\begin{array}{llll}
		 \GD u+k^2u=0 \mbox{ in }\RR_h^3\!\setminus \!\overline{D}_a\vspace{0.15cm},\\
		 \nabla u\cdot \Bn=v_n, (\mbox{ or } u=p_b), \mbox{ on }\partial D_a\vspace{0.15cm},\\
		 u\!=\!0 \mbox{ on } z\!=\!0, \;\;\displaystyle\frac{\partial u}{\partial z}\!=\!0 \mbox{ on } z\!=\!h,\\
		 \mbox{Radiation condition at infinity uniformly when } |\tilde{\Bx}|\rightarrow\infty,\end{array}\right.
\end{equation}
satisfies 
\begin{equation}\label{3-c}u\approx u_1, \mbox{ in } D_1, \mbox{ and } u\approx 0, \mbox{ in } D_2,\end{equation} where as above in Section \ref{fs} the condition on $D_2$ is not needed in the case when $D_2=\emptyset$.  

We mention that, the radiation condition at infinity in problem \eqref{3} is understood as in \cite{23,24}, i.e., for the solution $u$ represented in normal mode expansion 
$$u(\tilde{\Bx},z)=\displaystyle\sum_{n=0}^{\infty} \phi_n(z)\psi_n(\tilde{\Bx}), \mbox{ for } r>R,$$
(where $r=|\tilde{\Bx}|$) we have that 
\begin{equation}
\phi_n=\sin[k(1-a_n^2)^{\frac{1}{2}}z], \mbox{ with } a_n=\left[ 1-\frac{(2n+1)^2\pi^2}{4k^2h^2}\right]^{\frac{1}{2}} ,
\end{equation}
and 
each of the $\psi_n$ satisfy the following radiation condition when $r\rightarrow\infty$
$$\lim_{r\rightarrow\infty}r^{\frac{1}{2}}\left( \frac{\partial \psi_n}{\partial r}-ika_n \psi_n\right)=0,$$ uniformly for $\theta\in[0,2\pi]$ where $\theta=\tan^{-1}(\frac{x_2}{x_1})$. Next, we will describe how the results in \cite{20} can be extended for this case. Indeed, as above, assuming that $k$ is not a resonant frequency, we have that there exists an infinite class of smooth functions  $w$ so that normal velocity $v_n$ (or pressures $p_b$) given by,

\begin{eqnarray}
&v_n(\Bx)=&\displaystyle\frac{-i}{\rho c k}\frac{\partial}{\partial\Bn_\Bx}\int_{\partial D_{a'}}w(\By)\frac{\partial
		G(x,y)}{\partial \Bn_{y}}ds_{\By}, \mbox{ for } \Bx\in\partial D_a,\label{4}\\
&p_b(\Bx)=&\int_{\partial D_{a'}}w(\By)\frac{\partial
		G(x,y)}{\partial \Bn_{y}}ds_{\By}, \mbox{ for } \Bx\in\partial D_a,\label{4''}		
\end{eqnarray}
(where $\rho$ denotes the density of the surrounding media, $D_a'\Subset D_a$ is a smooth region, $\Bn_\By$ denotes the exterior normal to $\partial D_a'$ computed in $\By\in \partial D_a'$ and $G$ is the Green's function associated to problem \eqref{3}), will generate the required acoustic field $u$ satisfying \eqref{3} and \eqref{3-c}. Indeed, it is observed in \cite{24} that for $\By=(\zeta,\tilde{\By})$, with $\tilde{\By}=(y_1,y_2)$, the Green function $G$ associate to problem \eqref{3} is given by, 
$$G(z,\zeta,|\tilde{\Bx}-\tilde{\By}|)=\Phi(\Bx,\By)+\Phi_1(z,\zeta,|\tilde{\Bx}-\tilde{\By}|),$$
where $\Phi$ is the fundamental free space solution of Helmholtz equation and $\Phi_1$ above is bounded and continuous at $z=\zeta$ and $\tilde{\Bx}=\tilde{\By}$. Based on these onsiderations it is concluded in \cite{24}  that the double layer and the single layer operators associated with the Green's function $G$ have the same compactness properties and satisfy the same jump relations as the classical layer potentials associate to $\Phi$. Thus, by using this together with a few elementary technical adjustments it can be proved that the results presented in \cite{20} will extend to this case, i.e., normal velocities (or pressures) given by \eqref{4} (or \eqref{4''}) will generate acoustic fields described by double layer potentials associated to $G$ and satisfying \eqref{3} and \eqref{3-c}. Moreover, we make the observation that the expressions \eqref{4}, \eqref{4''} can be used in computations since the Green's function $G$ is computed explicitly in \cite{24}. 

In this regard we mention that the statement of Remark \ref{rem0} apply to the case of finite depth homogeneous oceans as well. The following remark is similar in spirit with Remark \ref{rem1} but we presented here for the sake of completeness.  
\begin{remark}
\label{rem2}
The normal velocity $v_n$ (or the pressure $p_b$) given at \eqref{4} (or \eqref{4''}) generates a solution $u$ of \eqref{3}, \eqref{3-c} represented as a double layer potential defined by
$$u(\Bx)=\int_{\partial D_{a'}}w(\By)\frac{\partial
		G(x,y)}{\partial \Bn_{y}}ds_{\By}, \mbox{ for } \Bx\in \RR^3\!\setminus \!{\overline D}_a,$$
but, in a similar manner as above, the results of \cite{20} could be easily extended to obtain solutions of \eqref{3}, \eqref{3-c} represented as single layer potentials or a linear combination between double layer and single layer potentials. 
\end{remark}

\section{Optimization schemes and Numerical simulations}
\label{numerical}

In this section we describe the mathematical ideas behind the optimization scheme used towards the approximation of solutions to \eqref{1}, \eqref{1-c} and respectively \eqref{3}, \eqref{3-c}. 

The $L^2$- optimization and sensitivity analysis for the 2D formulation of the problem \eqref{1}, \eqref{1-c} in the case \eqref{geom-setting} $iii)$ (with  $\{D_1\cup D_2\}\cap D_a=\emptyset$ as above), was performed in \cite{22} where it was numerically observed that good approximation of a stable solution with minimal power budget is achieved in the reactive near field of the source, i.e., when $D_1$ in \eqref{1} is located very close to the source $D_a$.  

Similarly as in the 2D case treated in \cite{22}, the 3D $L^2$- optimization scheme for problem \eqref{1}, \eqref{1-c} is based on Tikhonov regularization with Morozov discrepancy principle.  In this context, as in \cite{20}, \cite{22} regularity results and the well posedness of the interior and exterior acoustic boundary value problem (recall that $k$ was chosen non resonant) imply that in order to achieve approximate smooth controls in $D_1$ and $D_2$ it will be sufficient to have approximate $L^2$ controls on the boundaries of two slightly larger sets, $W_1\Subset\RR^3$ and $W_2\Subset\RR^3$, i.e., with $D_1\Subset W_1$ and $D_2\Subset W_2$. From \eqref{2} and Remark \ref{rem1} it follows that solutions of \eqref{1}, \eqref{1-c} can be approximated by a linear combination of double and single layer potentials, i.e. 

\begin{equation}
\mathcal{D}w_\alpha(\mathbf{x}) = \eta_{1}\int_{\partial D_{a'}} w_\alpha(\mathbf{y}) \frac{\partial \Phi(\mathbf{x}, \mathbf{y})}{\partial \nu_{\mathbf{y}}}\,dS_{\mathbf{y}} + i\eta_{2} \int_{\partial D_{a'}} w_\alpha(\mathbf{y}) \Phi(\mathbf{x},\mathbf{y})\,dS_{\mathbf{y}}, \label{lp}
\end{equation}
where $\eta_{1}, \eta_{2} \in \mathbb{R}$ are fixed parameters and where $w_\alpha$ is the Tikhonov regularization solution, i.e., minimizer of the following discrepancy functional,
\begin{equation}
F(w)=\frac{1}{\|f\|_{L^{2}(\partial D_{c})}^{2}}\| \mathcal{D}w - u_1\|_{L^{2}(\partial W_{1})}^{2} + \mu\|\mathcal{D}w\|_{L^{2}(\partial W_{2})}^{2} + \alpha\|w\|_{L^{2}(\partial D_{a'})}^{2}, \label{eq:tikhonovobjective}
\end{equation}
with  the regularization parameter $\alpha$ chosen according to the Morozov Discrepancy principle (see \cite{22} for a 2D implementation in the case \eqref{geom-setting} $iii)$ and \cite{25,26} for the general theoretical discussion) and the weight $\mu$ above given by
\begin{equation}\mu=\left\{\begin{array}{lll}
\label{mu}
0, & \mbox{ if } D_2=\emptyset,\\
1, & \mbox{ if } D_2 \mbox{ is bounded },\\
\frac{1}{4\pi R^2}, & \mbox{ if } D_2 =\RR^3\setminus B_R(\B0),
\end{array}\right.\end{equation}
with $B_R(\B0)$ denoting the ball centred in the origin with radius $R$ such that $D_a\cup D_1\Subset B_R(\B0)$.  

For the numerical simulations we make use of the spherical harmonic decomposition for $w_\alpha$ (the density of the layer potential operators used to represent the solution \eqref{lp}) and through the method of moments and Tikhonov regularization we approximate a solution of the problem \eqref{1} and \eqref{1-c} in all the geometrical situations described above at \eqref{geom-setting}.  In this regard in all of the numerical simulations below we assumed 30 spherical harmonic orders in the spherical harmonic decomposition of $w_\alpha$.

In all the simulations the fictitious domain $D_{a'}$ appearing in our strategy is the ball centered in the origin and radius $0.01$ and the actual physical source boundary $\partial D_a$ must be located in the annuli $0.01<r<dist(D_a, D_1)$ (where $dist$ here denotes the distance between the two sets) and for all the cases considered at \eqref{geom-setting} $D_1$ is given by
\begin{equation}
\label{D1}
D_1=\{(r,\theta, \phi) , r \in [0.011,0.015], \theta \in [-\frac{\pi}{4}, \frac{\pi}{4}] , \phi \in [\frac{3\pi}{4}, \frac{5\pi}{4}]  \}.
\end{equation} 
In the remainder of the paper will present numerical simulations of our strategy and, to simplify the exposition, we will focus only on problem \eqref{1} and \eqref{1-c}. Thus, the next three sections will present our numerical simulations for the Tikhonov regularization solution corresponding to problem \eqref{1} and \eqref{1-c} as follows: Section \ref{O} discuses the case \eqref{geom-setting} $i)$, Section \ref{I} discusses the case \eqref{geom-setting} $ii)$ while Section \ref{II} discusses the case \eqref{geom-setting} $iii)$.

\subsection{Synthesis of a prescribed pattern in a subregion of the source near-field }

\label{O}

In this section we present the Tikhonov regularization solution for the problem \eqref{1}, \eqref{1-c} introduced in Section \ref{fs} describing the applications to the synthesis of acoustic sources approximating a given field patterns in a near field bounded region $D_1$, see Figure \ref{Geom-synth}.
\begin{figure}[!htbp]\centering
\vspace{-3.5cm}
\includegraphics[width=0.6\textwidth, height=0.6\textwidth]{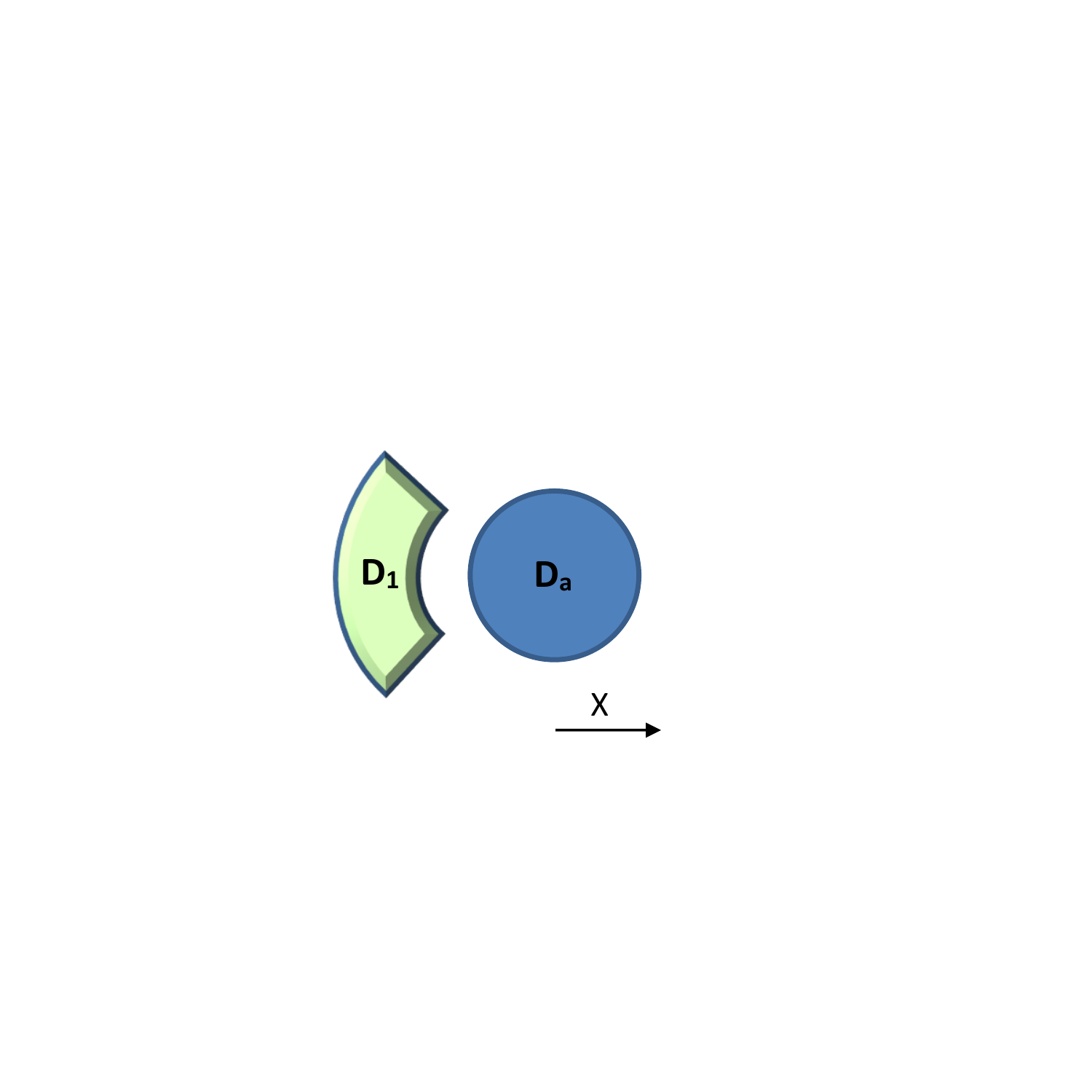}
\vspace{-2.5cm}
	\caption{Planar sketch of the control geometry.}
	\label{Geom-synth}
\end{figure}

\begin{figure}[!htbp] \centering
\vspace{0cm}
	\begin{subfigure}{\figsizeC}
		\includegraphics[width=1.5\textwidth]{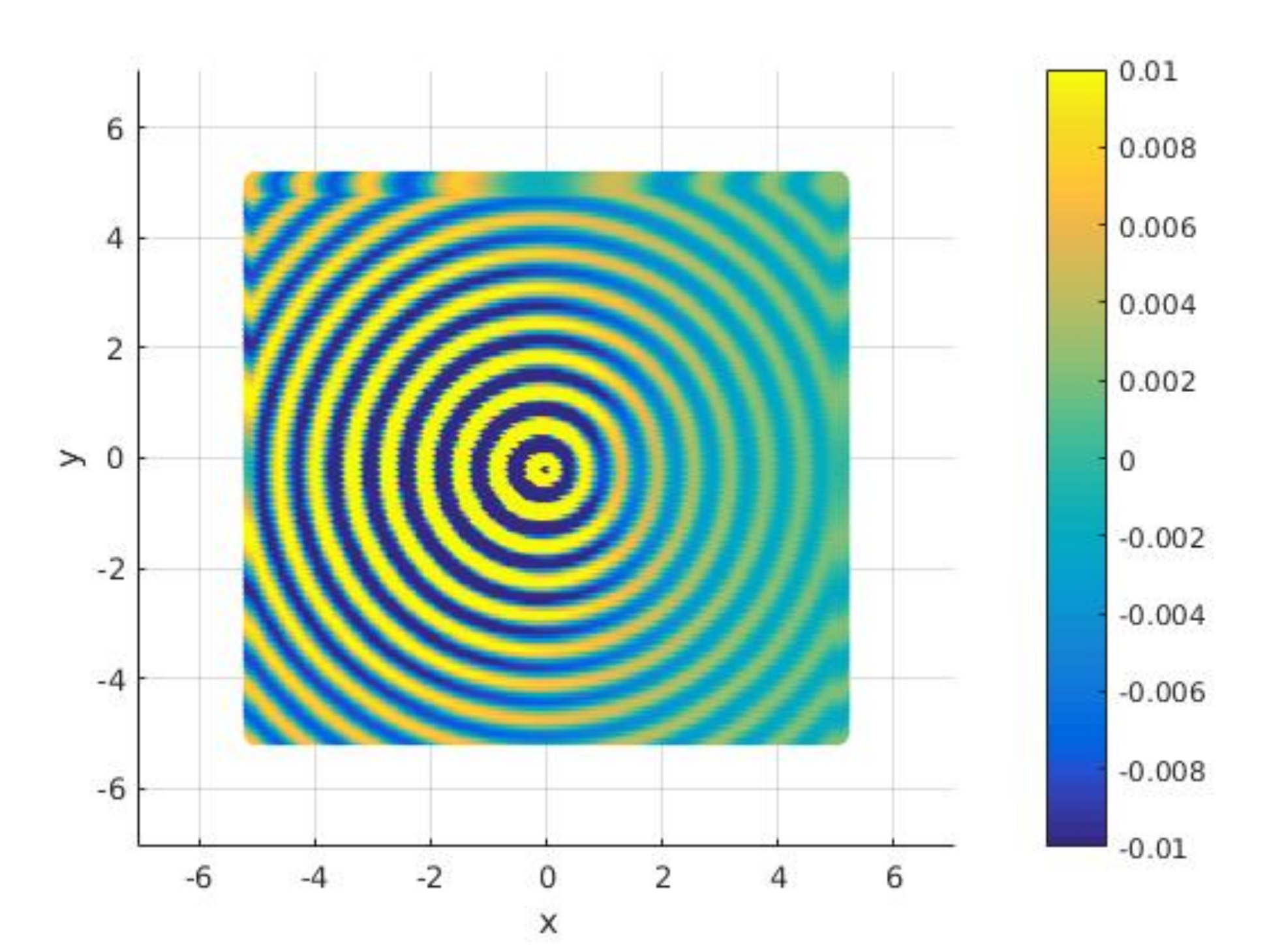} 
	\end{subfigure}
	\caption{Cross-section $z=0$ plot of the generated field} 
	\label{fig:k10h30synth_out}
\end{figure}

Thus, we will show the performance of the Tikhonov solution described in \eqref{lp}, \eqref{eq:tikhonovobjective}, \eqref{mu} in the case \eqref{geom-setting} $i)$. Without loss of generality we consider the case when the source to be synthesized approximates in region $D_1$, described at \eqref{D1}, an outgoing plane wave propagating along the negative $x$-axis, $u_1 = e^{-i x k} $ with wave number $k=10$.

First, in Figure \ref{fig:k10h30synth_out} we present a cross-sectional view of the generated field along $z=0$ in a region characterized by $(x,y)\in[-5,5]^2$. This plot indicates the synthesised source causality (i.e., the fact that the source field is outgoing). This fact can also be observed in the time domain simulation presented in %
\label{HL_s_left_t4}
\href{https://drive.google.com/open?id=0B7nf-pdU3Z4Dbm9IcU9NQkhSVjA}{animation 1}
where the propagating time-harmonic field generated by the synthesized source is shown. 

In Figure \ref{fig:synth_ue} we present the quality of our control results in the region of interest $D_1$ as required in \eqref{1-c}. The left and center plots in the figure describe respectively the field generated by the source, and the plane wave to be approximated $u_1 = e^{-i x k} $ in region $D_1$. The good accuracy of our approximation $O(10^{-3})$ can be observed in the right picture of Figure \ref{fig:synth_ue} where the relative pointwise error between $u$ (the solution of \eqref{1}) and $u_1 = e^{-i x k} $ (the field to be approximated) is presented. 
\begin{figure}[!htbp] \centering
\vspace{0cm}
	\begin{subfigure}{\figsizeC}
		\includegraphics[width=\textwidth]{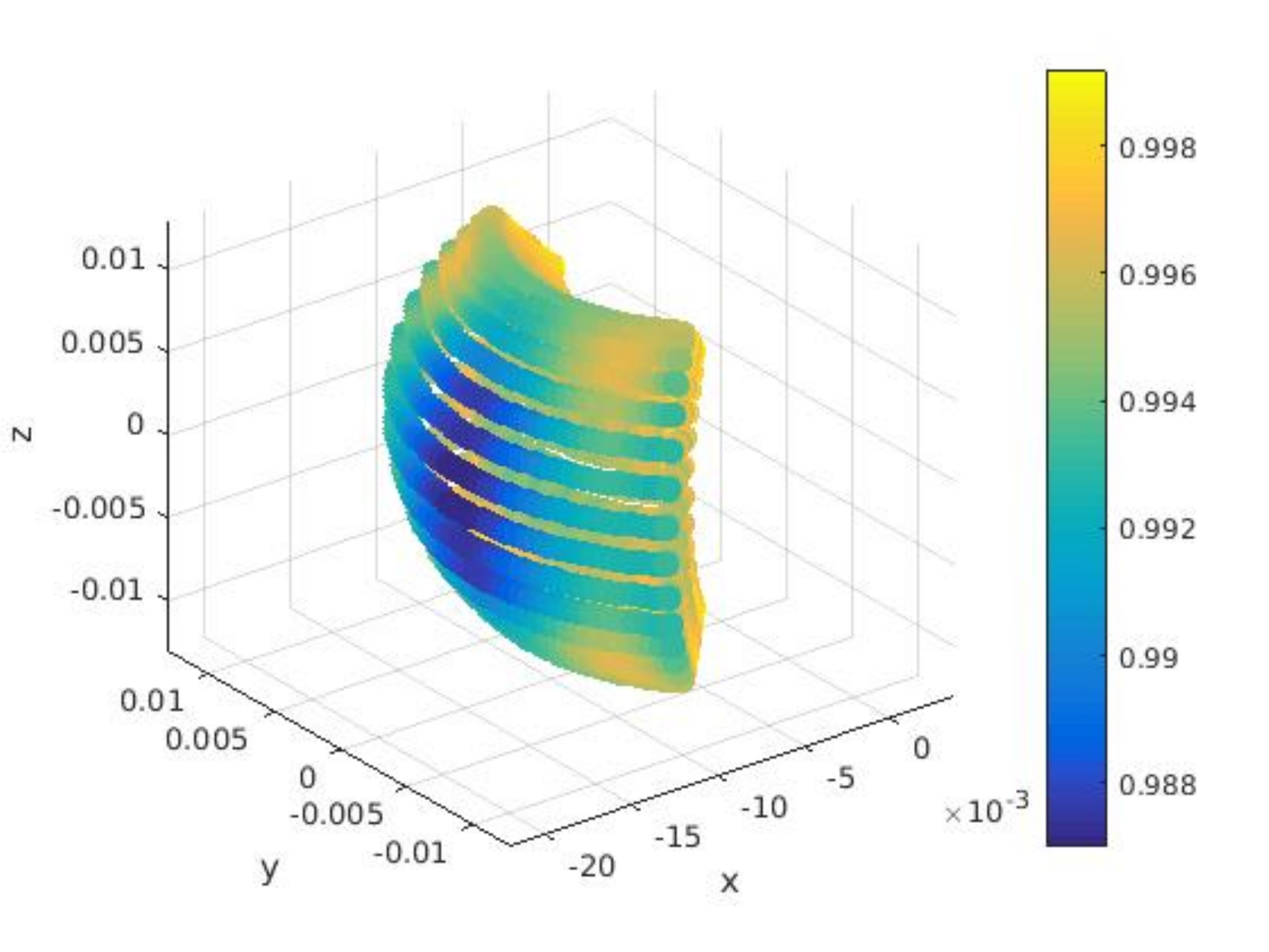}
		\caption{Generated field}
		\label{fig:double_u}
	\end{subfigure}
	\begin{subfigure}{\figsizeC}
		\includegraphics[width=\textwidth]{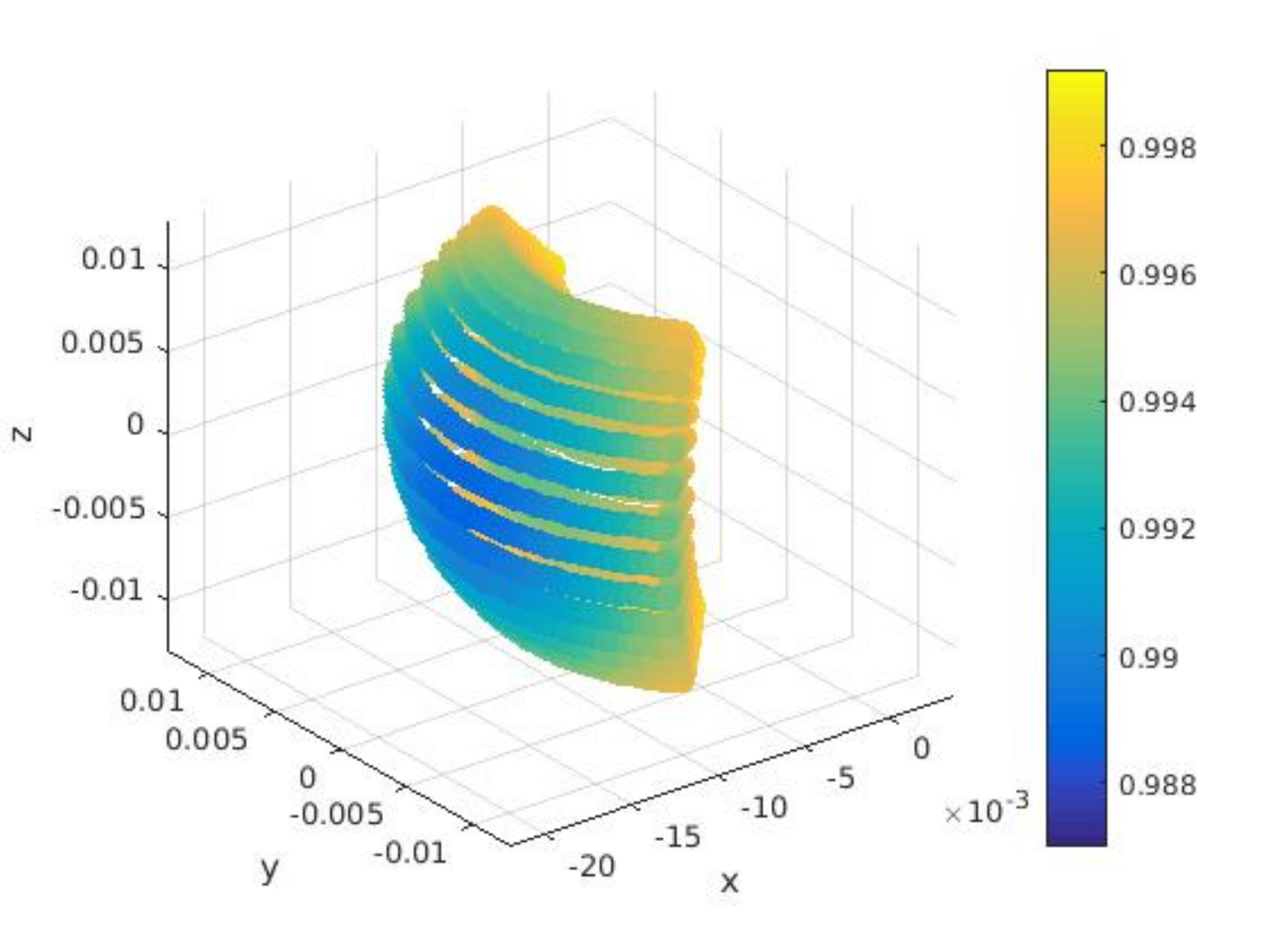}
		\caption{Field to match}
		\label{fig:double_f}
	\end{subfigure}
	\begin{subfigure}{\figsizeC}
		\includegraphics[width=\textwidth]{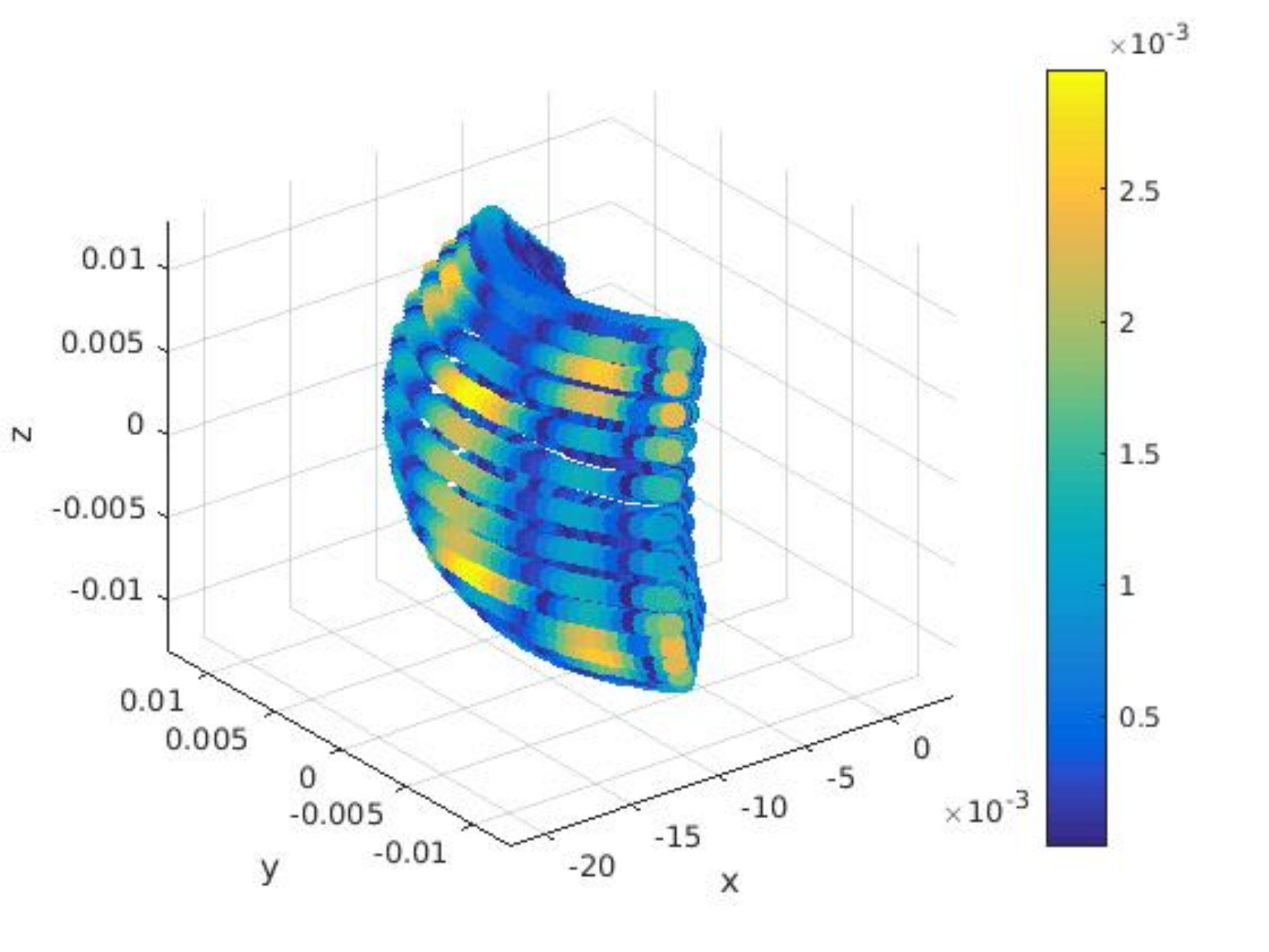}
		\caption{Pointwise relative error}
		\label{fig:double_e}
	\end{subfigure}	
	\caption{Control accuracy in region $D_1$}
	\label{fig:synth_ue}
\end{figure}
 \begin{figure}[!htbp] \centering
 \vspace{0cm}
 	\begin{subfigure}{\figsizeD}
 		\includegraphics[width=\textwidth]{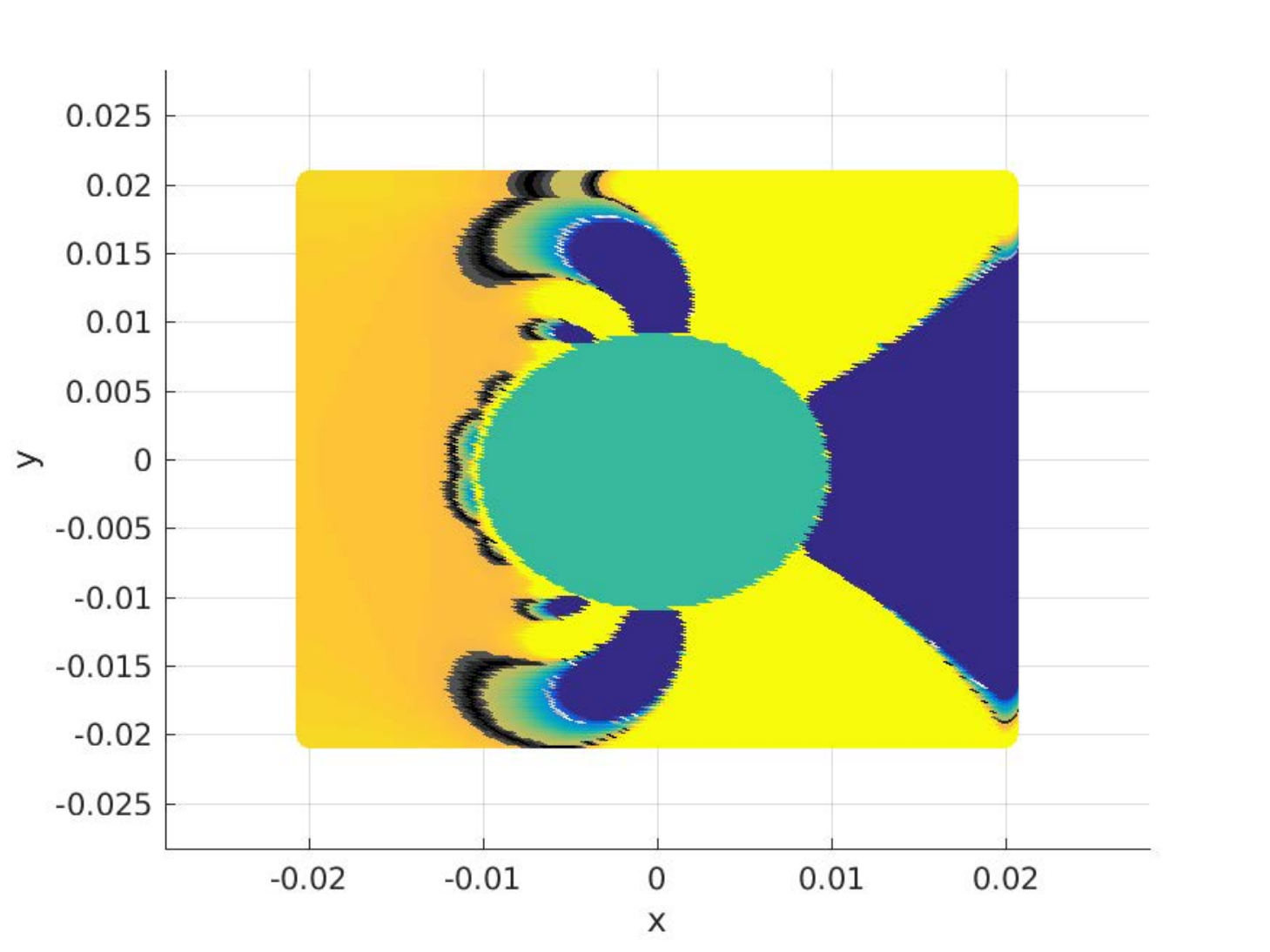}
 		\label{fig:k10h30synthSS15}
 		\caption{$\frac{15}{50}\pi$}
 	\end{subfigure}
 		\begin{subfigure}{\figsizeD}
 		\includegraphics[width=\textwidth]{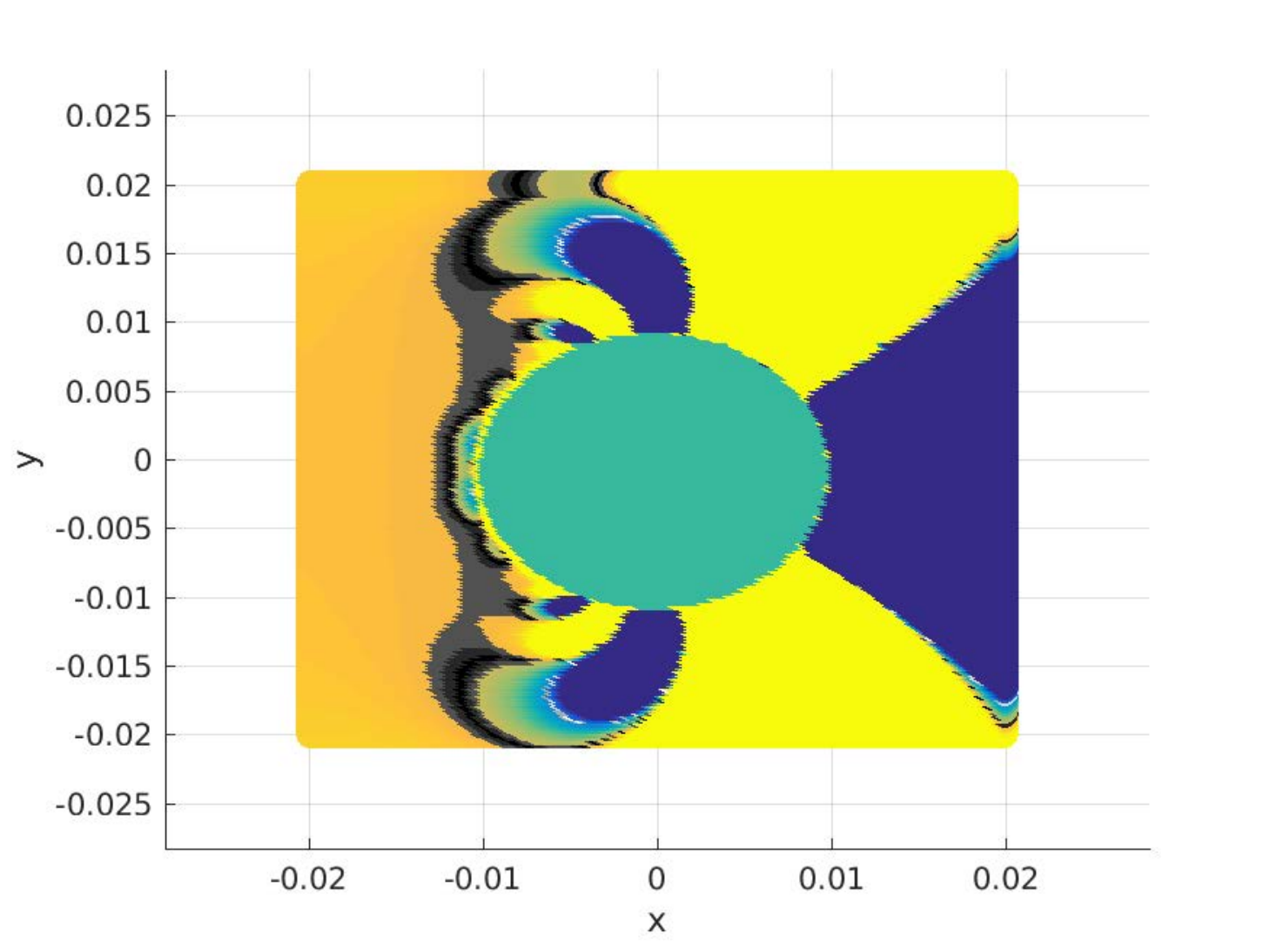}
 		\label{fig:k10h30synthSS16}
 		\caption{$\frac{16}{50}\pi$}
 	\end{subfigure}
 	\begin{subfigure}{\figsizeD}
 	\includegraphics[width=\textwidth]{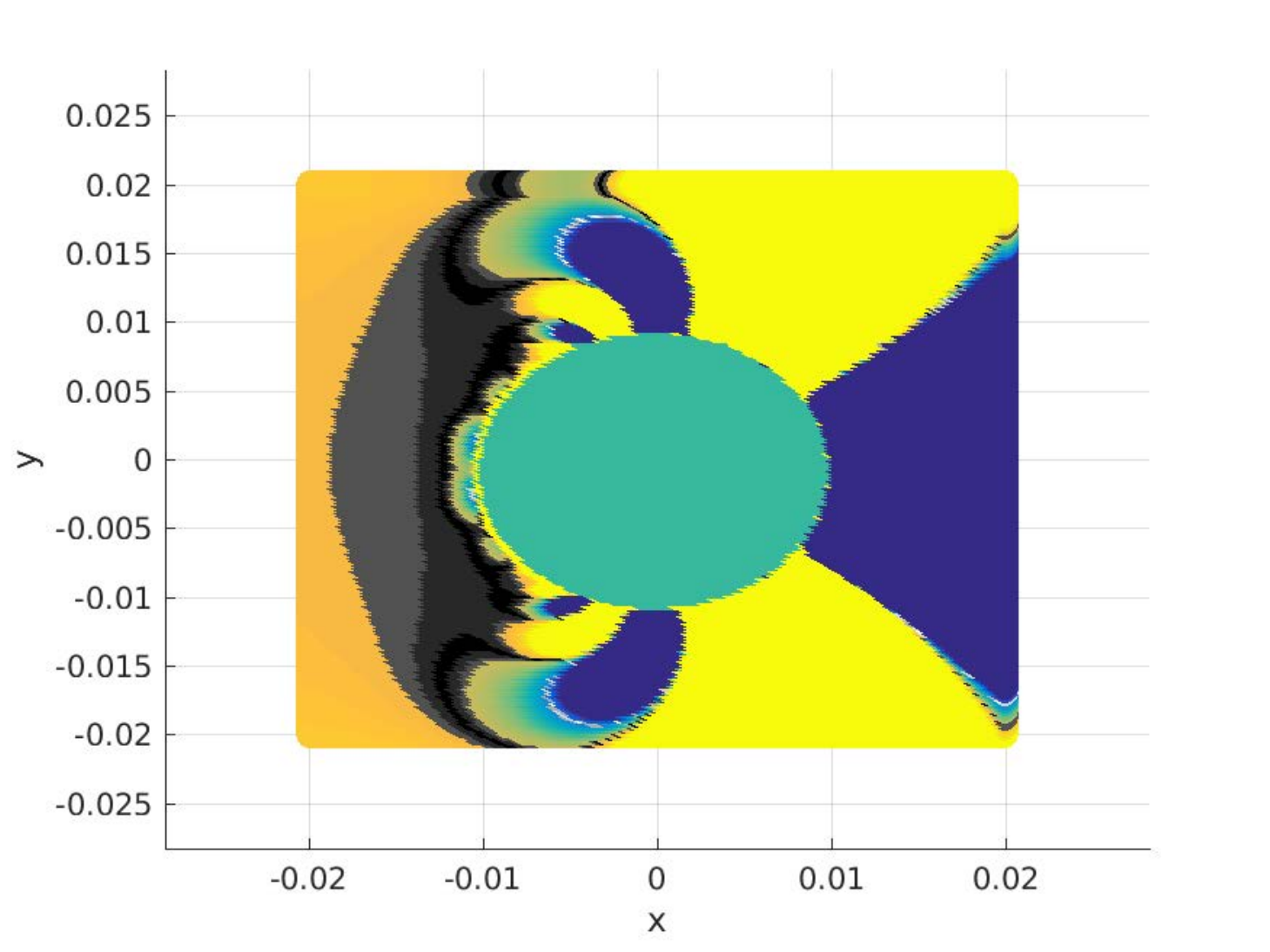}
 	\label{fig:k10h30synthSS17}
 	\caption{$\frac{17}{50}\pi$}
 	\end{subfigure}
 	\begin{subfigure}{\figsizeD}
 	\includegraphics[width=\textwidth]{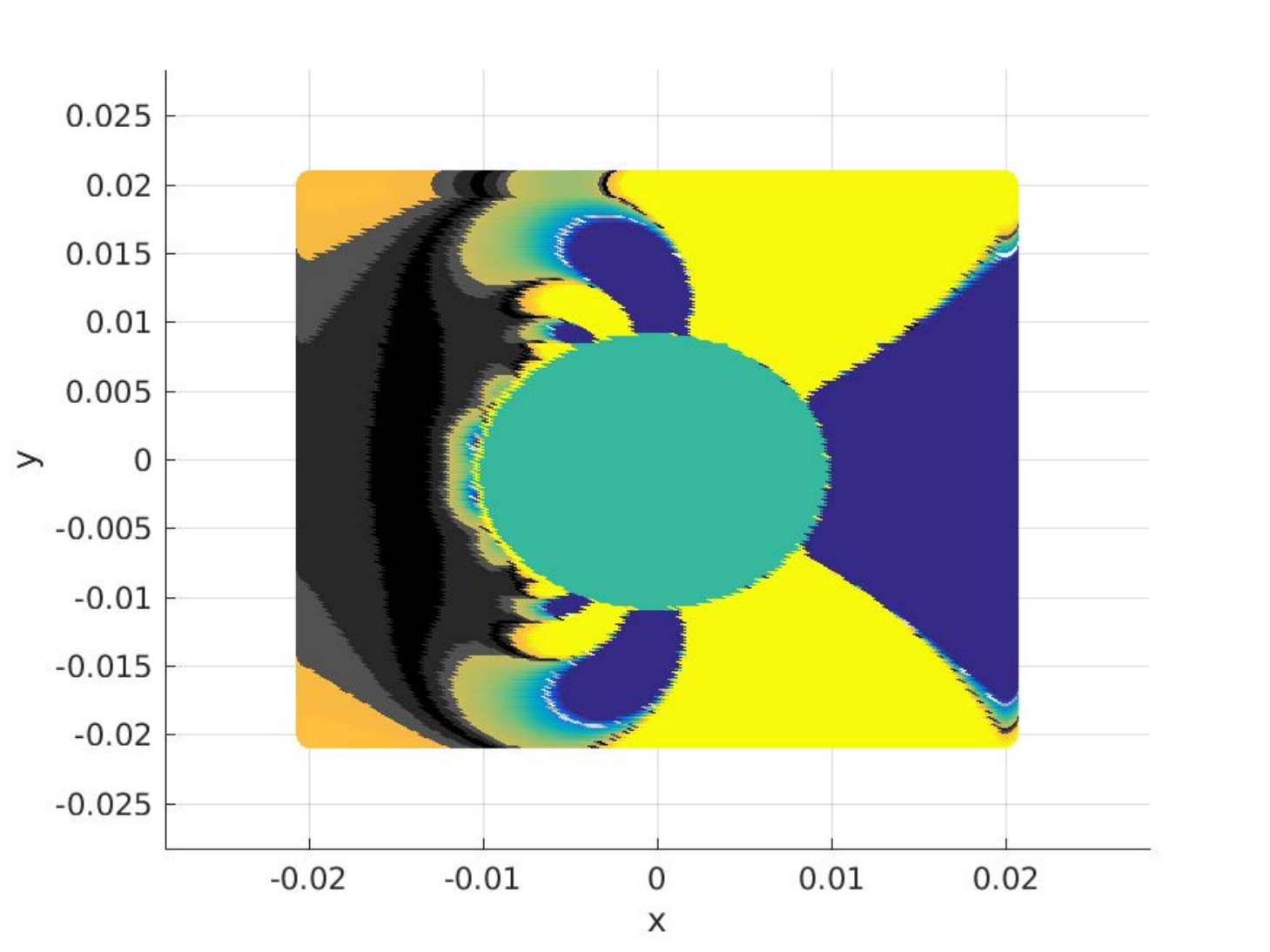}
 	\label{fig:k10h30synthSS18}
 	\caption{$\frac{18}{50}\pi$}
 	\end{subfigure}
 	\begin{subfigure}{\figsizeD}
 		\includegraphics[width=\textwidth]{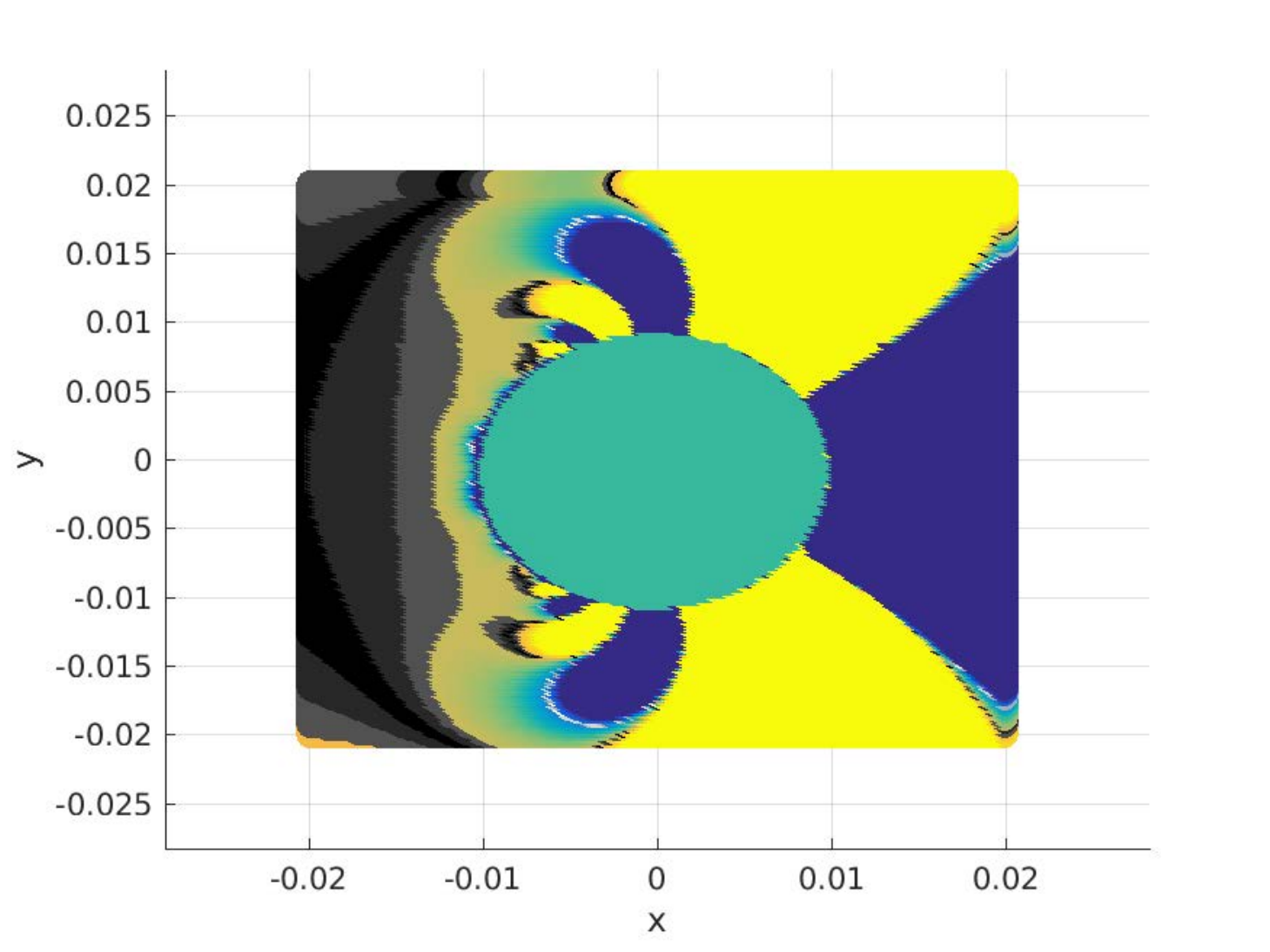}
 		\label{fig:k10h30synthSS19}
 		\caption{$\frac{19}{50}\pi$}
 	\end{subfigure}
 	\begin{subfigure}{\figsizeD}
 		\includegraphics[width=\textwidth]{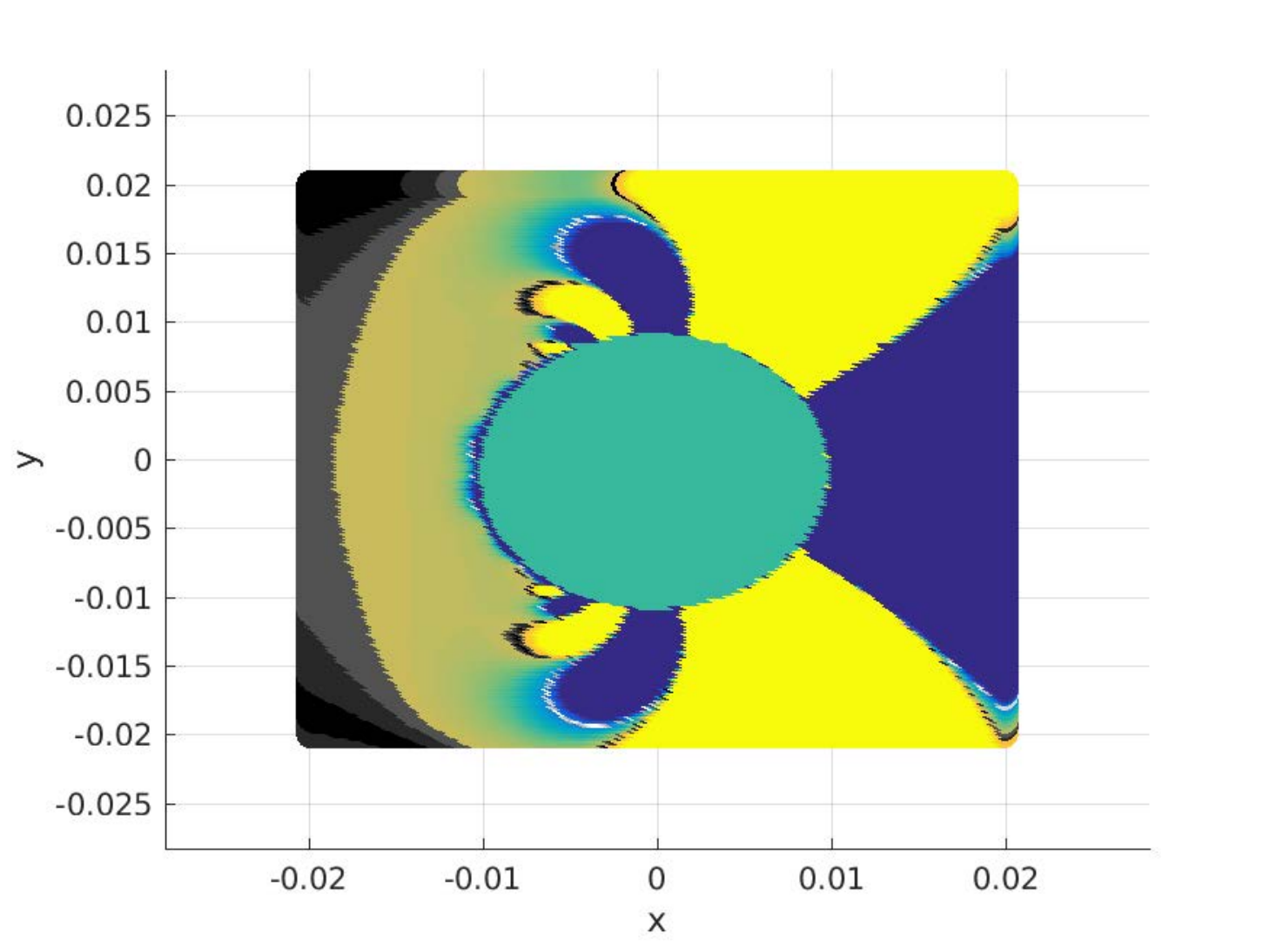}
 		\label{fig:k10h30synthSS20}
 		\caption{$\frac{20}{50}\pi$}
 	\end{subfigure}
 \vspace{0cm}
 	\caption{Cross-sectional ($z=0$) time snapshots of the propagating generated acoustic field for different values of $kct$.}
 	\label{fig:k10h30synth_ss}
 \end{figure}
 Figure \ref{fig:k10h30synth_ss} shows six a cross-sectional views of the generated field along $z=0$ in a near-field region characterized by $(x,y)\in[-0.02,0.02]^2$.  More explicitly, in order left to right from top left to bottom right plot, we present six cross-sectional ($z=0$) time-snapshots ( for $kct=\{\frac{15}{50}\pi,\frac{16}{50}\pi,\frac{17}{50}\pi,\frac{18}{50}\pi,\frac{19}{50}\pi,\frac{20}{50}\pi\}$) of the time-harmonic field generated by the synthesized  source in it's near field region, including the region of interest $D_1$ (where $c$ here was used as the speed of sound in air). The color scheme in the plots is (truncated to 1 light yellow and -1 dark blue) with the antenna region (colored cyan) not included in the numerical simulations and with the black stripe representing field amplitudes of $\approx 0.6$. Following the plots in order from top left to bottom right plot it can be observed how the source works to approximate an outgoing plane wave $u_1 = e^{-i x k}e^{-ikct} $  in region $D_1$ (e.g., corresponding rectilinear black strip outgoing propagating through the control region). Indeed, the plots of Figure \ref{fig:k10h30synth_ss} show the propagation of the generated field by focusing on the portion of the field with amplitude $\approx 0.6$ marked as a dark stripe. It can be observed how this portion of the field enters region $D_1$ at time $kct= \frac{17}{50}\pi$ in a nearly rectilinear shape and continues keeping the same rectilinear form (indicating plane wave character of the approximated field in the control region) outgoing throughout a neighbourhood of region $D_1$. 
 \begin{figure}[!htbp] \centering
   \vspace{0cm}
   	\begin{subfigure}{\figsizeC}
   		\includegraphics[width=\textwidth]{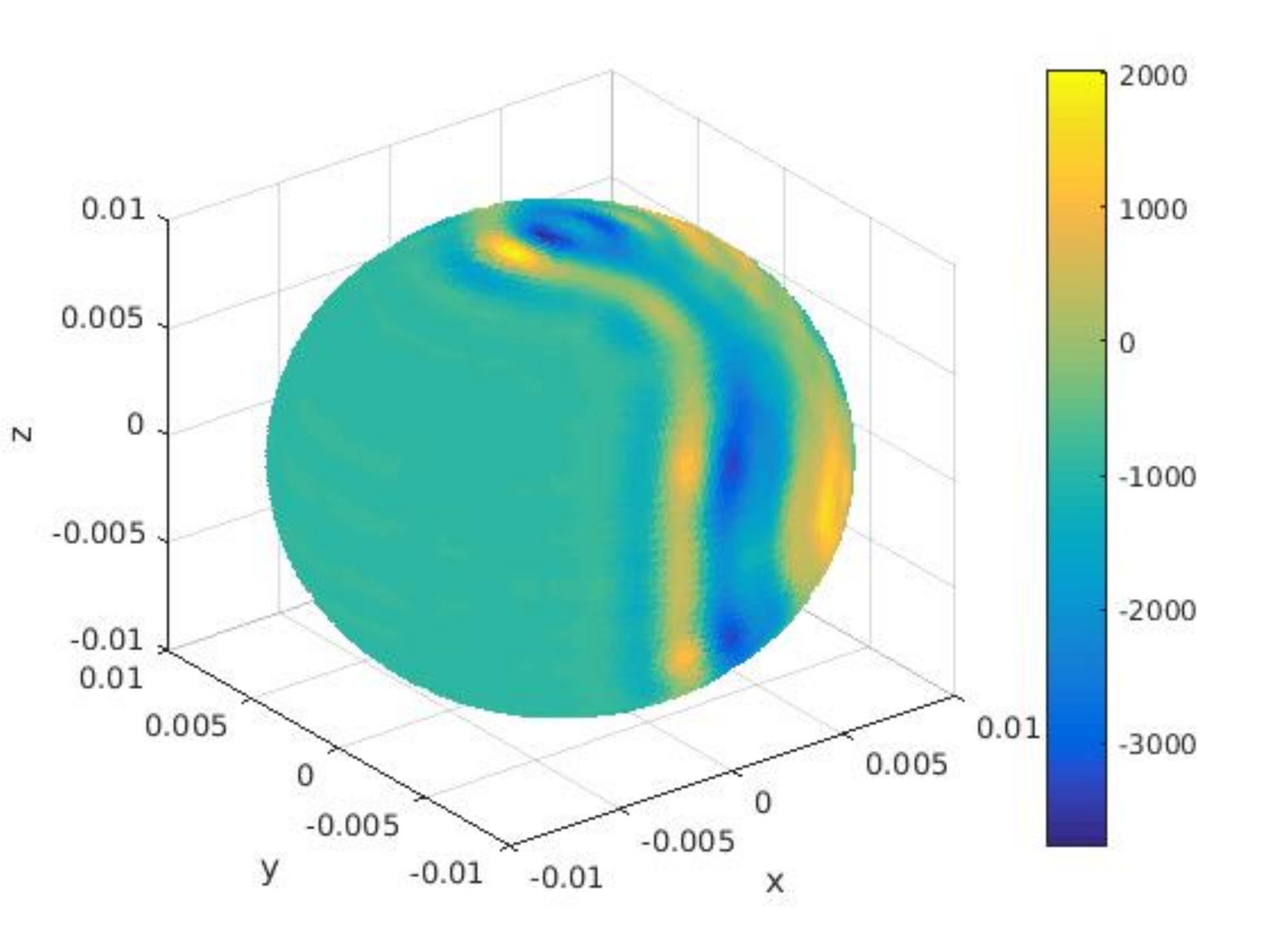}
   		\label{fig:synth_ad_s1}
   		\caption{side}
   	\end{subfigure}
   	\begin{subfigure}{\figsizeC}
   		\includegraphics[width=\textwidth]{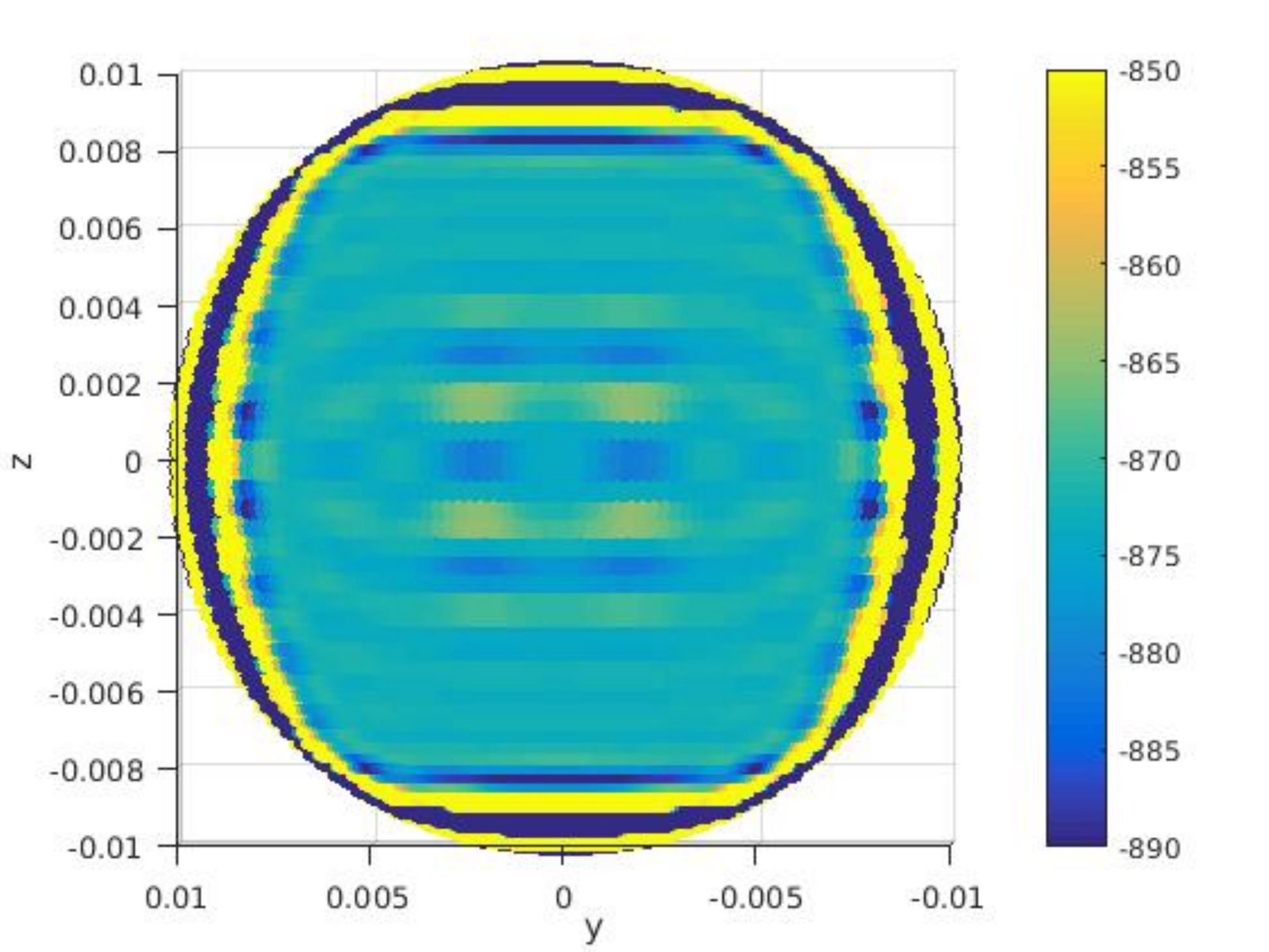}
   		\label{fig:synth_ad_b1}
   		\caption{front}
   	\end{subfigure}
   	\begin{subfigure}{\figsizeC}
   		\includegraphics[width=\textwidth]{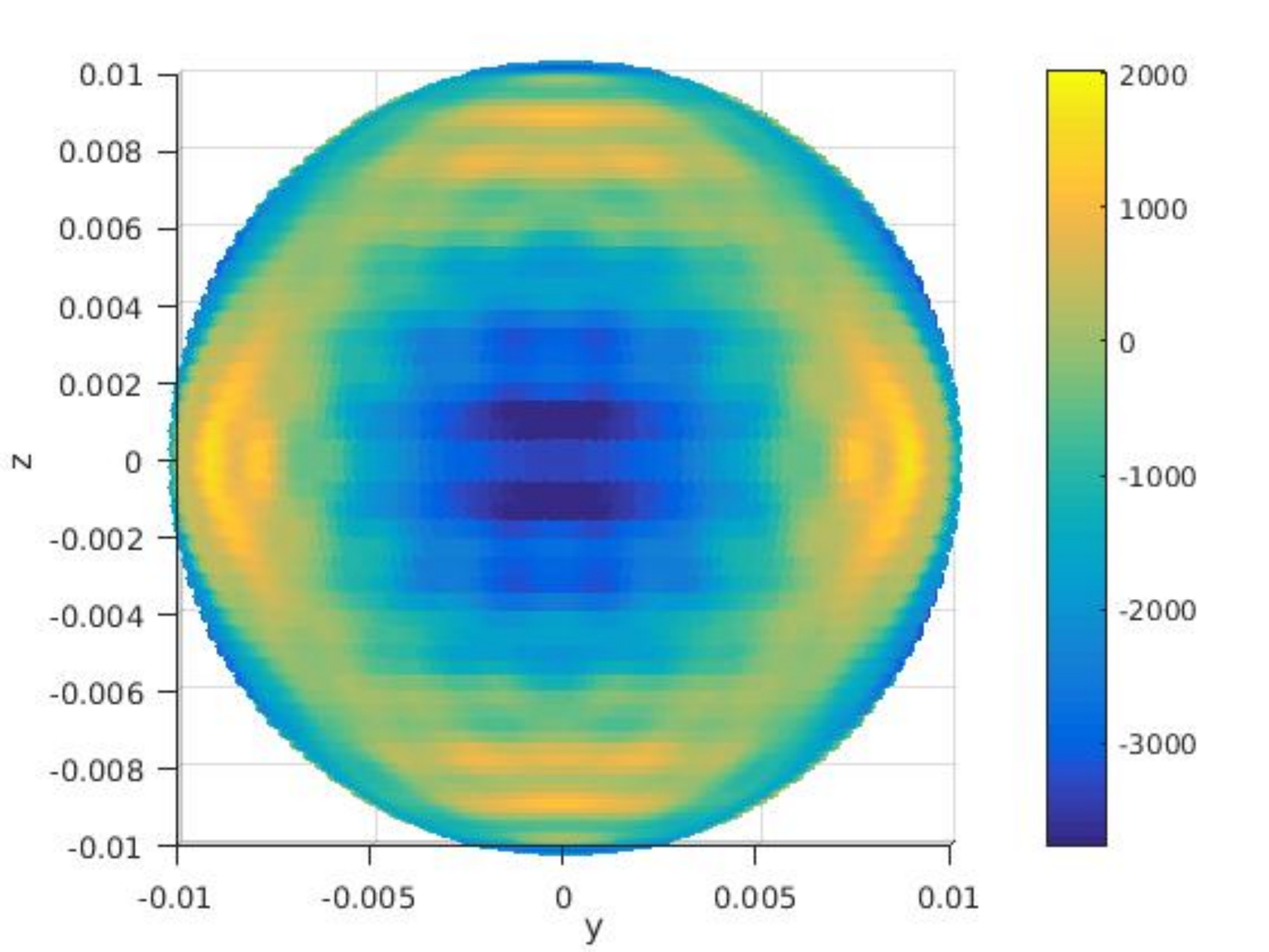}
   		\label{fig:synth_ad_a2}
   		\caption{back}
   	\end{subfigure}
   	\caption{Density $w_\alpha$ with various colour maps}
   	\label{fig:synthdensity}
   \end{figure}
  The time domain animation 
  \label{HL_s_left_t2}
  \href{https://drive.google.com/open?id=0B7nf-pdU3Z4DUjJlUHF4NmRIc2M}{animation 2}
  presents the cross-sectional view along $z=0$ of the time-harmonic evolution of the field generated by the synthesized source and respectively the propagating plane wave $u_1= e^{-i x k}e^{-ikct}$ in a near field region given by $(x,y)\in[-0.1,0.1]^2$. The multimedia file shows two animations: the top one describing the time propagation of the generated field and the bottom one describing the time propagation of the plane wave $u_1= e^{-i x k}e^{-ikct}$. The color scheme in the movies is (truncated to 1 light yellow and -1 dark blue) with the antenna region removed from the simulations (colored cyan) and with the black stripe representing amplitude values around $\approx 0.6$, and the white stripe representing amplitude values around $\approx -0.6$ respectively. Observing that there will be two black stripes and respectively two white stripes per period for the approximated plane wave one can see in the animation the good accuracy of the approximation of the outgoing plane wave $u_1= e^{-i x k}e^{-ikct}$ in region $D_1$.
  
 Figure \ref{fig:synthdensity} describes the density $w_\alpha$ (see \eqref{lp}) on the boundary of the fictitious domain $D_{a'}$ as an indication of the possible complexity of the required source inputs $v_n$ or $p_b$ described at \eqref{2} or \eqref{2''}. In the left plot of the figure we present the density values on the surface of $D_{a'}$ viewed in a side 3D perspective and for better visualization we show two more plots in the figure: the center plot shows the density values on the part of the surface facing region $D_1$ while the right plot of the figure presents the density values on the oposite part of the surface.

\subsection{Synthesis of different prescribed patterns in disjoint subregions of the source near-field}
\label{I}

In this section we present the Tikhonov regularization solution for the problem \eqref{1}, \eqref{1-c} described in Section \ref{fs} describing the applications to the synthesis of acoustic sources approximating two different field patterns in two prescribed disjoint  near field regions. Thus, we show next the performance of the Tikhonov solution described in \eqref{lp}, \eqref{eq:tikhonovobjective}, \eqref{mu} in the case \eqref{geom-setting} $ii)$.

As in Section \ref{O}, we consider the case when the synthesized source approximates in region $D_1$, described at \eqref{D1}, an outgoing plane wave propagating along the negative $x_1$-axis, $u_1 = e^{-i x k} $ with wave number $k=10$ while having a null in region $D_2 = \{(r,\theta, \phi) , r \in [0.011,0.015], \theta \in [-\frac{\pi}{4}, \frac{\pi}{4}] , \phi \in [-\frac{\pi}{4}, \frac{\pi}{4}]  \} + (0.018,0,0)  $, (i.e. region $D_2$ is the same as region $D_1$ but shifted to the right along $x_1$ axis 0.018 units). A sketch of the geometries are presented with Figure \ref{fig:double_g}.
\begin{figure}[!htbp]
\centering
\vspace{-3.5cm}
		\includegraphics[width=0.6\textwidth, height=0.6\textwidth]{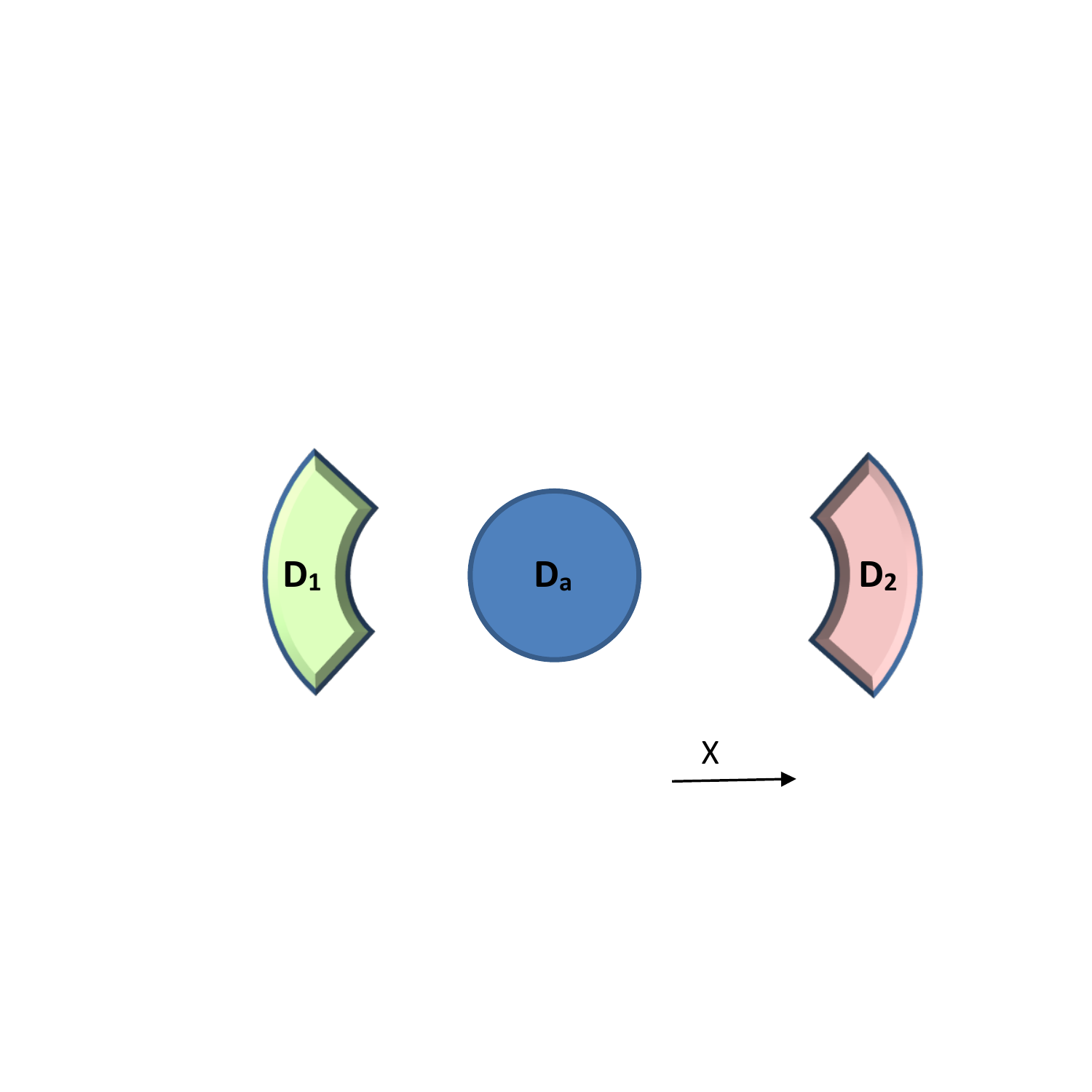}
		\vspace{-2.3cm}
	\caption{Planar sketch of the control geometry.}
	\label{fig:double_g}
\end{figure}
In Figure \ref{fig:k10h30double_out} we present a cross-sectional view of the generated field along $z=0$ in a region characterized by $(x,y)\in[-5,5]^2$. As above, this plot supports the claim about the source causality (i.e., the fact that the source field is outgoing). This fact can be better observed in the time domain simulation presented in 
\label{HL_d_left_t4}
\href{https://drive.google.com/open?id=0B7nf-pdU3Z4DdHpuYWVJU0NjbkU}{animation 3} 
where the propagating time-harmonic field generated by the synthesized source is shown. 
\begin{figure}[!htbp] \centering
\vspace{0cm}
	\begin{subfigure}{\figsizeC}
		\includegraphics[width=1.5\textwidth]{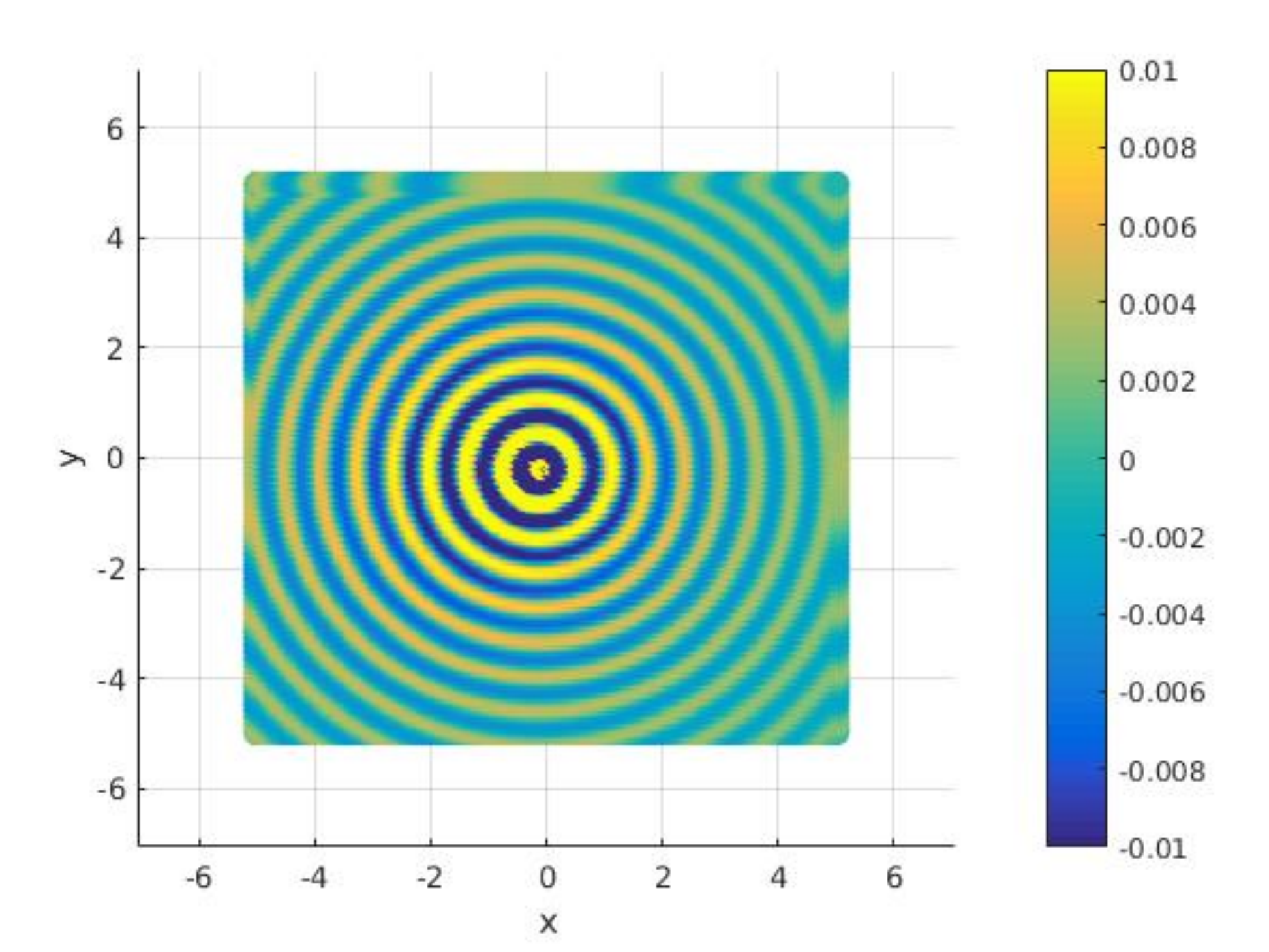} 
	\end{subfigure}
	\caption{Cross-section $z=0$ plot of the generated field} 
	\label{fig:k10h30double_out}
\end{figure}
\begin{figure}[!htbp] \centering
\vspace{0cm}
	\begin{subfigure}{\figsizeC}
		\includegraphics[width=\textwidth]{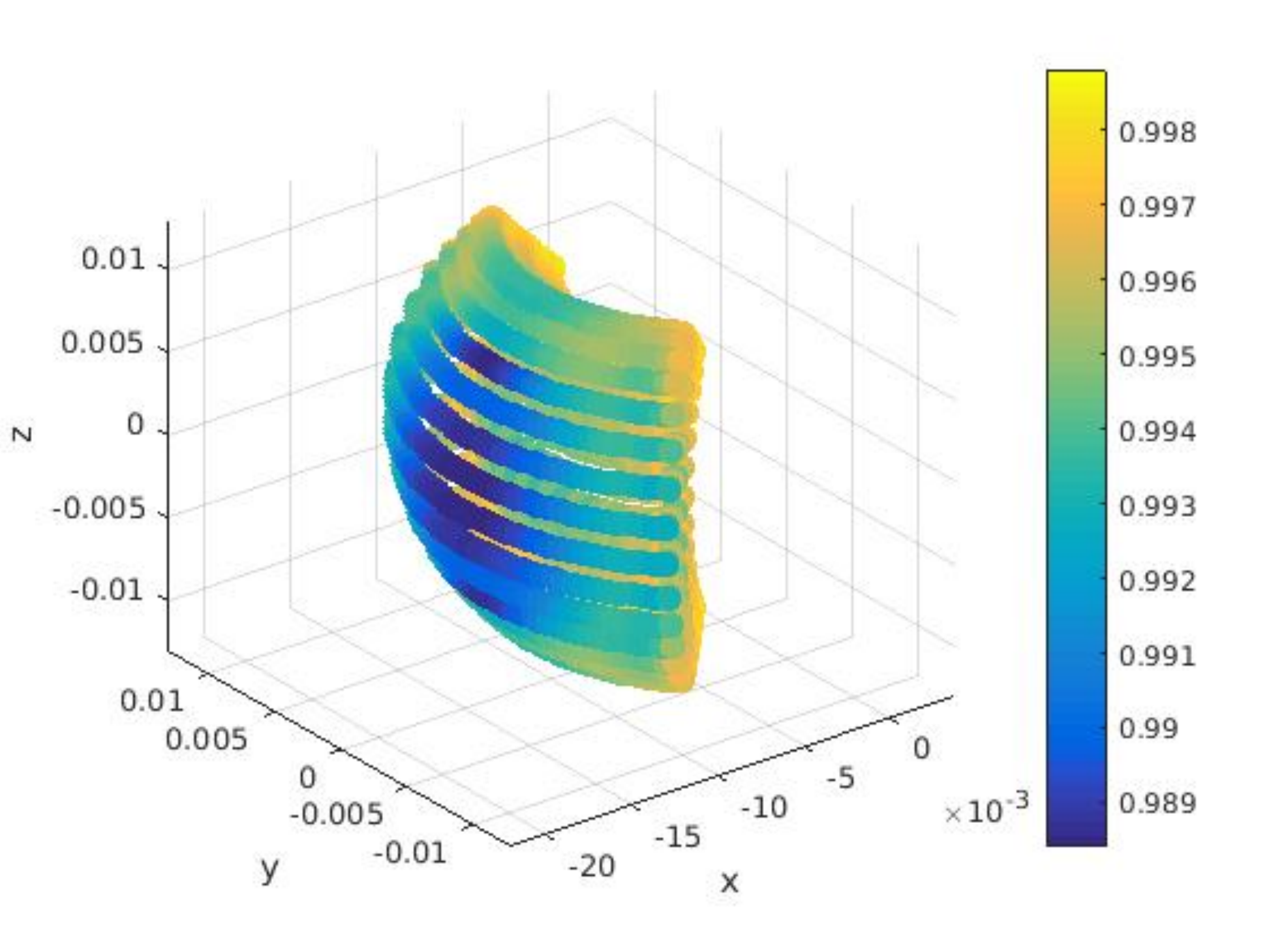}
		\caption{Generated field in $D_1$}
		\label{fig:double_u}
	\end{subfigure}
	\begin{subfigure}{\figsizeC}
		\includegraphics[width=\textwidth]{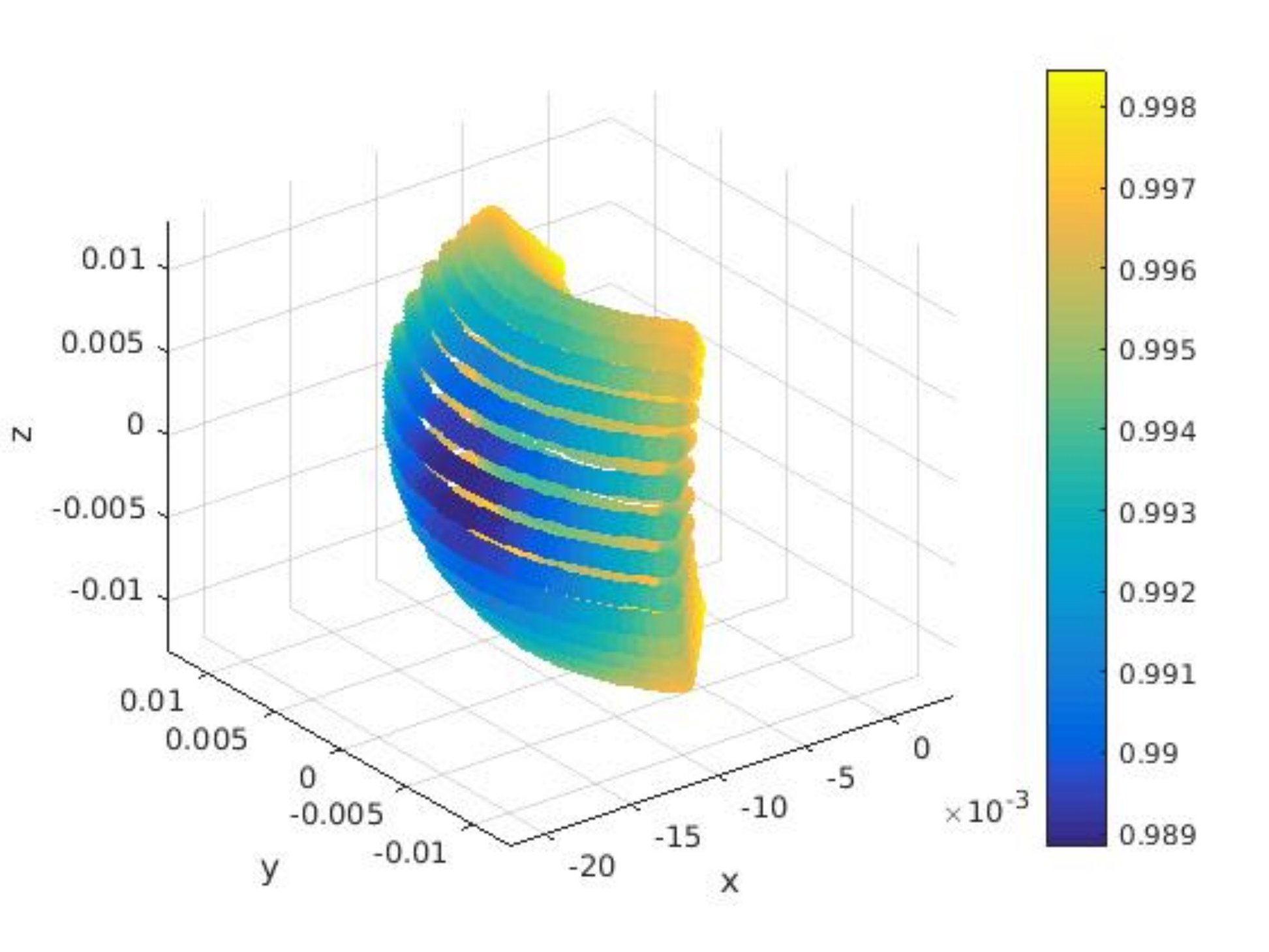}
		\caption{Field to match in $D_1$}
		\label{fig:double_f}
	\end{subfigure}
	\begin{subfigure}{\figsizeC}
		\includegraphics[width=\textwidth]{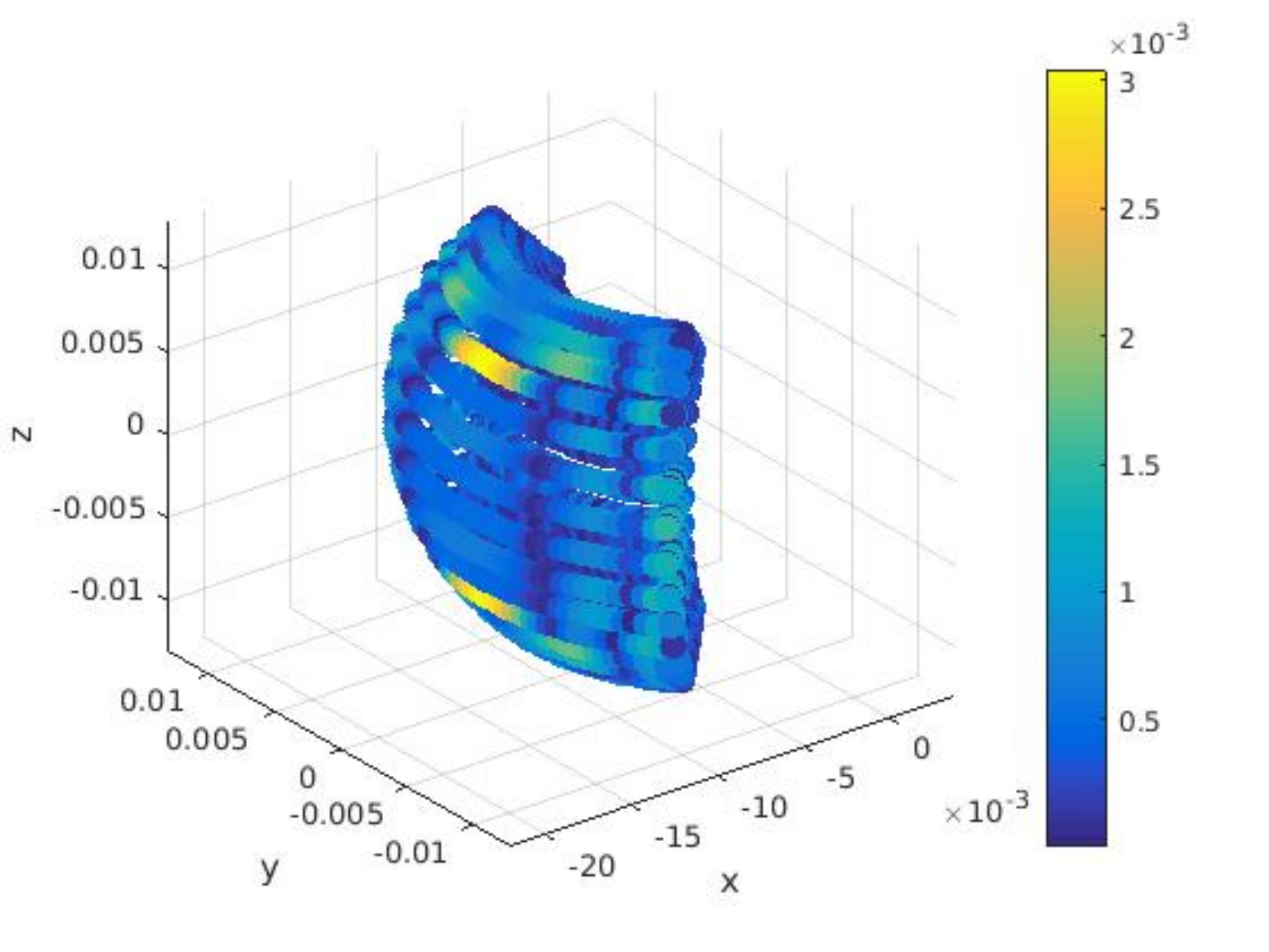}
		\caption{Pointwise relative error}
		\label{fig:double_e}
	\end{subfigure}
	\begin{subfigure}{\figsizeC}
		\includegraphics[width=\textwidth]{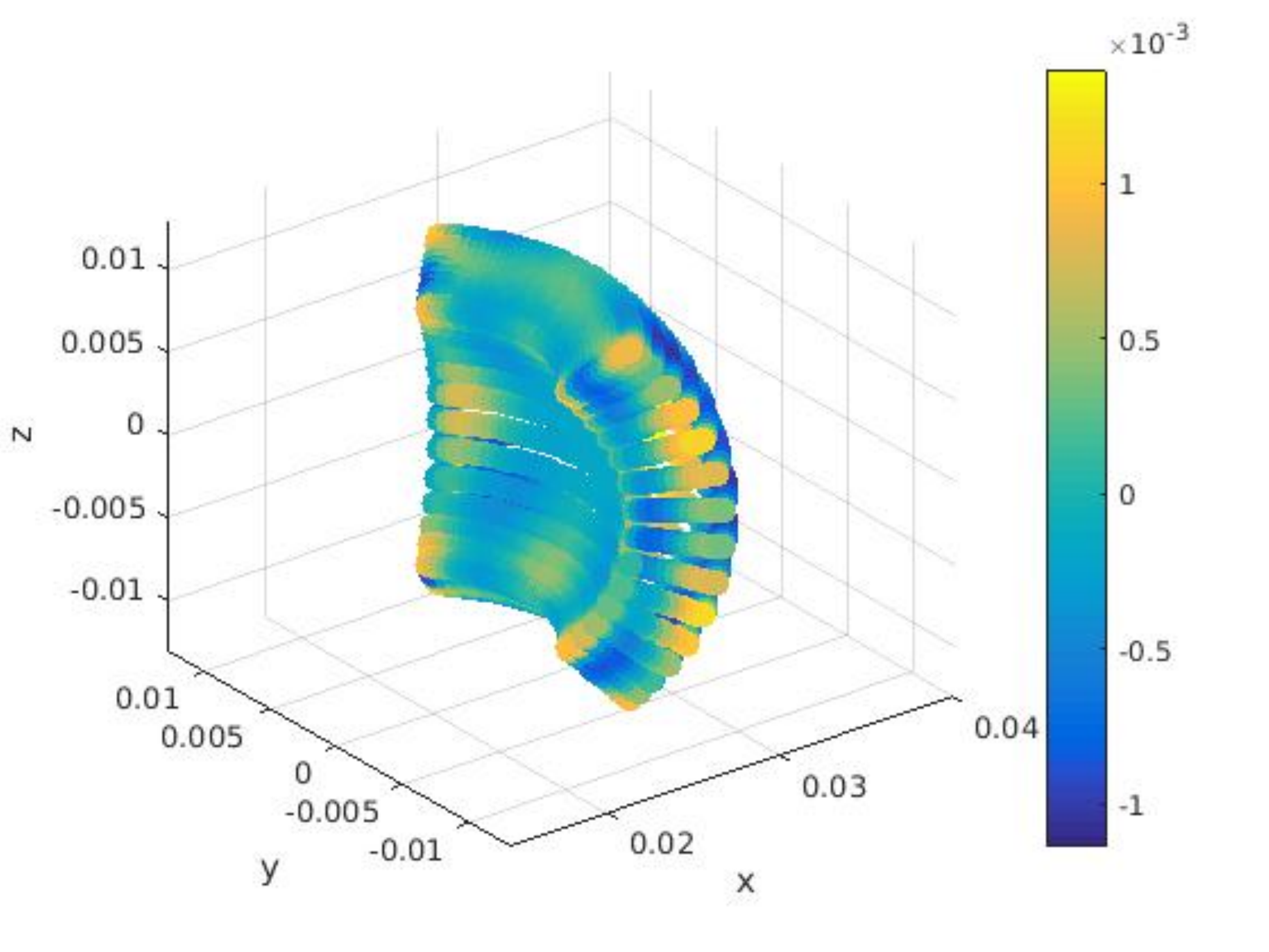}
		\caption{Generated field in $D_2$}
		\label{fig:double_u2}
	\end{subfigure}
	\caption{Demonstration of field matching in $D_1$ and $D_2$}
	\label{fig:double_ue}
\end{figure} 
 \begin{figure}[!htbp] \centering
 \vspace{0cm}
 	\begin{subfigure}{\figsizeD}
 		\includegraphics[width=\textwidth]{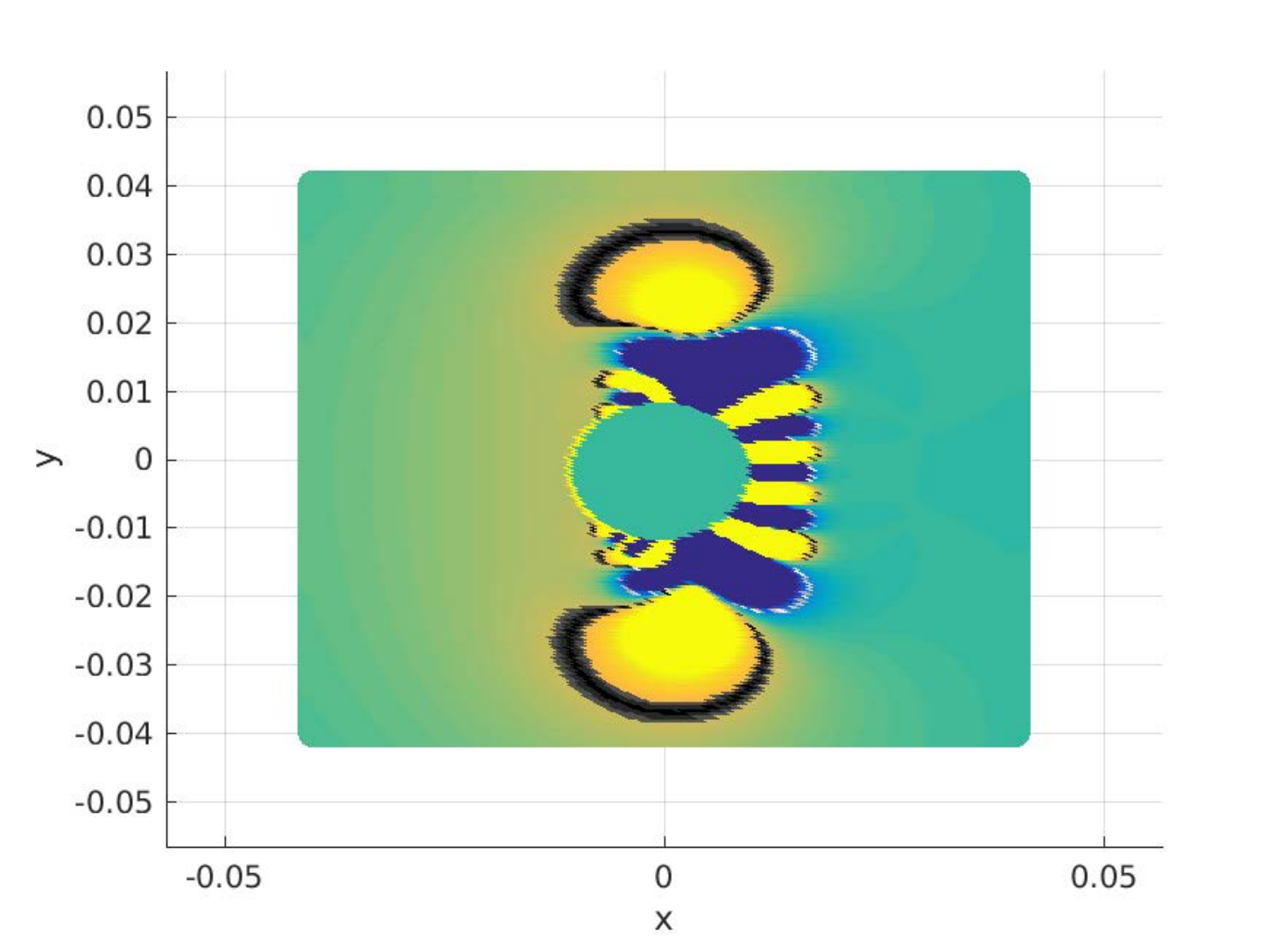}
 		\label{fig:k10h30double_t2_ss83}
 		\caption{$\frac{83}{50}\pi$}
 	\end{subfigure}
 		\begin{subfigure}{\figsizeD}
 		\includegraphics[width=\textwidth]{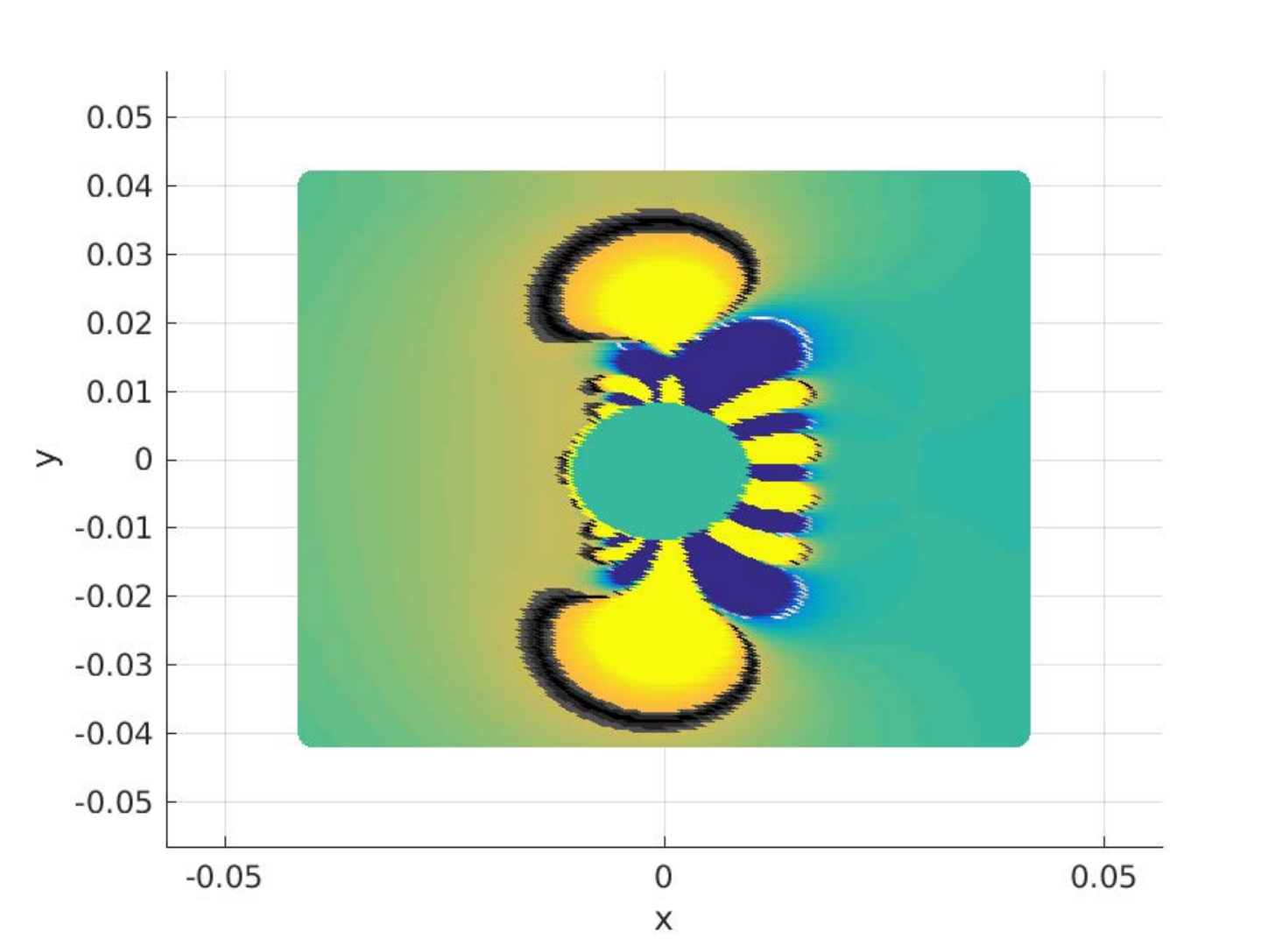}
 		\label{fig:k10h30double_t2_ss84}
 		\caption{$\frac{84}{50}\pi$}
 	\end{subfigure}
 	\begin{subfigure}{\figsizeD}
 	\includegraphics[width=\textwidth]{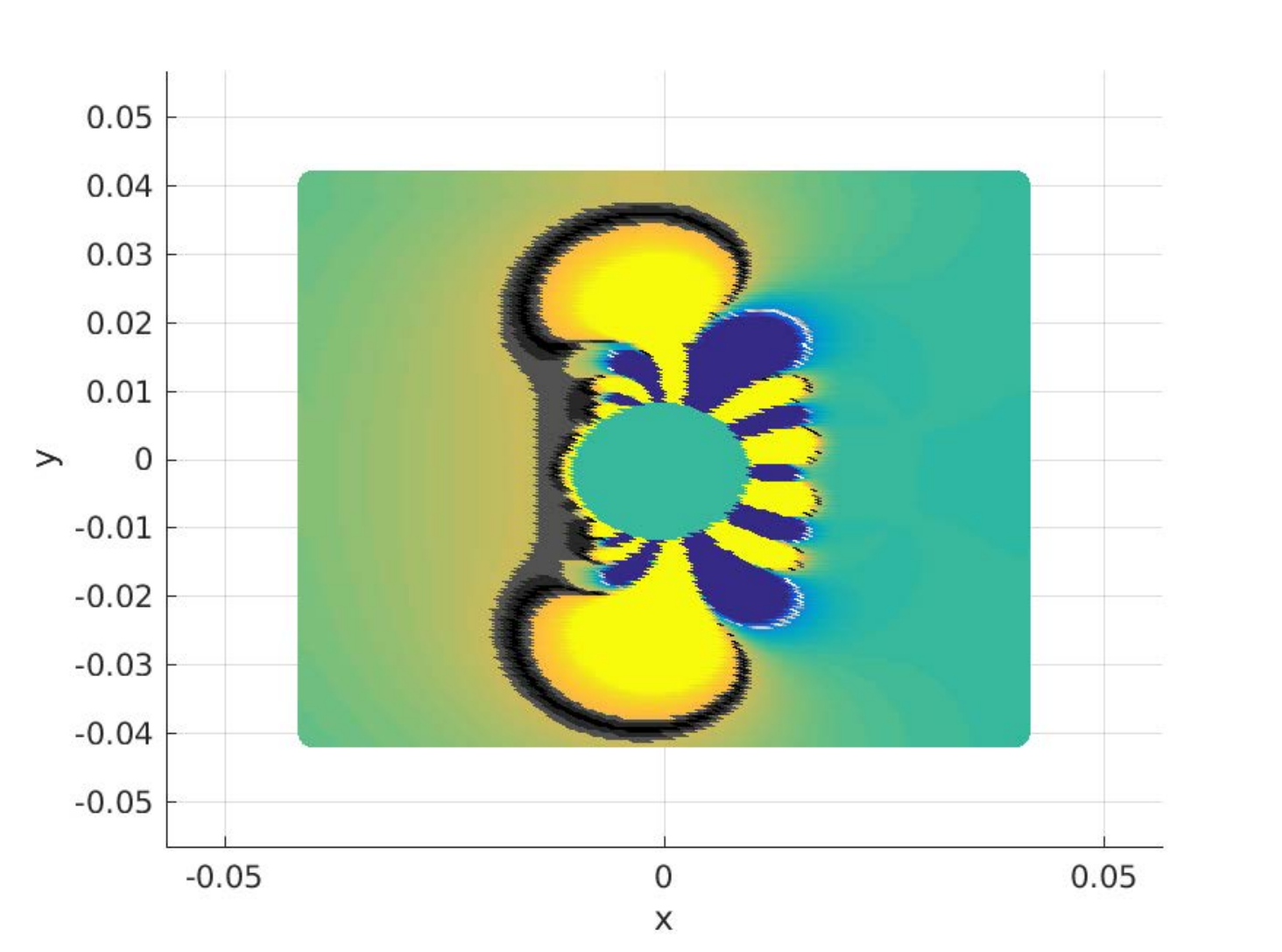}
 	\label{fig:k10h30double_t2_ss85}
 	\caption{$\frac{85}{50}\pi$}
 	\end{subfigure}
 	\begin{subfigure}{\figsizeD}
 	\includegraphics[width=\textwidth]{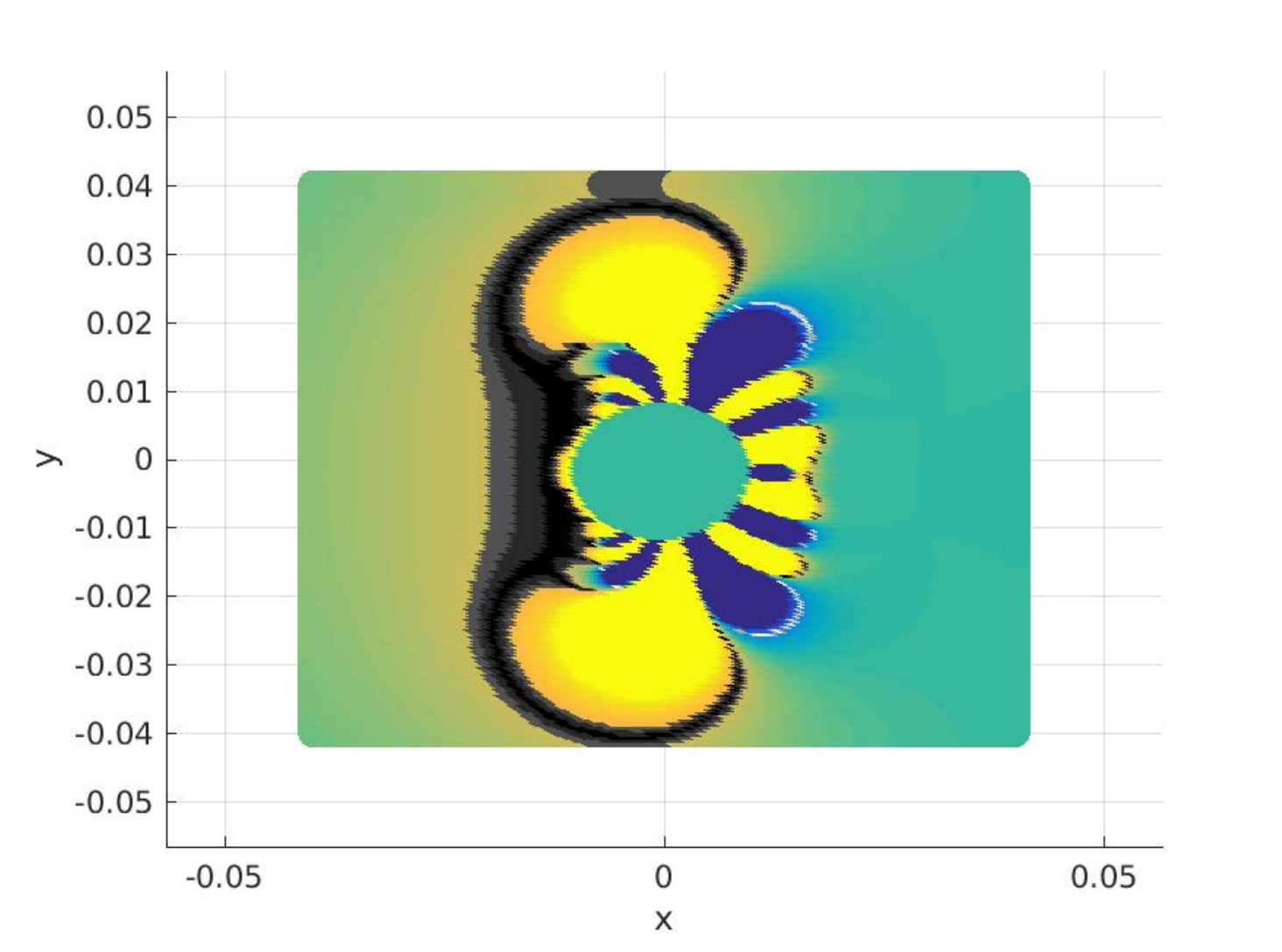}
 	\label{fig:k10h30double_t2_ss86}
 	\caption{$\frac{86}{50}\pi$}
 	\end{subfigure}
 	\begin{subfigure}{\figsizeD}
 		\includegraphics[width=\textwidth]{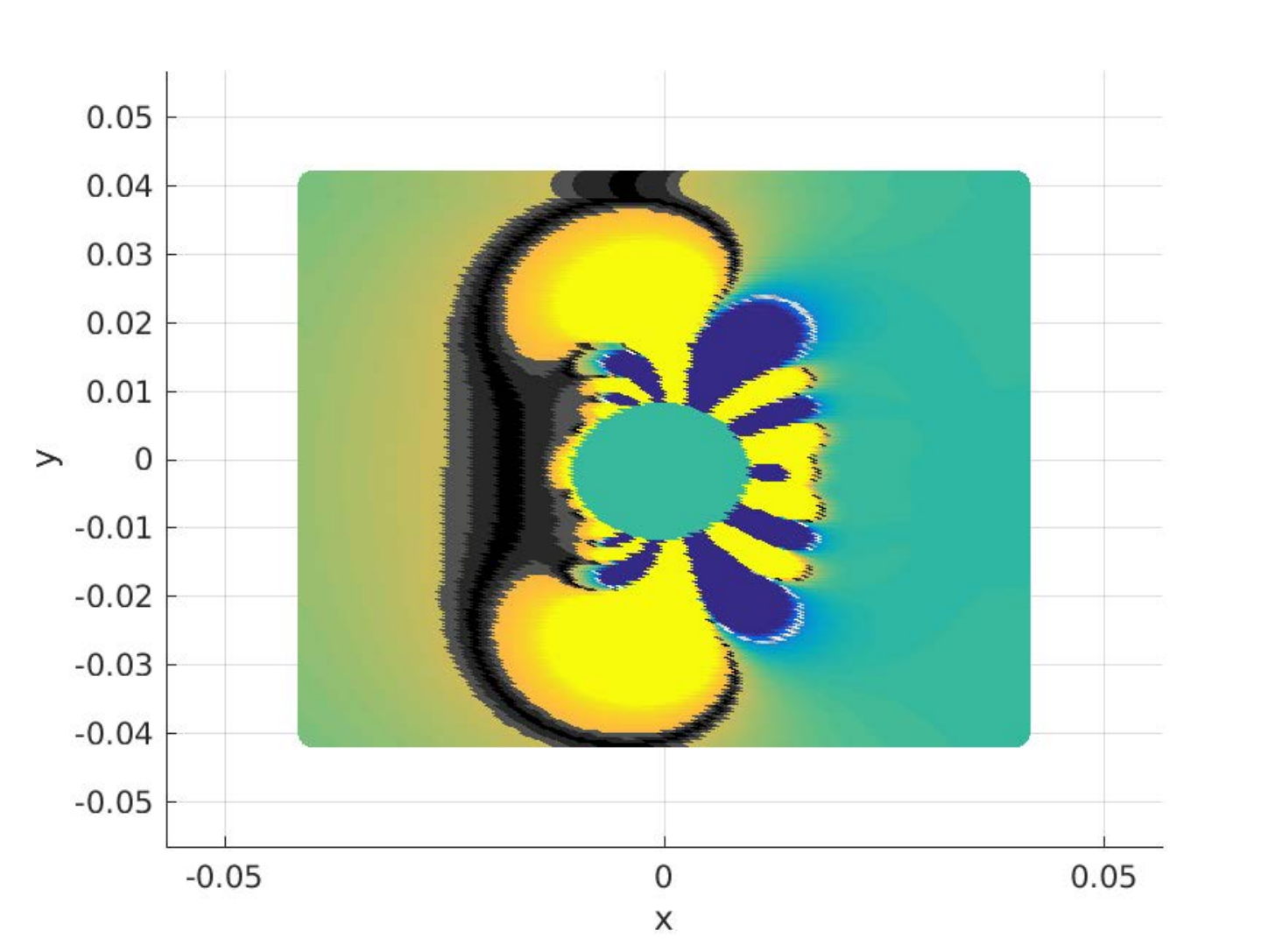}
 		\label{fig:k10h30double_t2_ss87}
 		\caption{$\frac{87}{50}\pi$}
 	\end{subfigure}
 	\begin{subfigure}{\figsizeD}
 		\includegraphics[width=\textwidth]{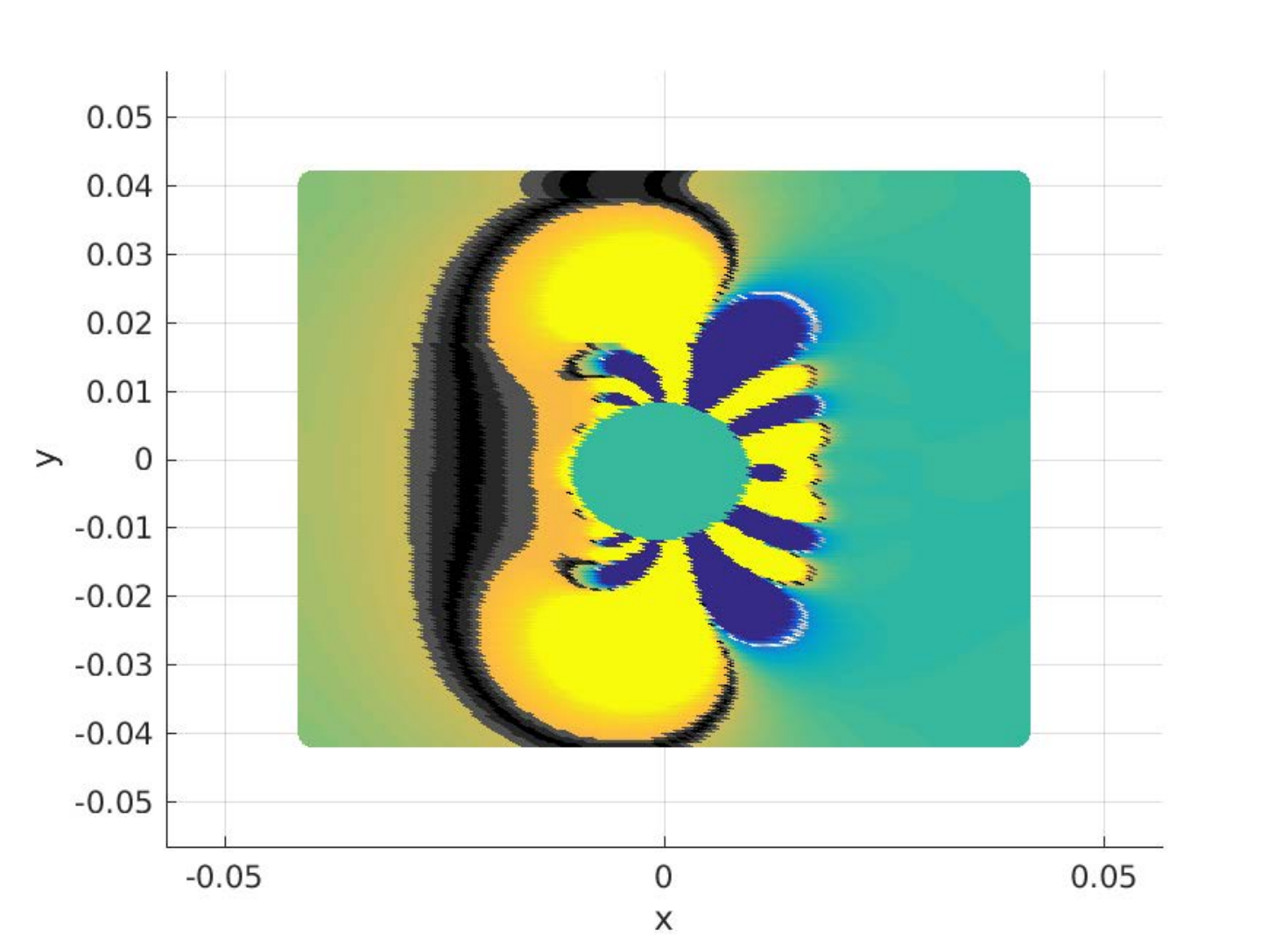}
 		\label{fig:k10h30double_t2_ss88}
 		\caption{$\frac{88}{50}\pi$}
 	\end{subfigure}
 	\caption{Cross-sectional ($z=0$) time snapshots of the propagating generated acoustic field for different values of $kct$.}
 	\label{fig:k10h30double_ss}
 \end{figure}
 
 In Figure \ref{fig:double_ue} we present the quality of our control results in the two regions of interest $D_1, D_2$ as required in \eqref{1-c}. The left and center plots on the top row in the figure describe respectively the field generated by the source, and the outgoing plane wave to be approximated $u_1 = e^{-i x k} $. The good accuracy of our approximation $O(10^{-3})$ can be observed in the right picture on the top row of Figure \ref{fig:double_ue} where the relative pointwise error between $u$ (the field generated solution of \eqref{1}) and $u_1 = e^{-i x k} $ (the field to be approximated) is presented. The fourth picture in Figure \ref{fig:double_ue} (bottom row of the figure)  presents a scattered plot with the values of the generated field in region $D_2$ where very small values of the field (associated with the null effect as required from the optimization procedure) can be observed.
 
 Figure \ref{fig:k10h30double_ss}, shows six a cross-sectional views of the generated field along $z=0$ in a near-field region characterized by $(x,y)\in[-0.1,0.1]^2$. More explicitly, in order left to right and from top left to bottom right plot, we present six cross-sectional ($z=0$) time-snapshots ($kct=\{\frac{83}{50}\pi,\frac{84}{50}\pi,\frac{85}{50}\pi,\frac{86}{50}\pi,\frac{87}{50}\pi,\frac{88}{50}\pi\}$) of the time-harmonic field generated by the synthesized  source in it's near field region, including the two regions of interest $D_1, D_2$. The color scheme in the plots is (truncated to 1 light yellow and -1 dark blue) with the antenna region not included in the numerical simulations and thus corresponding to zero values field (cyan color) and with the black stripe representing field amplitudes of $\approx 0.6$. Following the plots in order from top left to bottom right plot it can be observed how the source works to approximate a plane wave corresponding to a straight black strip in region $D_1$ to its left while maintaining a null in region $D_2$ to its right. Indeed, the plots of Figure \ref{fig:k10h30double_ss} show the time propagation of the generated field by focusing on the portion of the field with amplitude $\approx 0.6$ marked as a dark stripe. It can be observed how this portion of the field enters region $D_1$ at time $kct= \frac{85}{50}\pi$ in a nearly rectilinear shape and continues keeping the same form outgoing throughout a neighbourhood of region $D_1$ (thus indicating the good plane wave structure of the approximated field in the control region) while in all the plots the field is approximately zero in region $D_2$.
 
 The time domain animation
 \label{HL_d_left_t2}
 \href{https://drive.google.com/open?id=0B7nf-pdU3Z4DV3N5cmxsS1l3UE0}{animation 4}
 enhances the message of Figure \ref{fig:k10h30double_ss}. The multimedia file presents the cross-sectional view along $z=0$ of the time-harmonic evolution of the field generated by the synthesized source and respectively the propagating plane wave $u_1= e^{-i x k}e^{-ikct}$ in a near field region given by $(x,y)\in[-0.1,0.1]^2$ in two simultaneous animations: the top one describing the time propagation of the generated field and the bottom one describing the time propagation of the plane wave $u_1= e^{-i x k}e^{-ikct}$. The color scheme in the movies is (truncated to 1 light yellow and -1 dark blue) with the antenna region removed from the simulations (colored cyan) and with the black stripe representing amplitude values around $\approx 0.6$, and the white stripe representing amplitude values around $\approx -0.6$ respectively. As above, we point out that there are two black stripes and respectively two white stripes per period for the approximated plane wave. The animation clearly shows the good accuracy of the approximation in region $D_1$ as well as the null in region $D_2$.

 Figure \ref{fig:double_ad} describes the density $w_\alpha$ (see \eqref{lp}) on the boundary of the fictitious domain $D_{a'}$ as an indication of the possible complexity of the required source inputs $v_n$ or $p_b$ described at \eqref{2} or \eqref{2''}. In the left plot in the figure we present the density values on the surface of $D_{a'}$ viewed from a 3D side perspective and for better visualization we show two more plots in the figure; the center plot shows the density values on the part of the surface facing region $D_1$; the right plot of the figure presents the density values on the part of the surface facing $D_2$.
 \begin{figure}[!htbp] \centering
   \vspace{0cm}
   	\begin{subfigure}{\figsizeC}
   		\includegraphics[width=\textwidth]{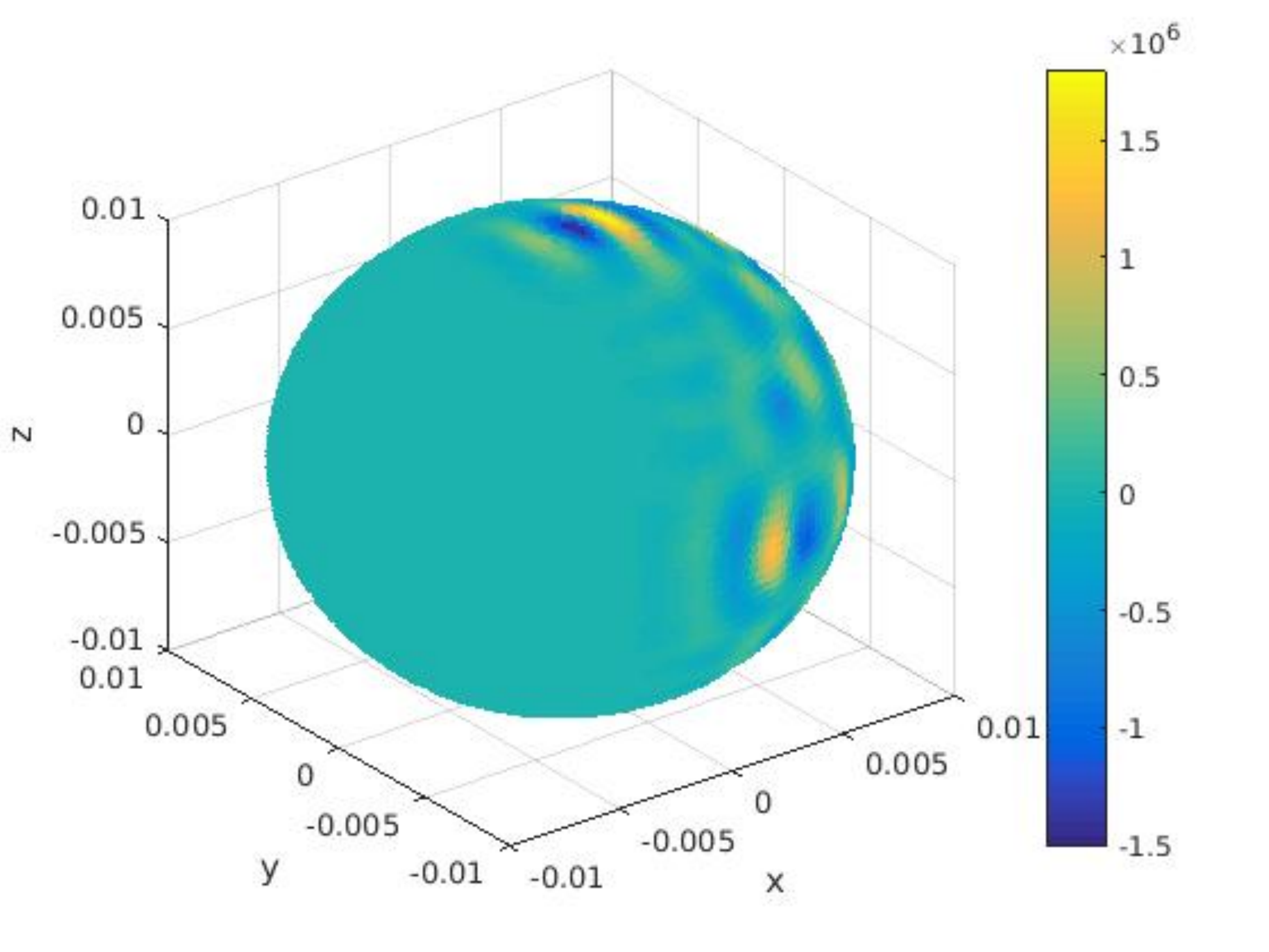}
   		\label{fig:double_ad_s1}
   		\caption{side}
   	\end{subfigure}
   	\begin{subfigure}{\figsizeC}
   		\includegraphics[width=\textwidth]{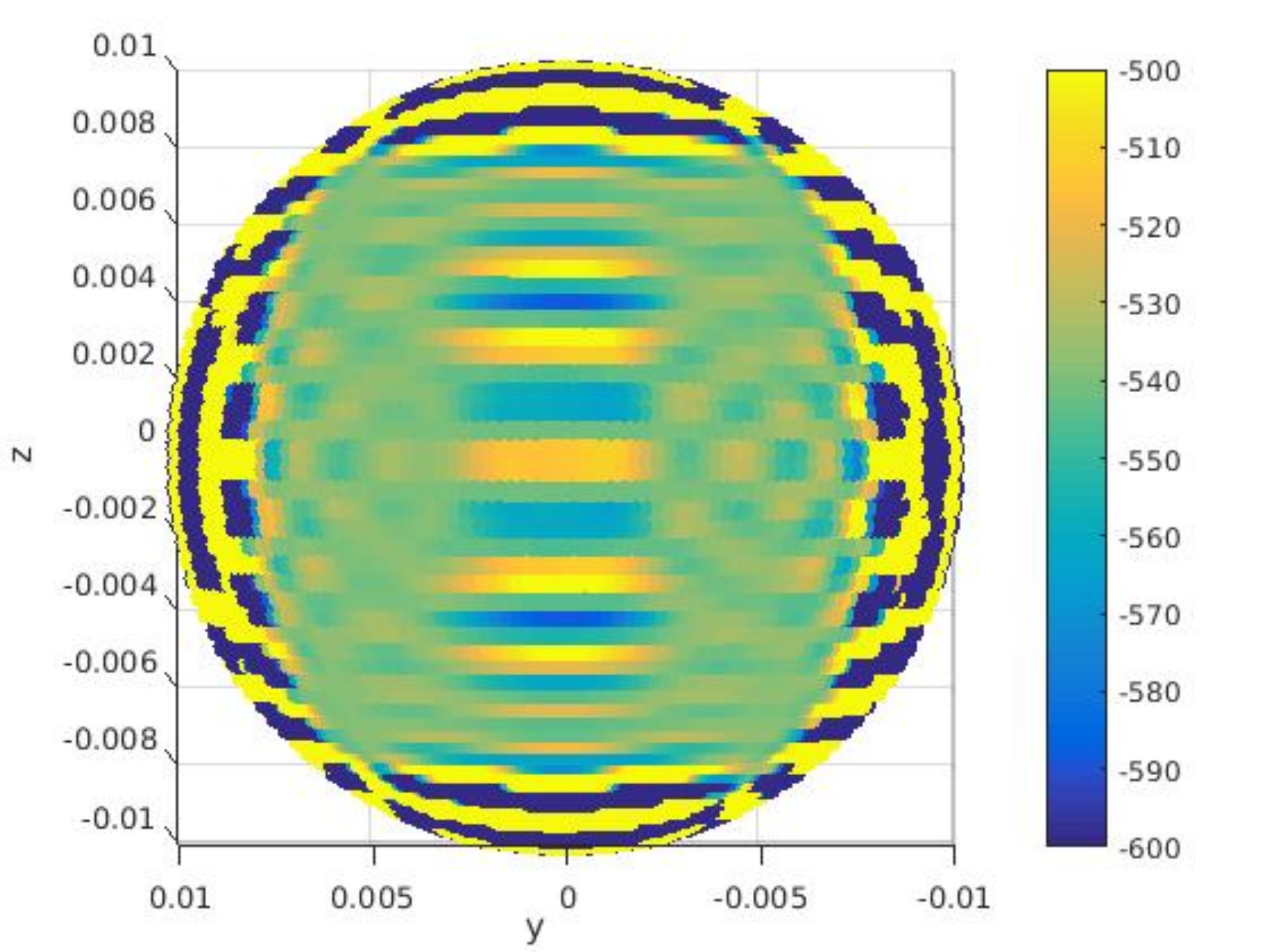}
   		\label{fig:double_ad_b1}
   		\caption{front}
   	\end{subfigure}
   	\begin{subfigure}{\figsizeC}
   		\includegraphics[width=\textwidth]{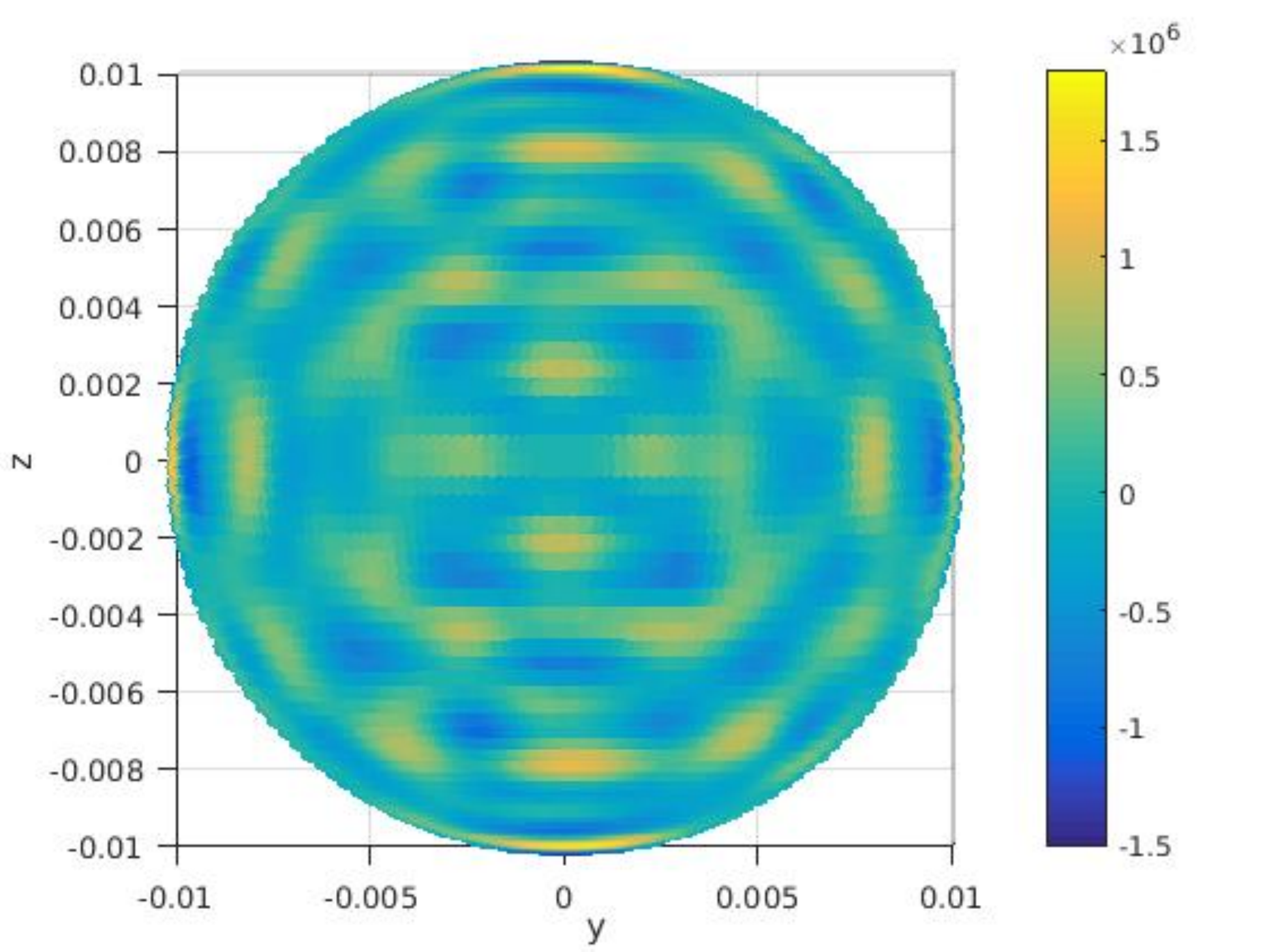}
   		\label{fig:double_ad_a2}
   		\caption{back}
   	\end{subfigure}
   	\caption{Density $w_\alpha$ with various colour maps}
   	\label{fig:double_ad}
   \end{figure}

\subsection{Almost non-radiating acoustic sources with controllable near fields}
\label{II}

In this section we present a third application of the results discussed in Section \ref{fs} and study the problem in the case \eqref{geom-setting} $iii)$. As an extreme example, we show how the source synthesised by our scheme approximates an {\textit{incoming}} plane wave in region $D_1$ while having a very small field beyond radius $R=10$, (i.e. region $D_2$ in this case is the exterior of , $B_{10}(\B0)$, i.e., the ball centered at the origin with radius 10).

Thus, in what follows we consider the case when the source is required to approximate $u_1 = e^{i x k} $ (plane wave propagating towards the source in the positive $x$ direction) with wave number $k=10$ in region $D_1$ described at \eqref{D1} while having a very small field in region $D_2=\RR^3\!\setminus \!B_{10}(\B0)$. A sketch of the geometries are presented in Figure \ref{fig:k10h30F10_g}.
\begin{figure}[!htbp] \centering
\vspace{0cm}
		\includegraphics[width=0.4\textwidth]{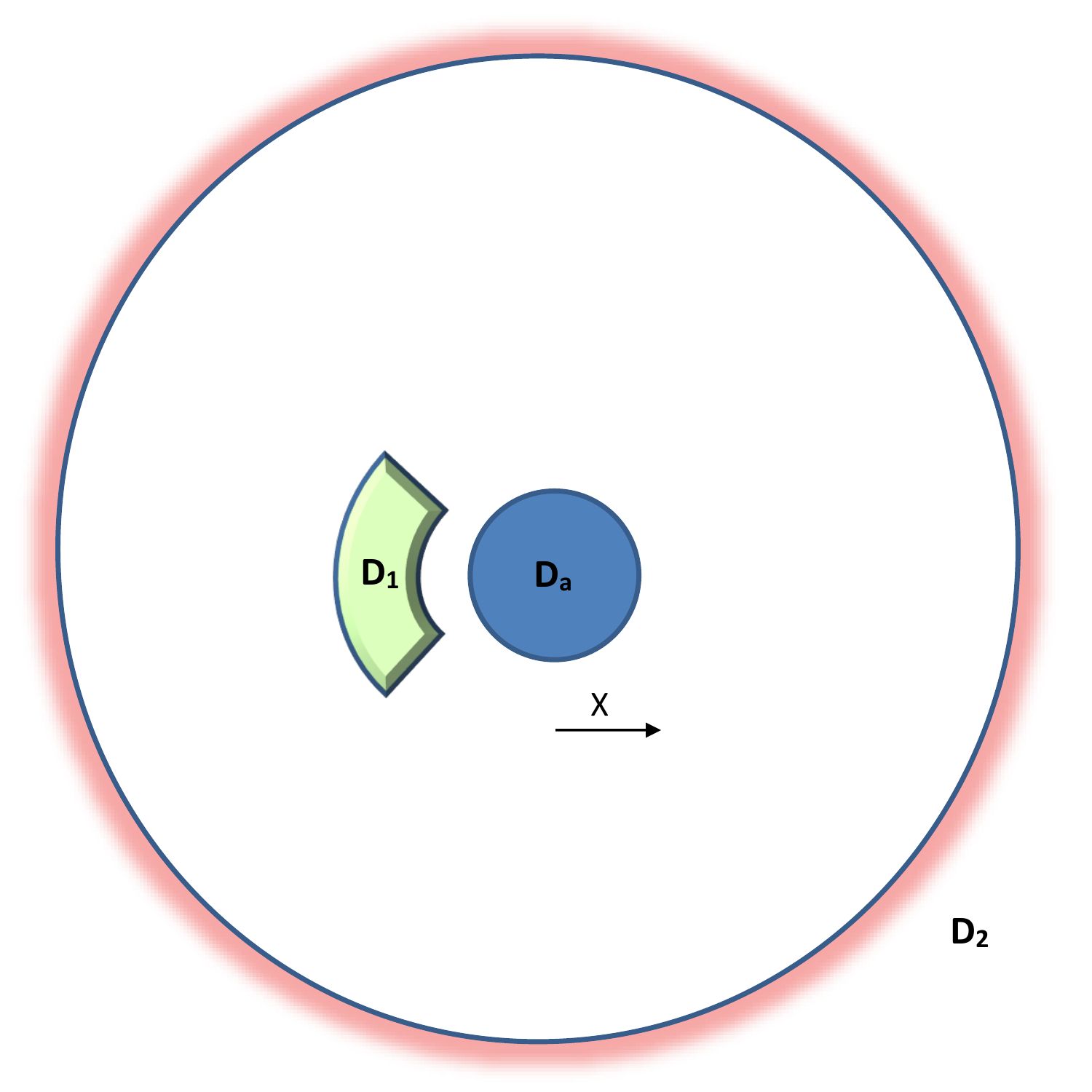}
		\vspace{0cm}
		\caption{Planar sketch of control geometry}
		 \label{fig:k10h30F10_g}
\end{figure}

Figure \ref{fig:k10h30F10_out} shows a cross-sectional view of the generated field along $z=0$ in a region characterized by $(x,y)\in[-5,5]^2$. This plot proves the source causality (i.e., the fact that the synthesized source field is outgoing). This fact can also be observed in the time domain simulation presented in 
\label{HL_anr_left_t4}
\href{https://drive.google.com/open?id=0B7nf-pdU3Z4DMFpaTURnUGtfemc}{animation 5}
where the propagating time-harmonic field generated by the synthesized source is shown. 
\begin{figure}[!htbp] \centering
	\begin{subfigure}{\figsizeC}
		\includegraphics[width=1.5\textwidth]{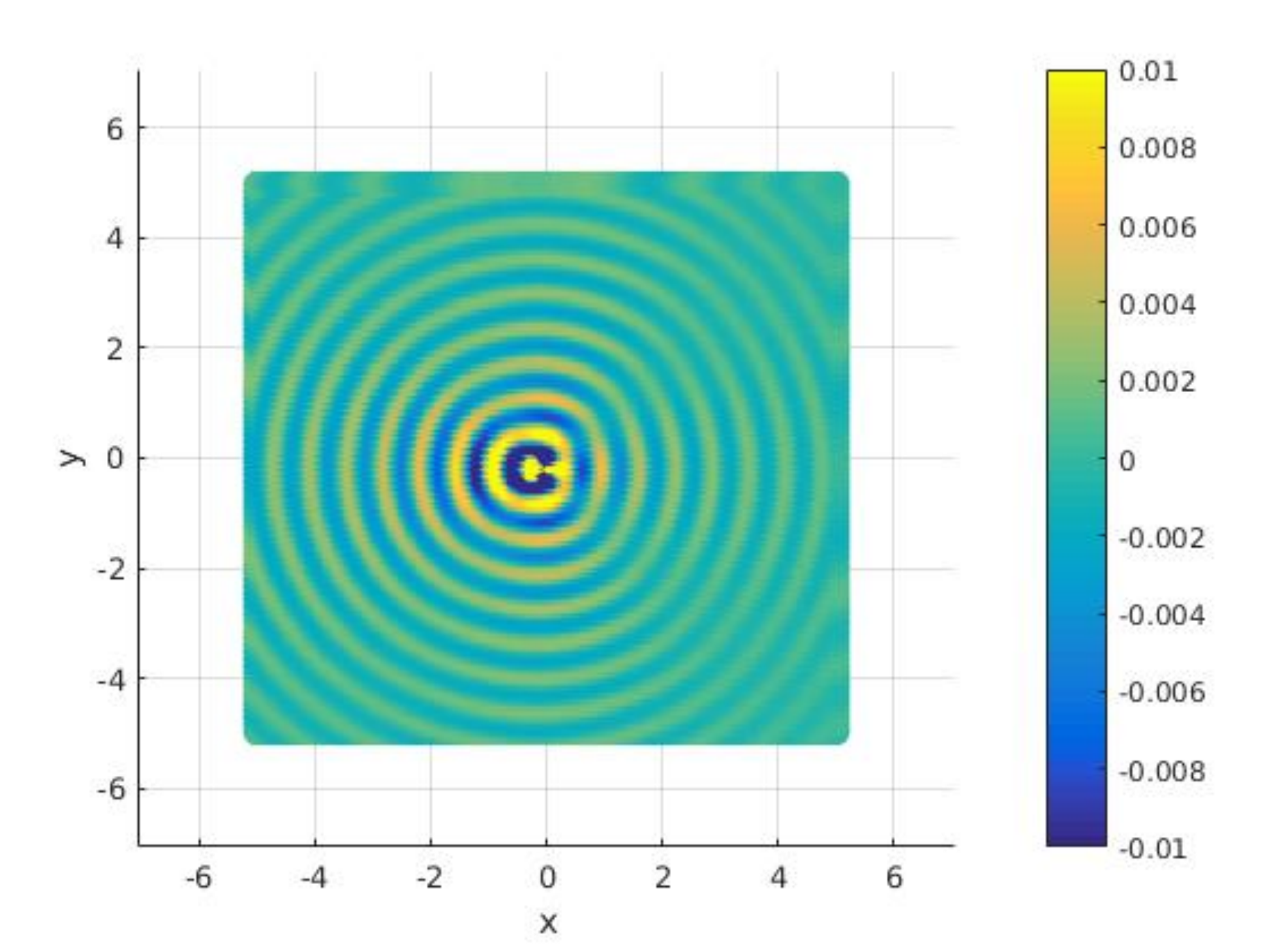} 
	\end{subfigure}
	\caption{Cross-section, with demonstration of radiation decay } 
	\label{fig:k10h30F10_out}
\end{figure}

In Figure \ref{fig:k10h30F10_ue} we show the quality of our control results in region $D_1$ as required by \eqref{1-c}. The left and center plots in the figure describe respectively the field generated by the source $u$, and the plane wave to be approximated $u_1 = e^{i x k}$. The good accuracy of our approximation ( ($O(10^{-3})$) can be observed in the right plot of Figure \ref{fig:k10h30F10_ue} where the relative pointwise error between the synthesised field $u$ and $u_1 = e^{i x k}$ is presented.
\begin{figure}[!htbp] \centering
\vspace{0cm}
	\begin{subfigure}{\figsizeC}
		\includegraphics[width=\textwidth]{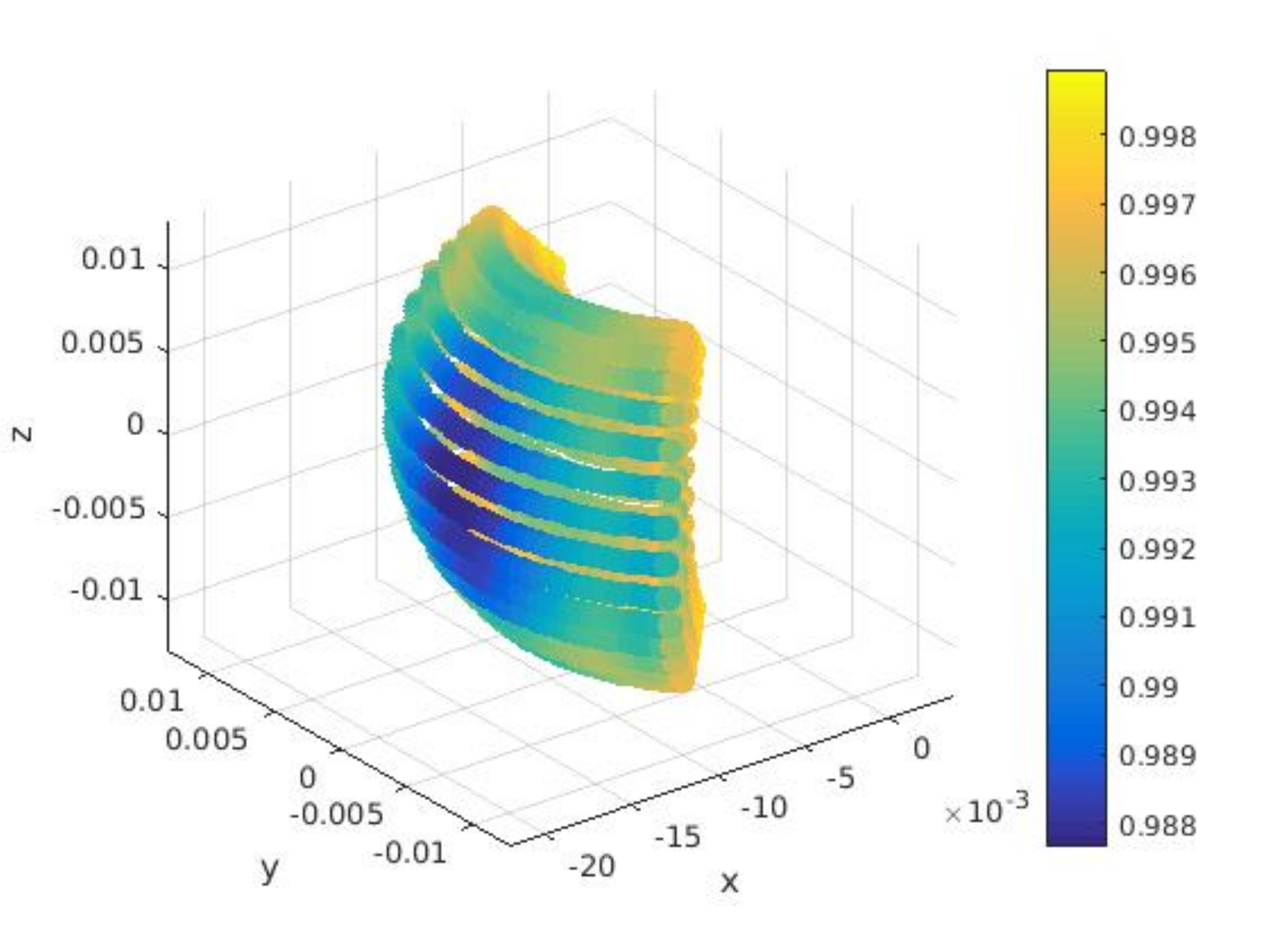}
		\caption{Generated field}
		\label{fig:k10h30F10_nu}
	\end{subfigure}
	\begin{subfigure}{\figsizeC}
		\includegraphics[width=\textwidth]{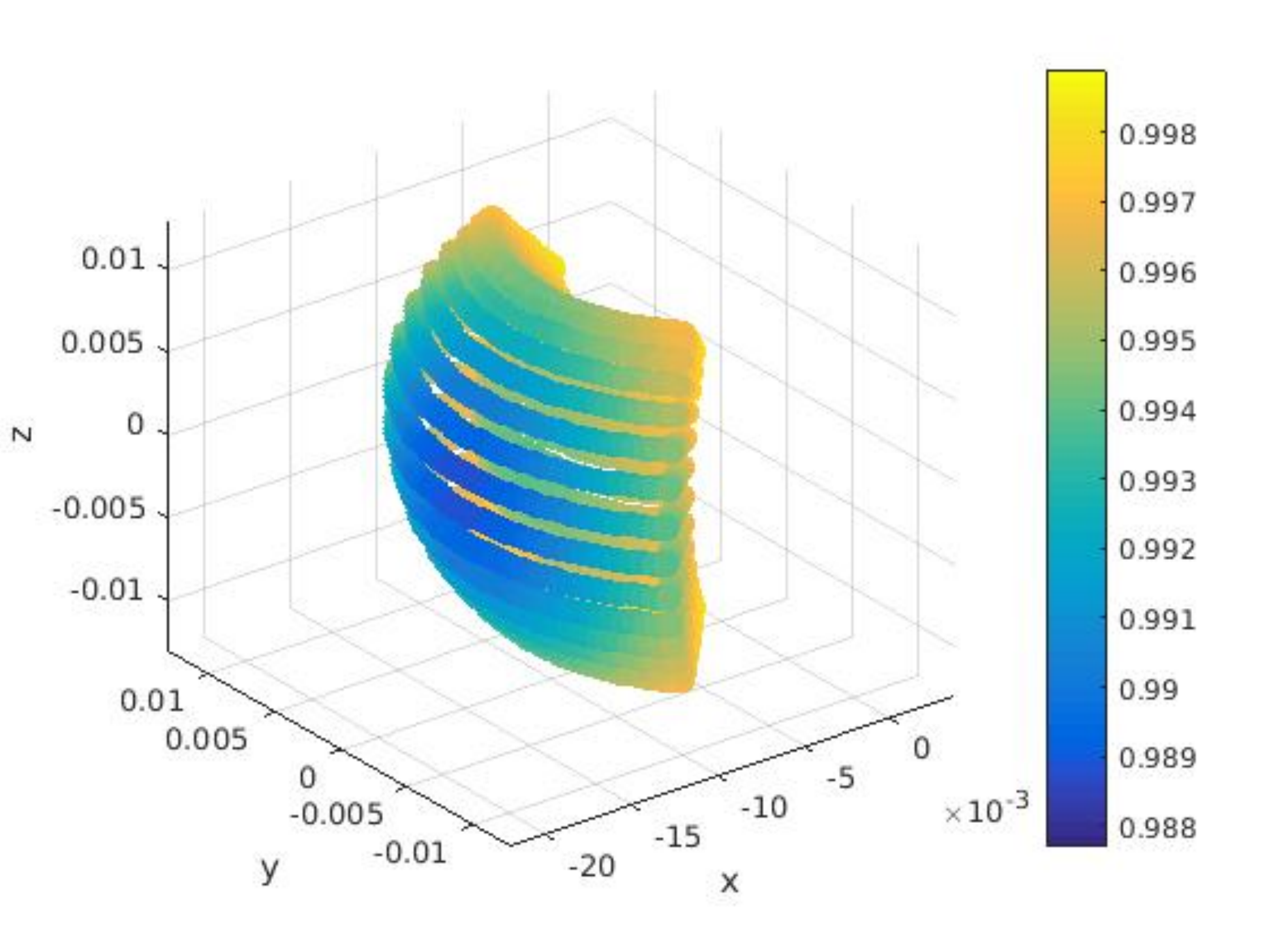}
		\caption{Incident field}
		\label{fig:k10h30F10_nf}
	\end{subfigure}
	\begin{subfigure}{\figsizeC}
		\includegraphics[width=\textwidth]{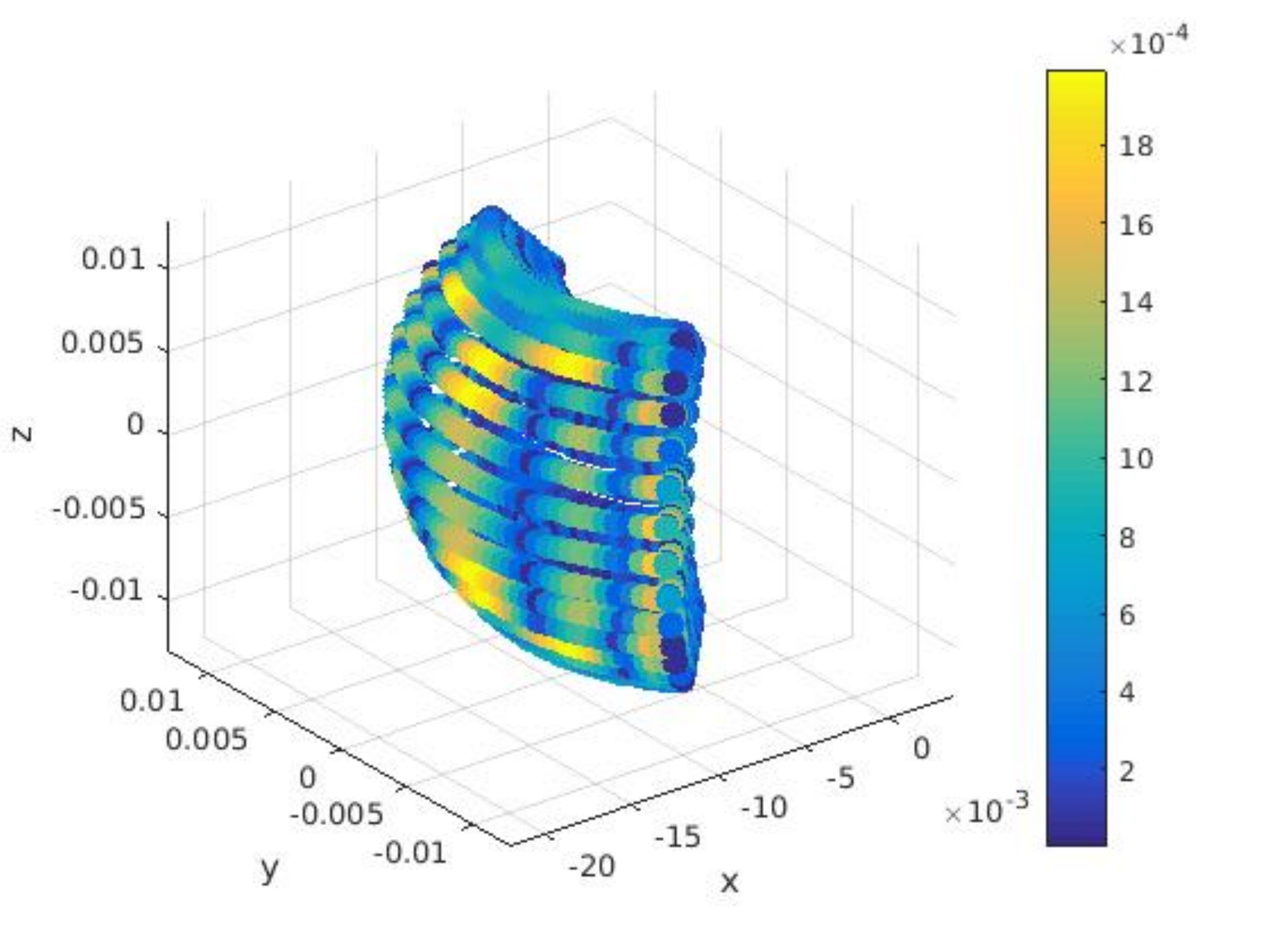}
		\caption{Pointwise relative error}
		\label{fig:k10h30F10_ne}
	\end{subfigure}
	\vspace{0cm}
	\caption{Accuracy of control in region $D_1$.}
	\label{fig:k10h30F10_ue}
\end{figure}

Figure \ref{fig:k10h30F10_far} shows the fast decay in region $D_2$ as required by \eqref{1-c}. Indeed, the left plot of the figure describes the very small values of the generated field computed on the sphere of radius 10. On the other hand, the right plot in Figure \ref{fig:k10h30F10_far} describes the absolute values of $\displaystyle r\cdot\sup_{B_r(\B0)}|u|$ as a function of $r\in(10,1000)$. The asymptotic limit of this function, $O(10^{-2})$, is the supremum value of far field pattern  and this once more confirms the fact that the source synthesized by our scheme is a weak radiator. In fact, we also computed the actual power radiated by this source, i.e., $P=Re(\int_S u^*\Bv\cdot \Bn)$ where $u$ represents the pressure field solution of \eqref{1}, \eqref{1-c} , $A^*$ denotes the complex conjugate of complex quantity $A$, and $\Bv\cdot\Bn$ denotes the normal velocity on a sphere $S$ surrounding the source $D_a$, and we found that it is of order $O(10^{-7})$, once more indicating a very weak radiator.
\begin{figure}[!htbp] \centering
	\begin{subfigure}{\figsizeC}
		\includegraphics[width=1.1\textwidth]{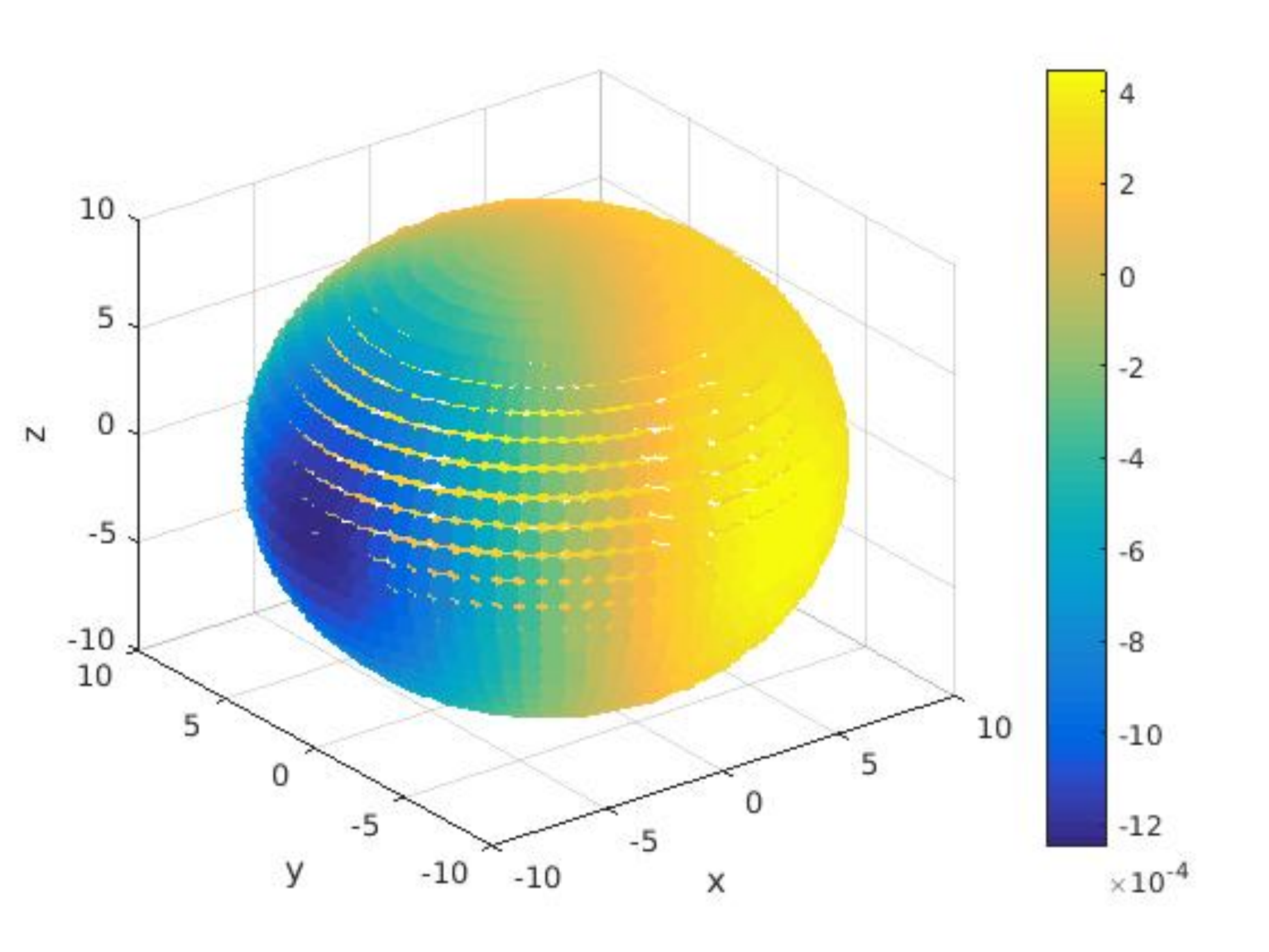}
		\caption{Generated field on $R=10$}
		\label{fig:k10h30F10_fu}
	\end{subfigure}
\hspace{2cm}
	\begin{subfigure}{\figsizeC}
		\includegraphics[width=1.1\textwidth]{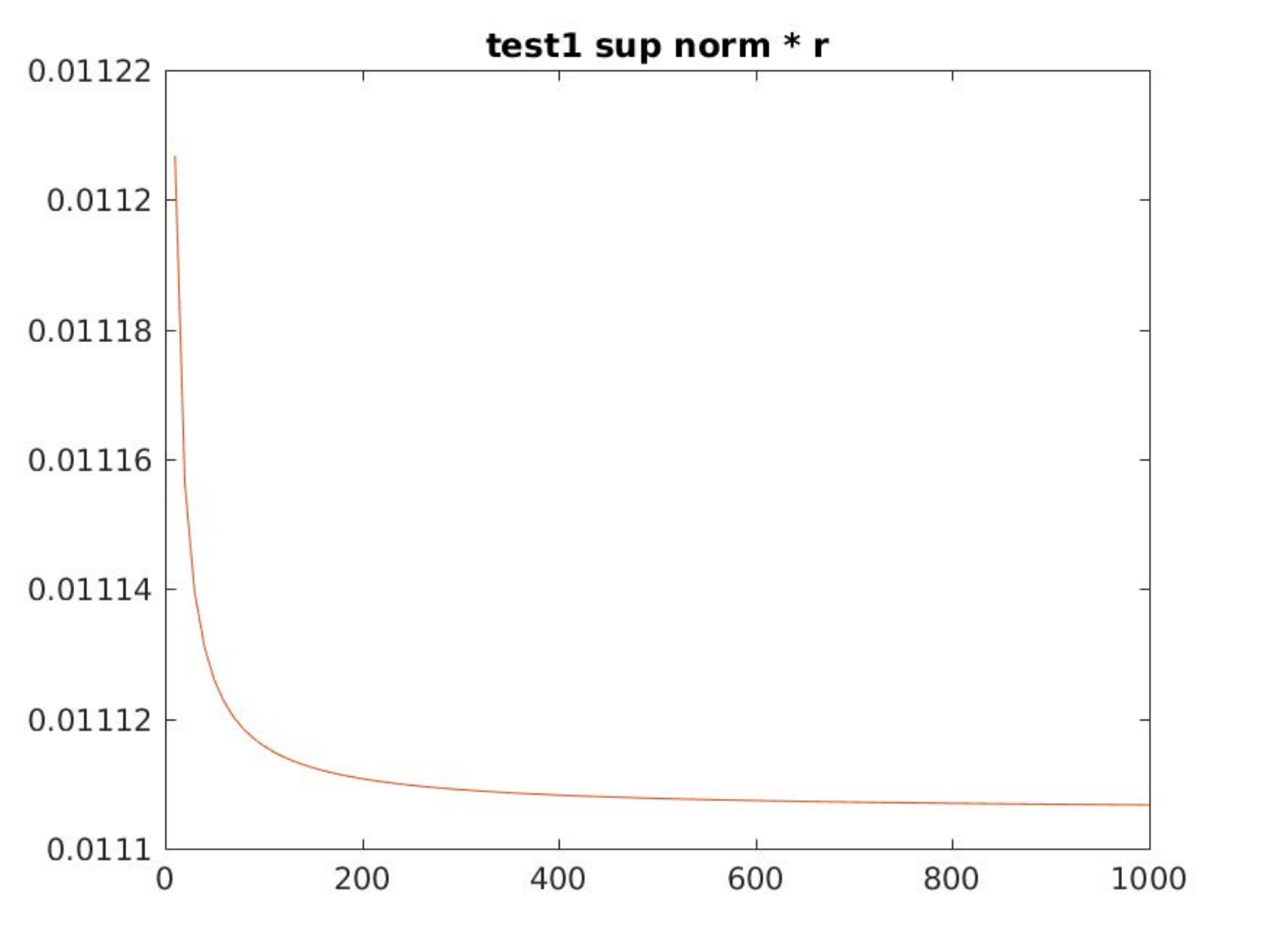}
		\caption{Far-field pattern convergence}
		\label{fig:k10h30F10_decayB}
	\end{subfigure}
	\caption{Far field pattern}
	\label{fig:k10h30F10_far}
\end{figure}
\begin{figure}[!htbp] \centering
\vspace{0cm}
	\begin{subfigure}{\figsizeD}
		\includegraphics[width=\textwidth]{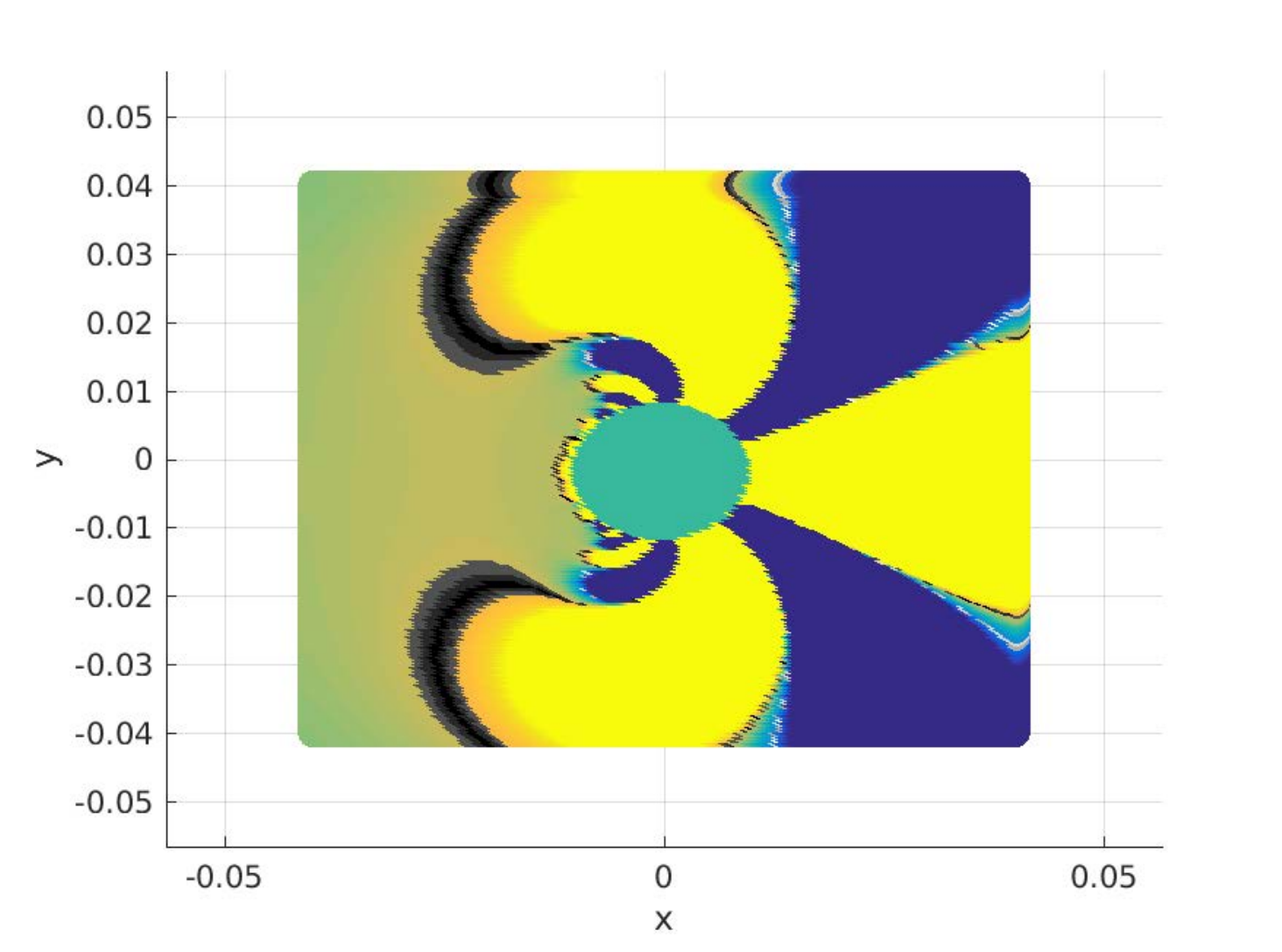}
		\label{fig:k10h30F10_t2_ss79}
		\caption{$\frac{79}{50}\pi$}
	\end{subfigure}
	\begin{subfigure}{\figsizeD}
		\includegraphics[width=\textwidth]{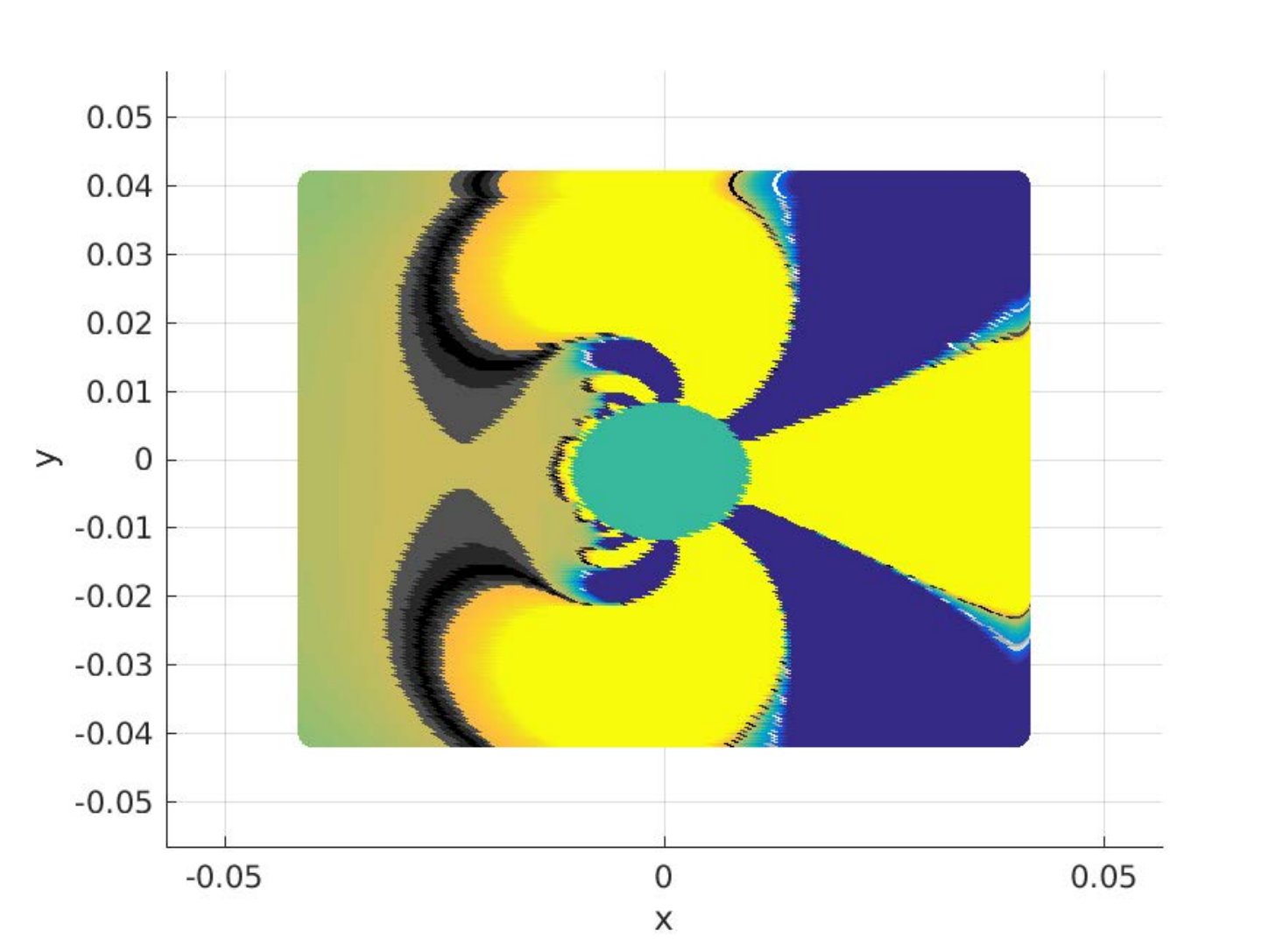}
		\label{fig:k10h30F10_t2_ss80}
		\caption{$\frac{80}{50}\pi$}
	\end{subfigure}
	\begin{subfigure}{\figsizeD}
		\includegraphics[width=\textwidth]{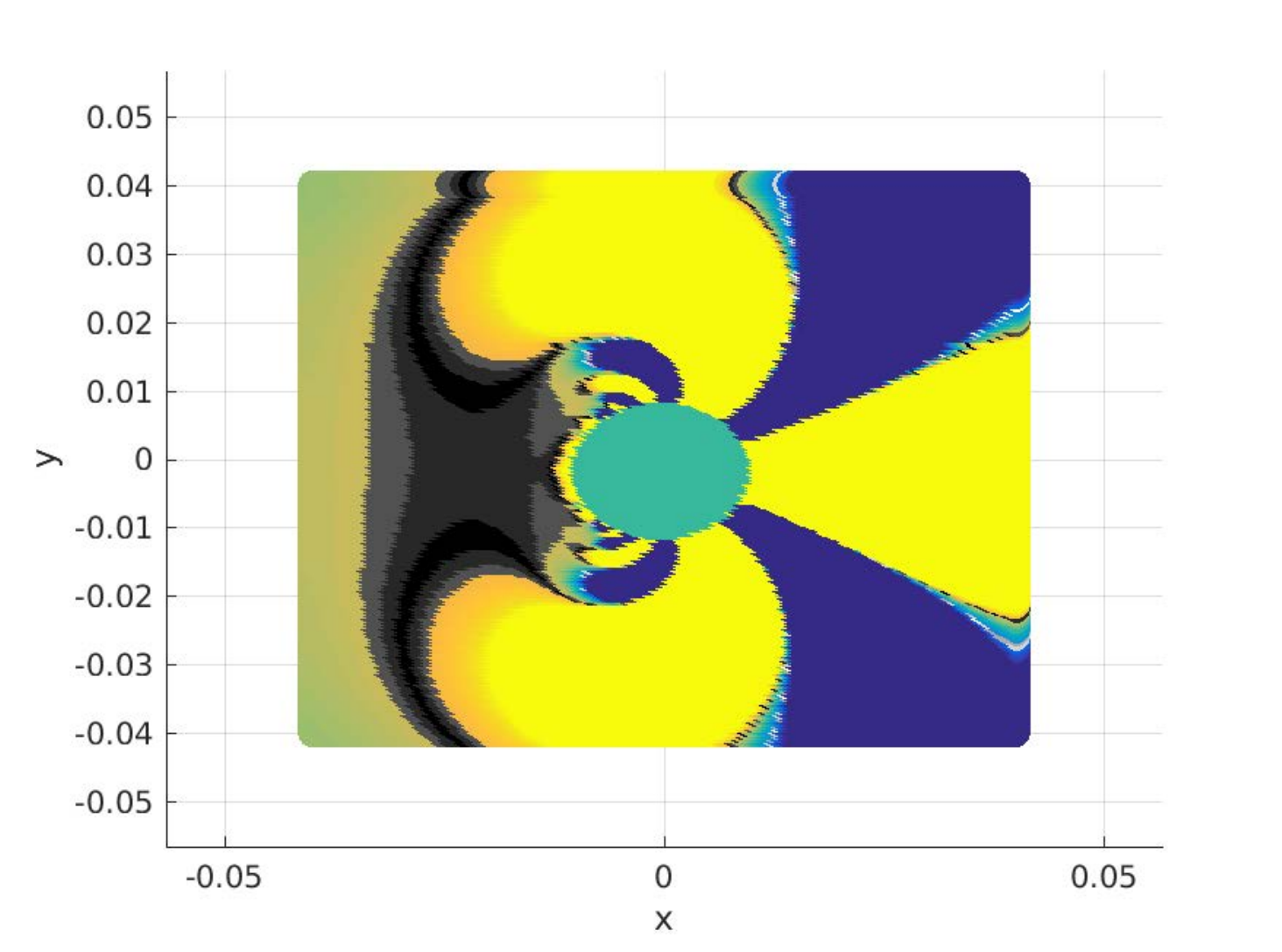}
		\label{fig:k10h30F10_t2_ss81}
		\caption{$\frac{81}{50}\pi$}
	\end{subfigure}
	\begin{subfigure}{\figsizeD}
		\includegraphics[width=\textwidth]{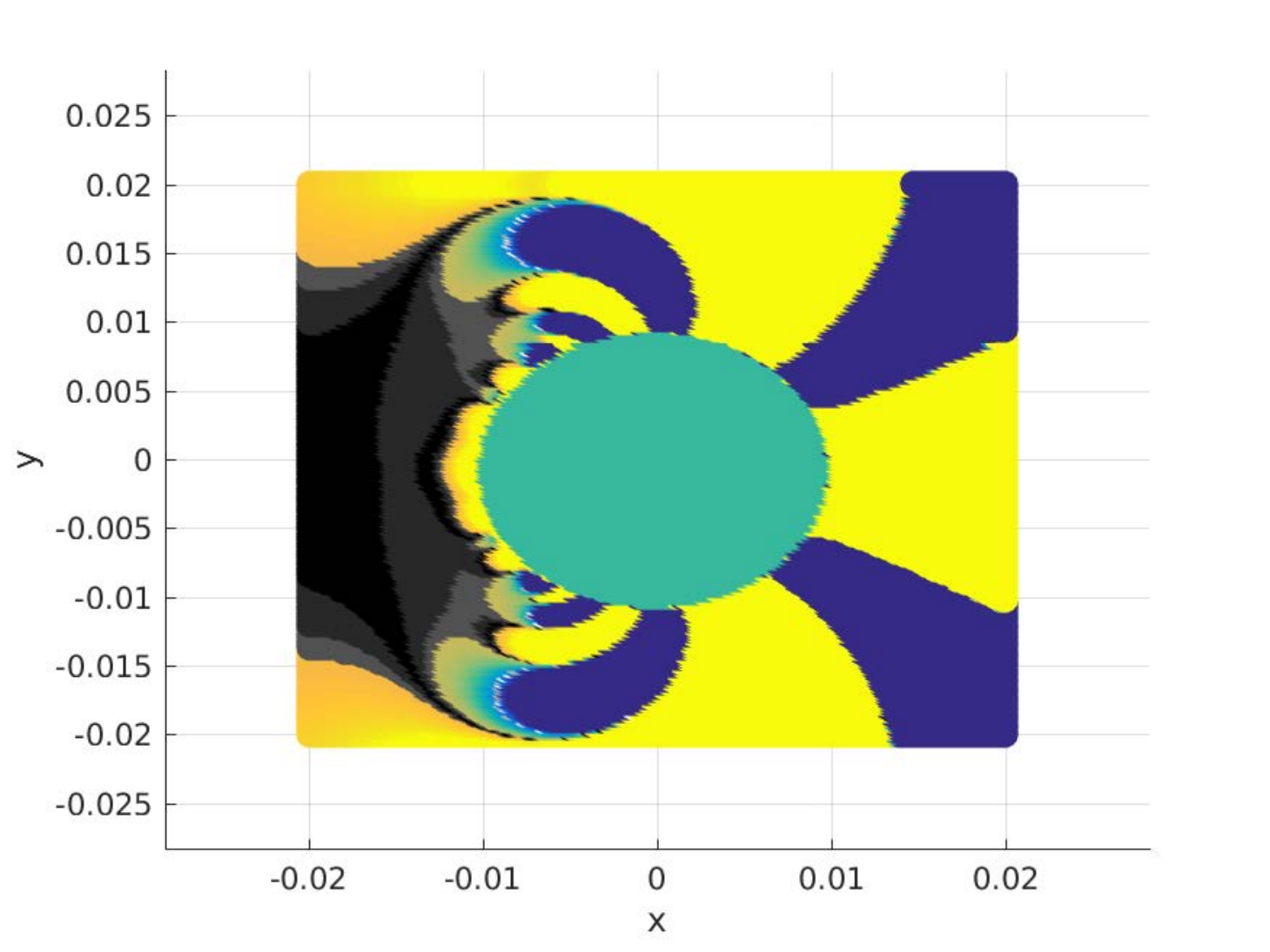}
		\label{fig:k10h30F10_t2_ss81.5}
		\caption{$\frac{81.5}{50}\pi$}
	\end{subfigure}
	\begin{subfigure}{\figsizeD}
		\includegraphics[width=\textwidth]{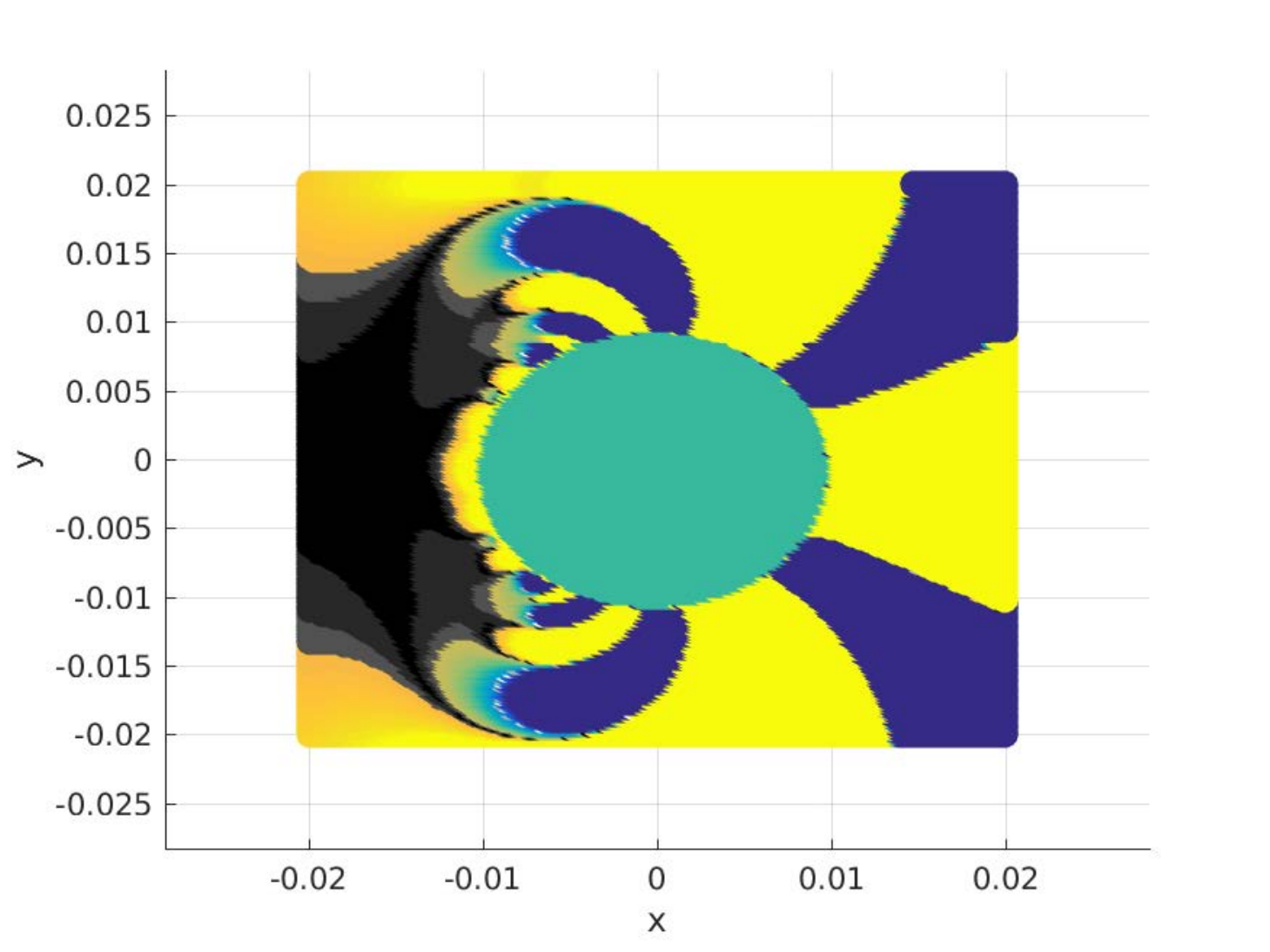}
		\label{fig:k10h30F10e_t2_ss81.7}
		\caption{$\frac{81.7}{50}\pi$}
	\end{subfigure}
	\begin{subfigure}{\figsizeD}
		\includegraphics[width=\textwidth]{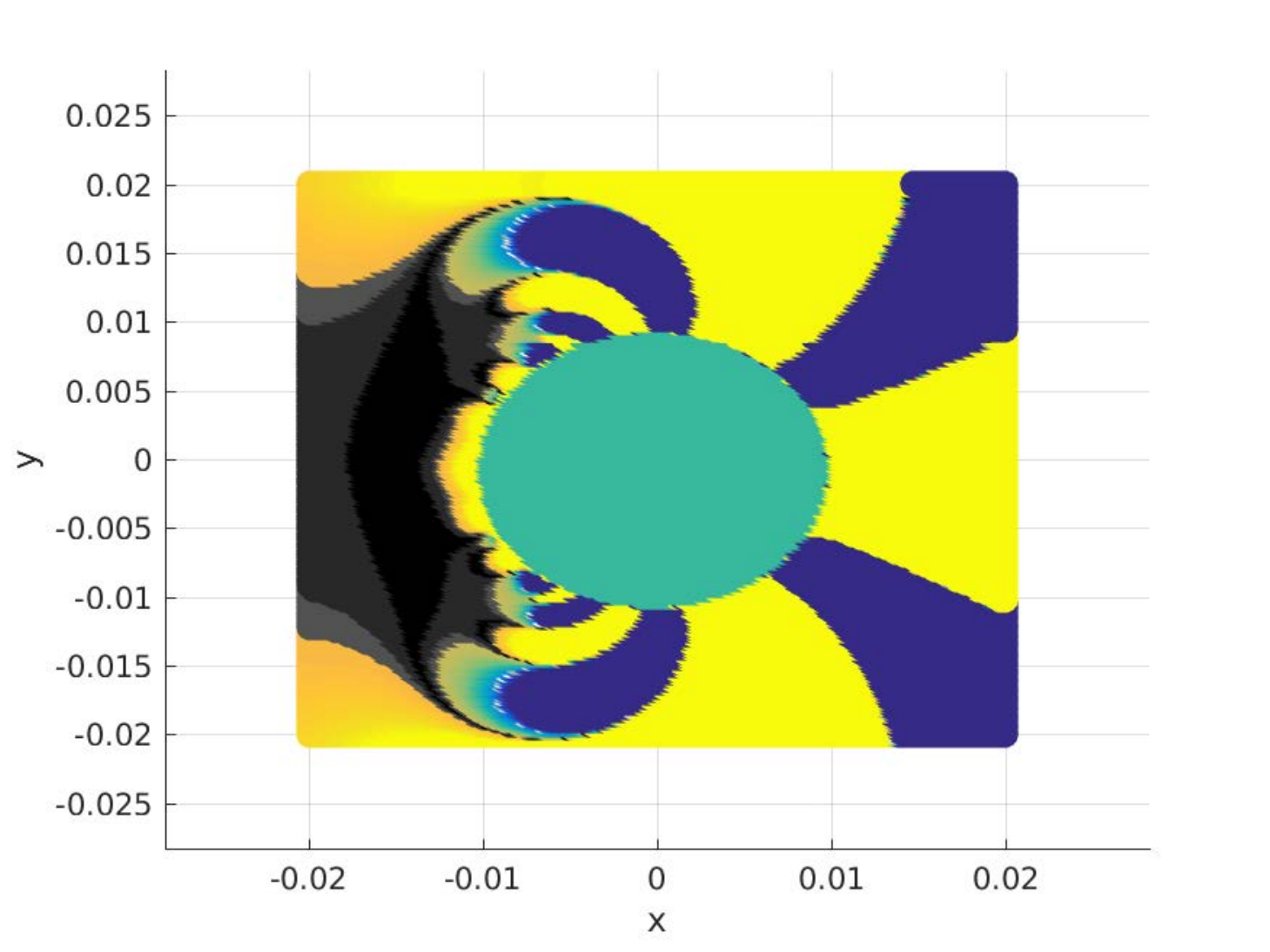}
		\label{fig:k10h30F10_t2_ss82}
		\caption{$\frac{82}{50}\pi$}
	\end{subfigure}
	\begin{subfigure}{\figsizeD}
	\includegraphics[width=\textwidth]{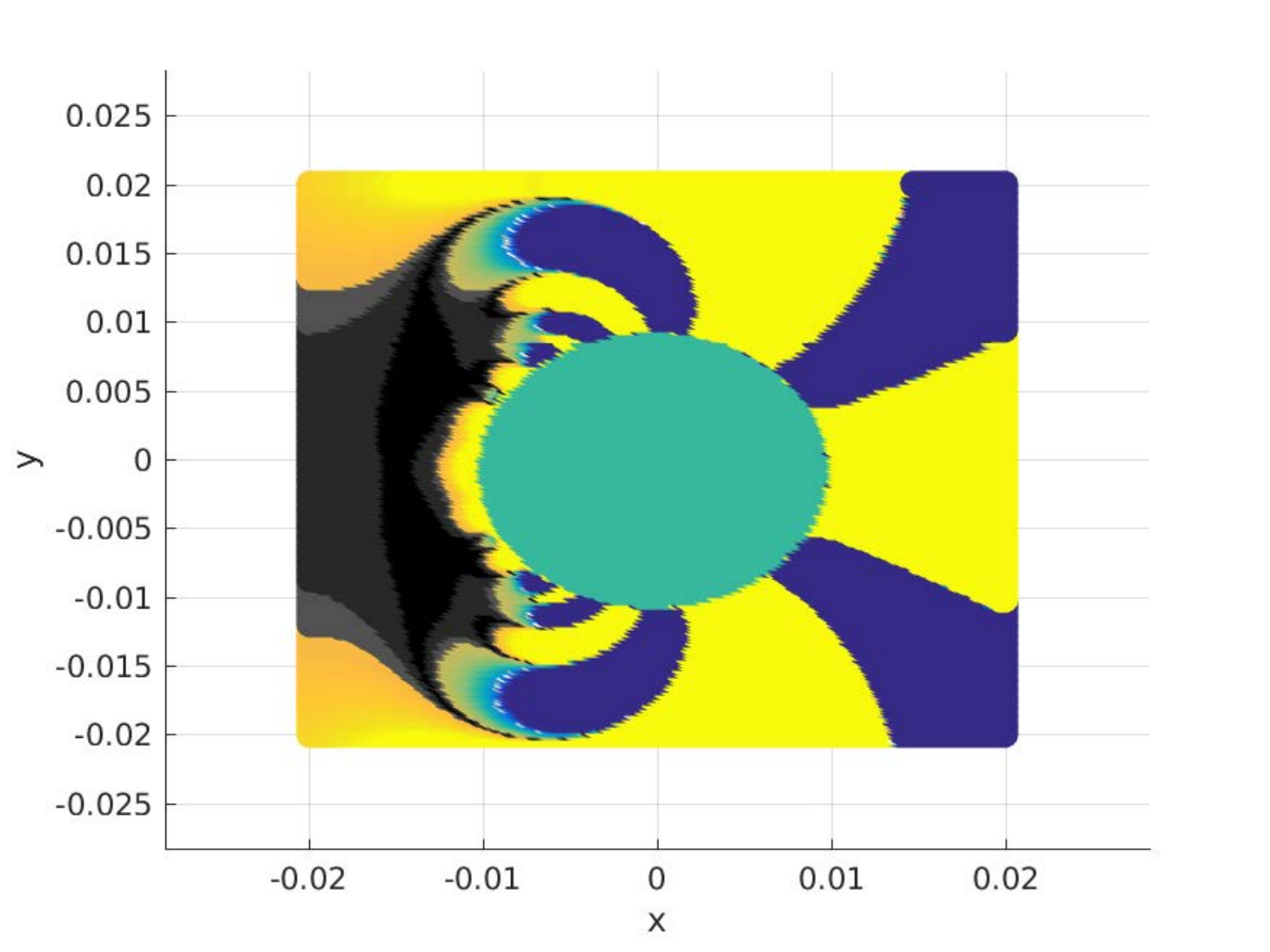}
	\label{fig:k10h30F10_t2_ss82.1}
	\caption{$\frac{82.1}{50}\pi$}
		\end{subfigure}
	\begin{subfigure}{\figsizeD}
	\includegraphics[width=\textwidth]{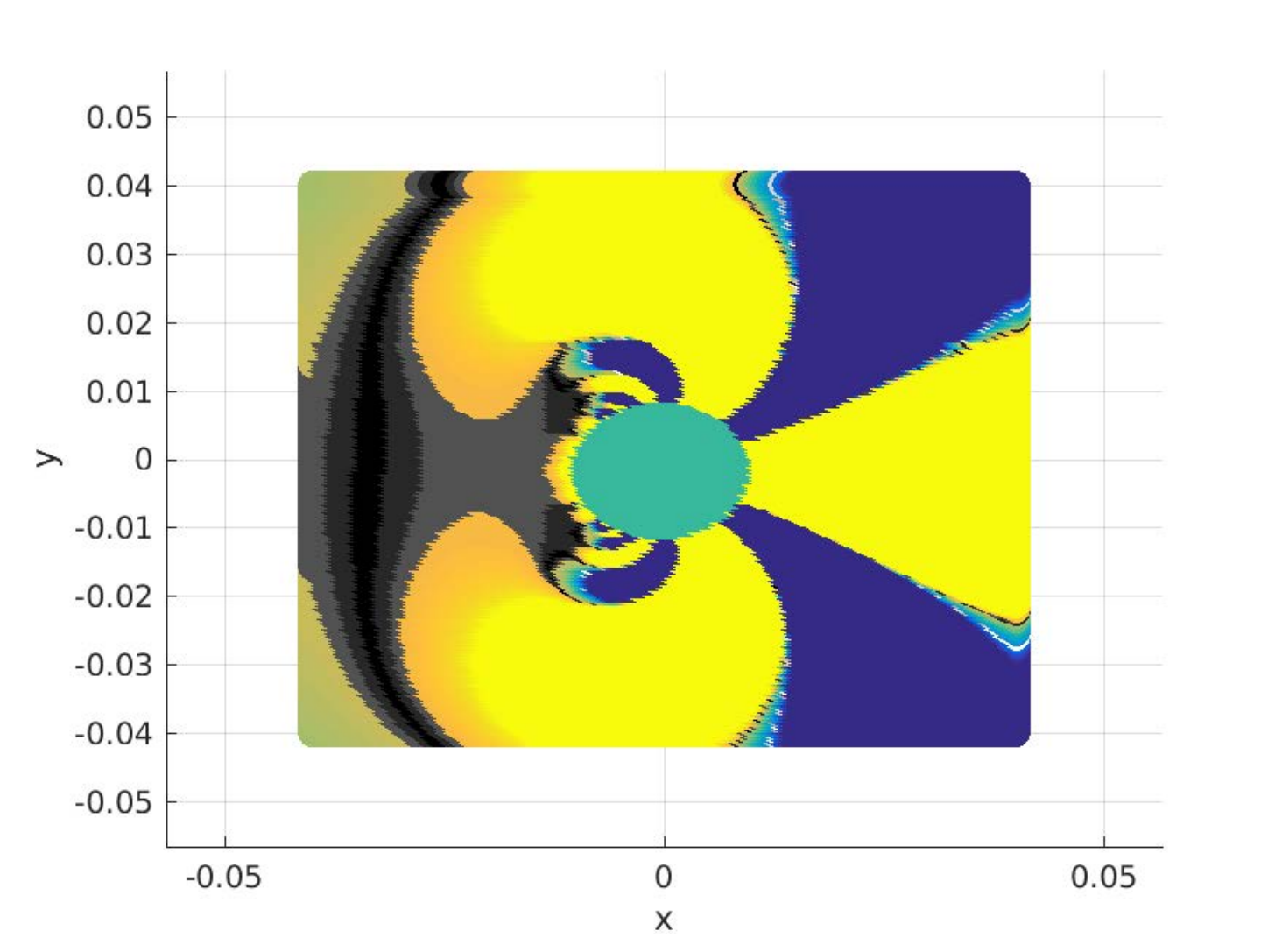}
	\label{fig:k10h30F10_t2_ss83}
	\caption{$\frac{83}{50}\pi$}
	\end{subfigure}
	\begin{subfigure}{\figsizeD}
		\includegraphics[width=\textwidth]{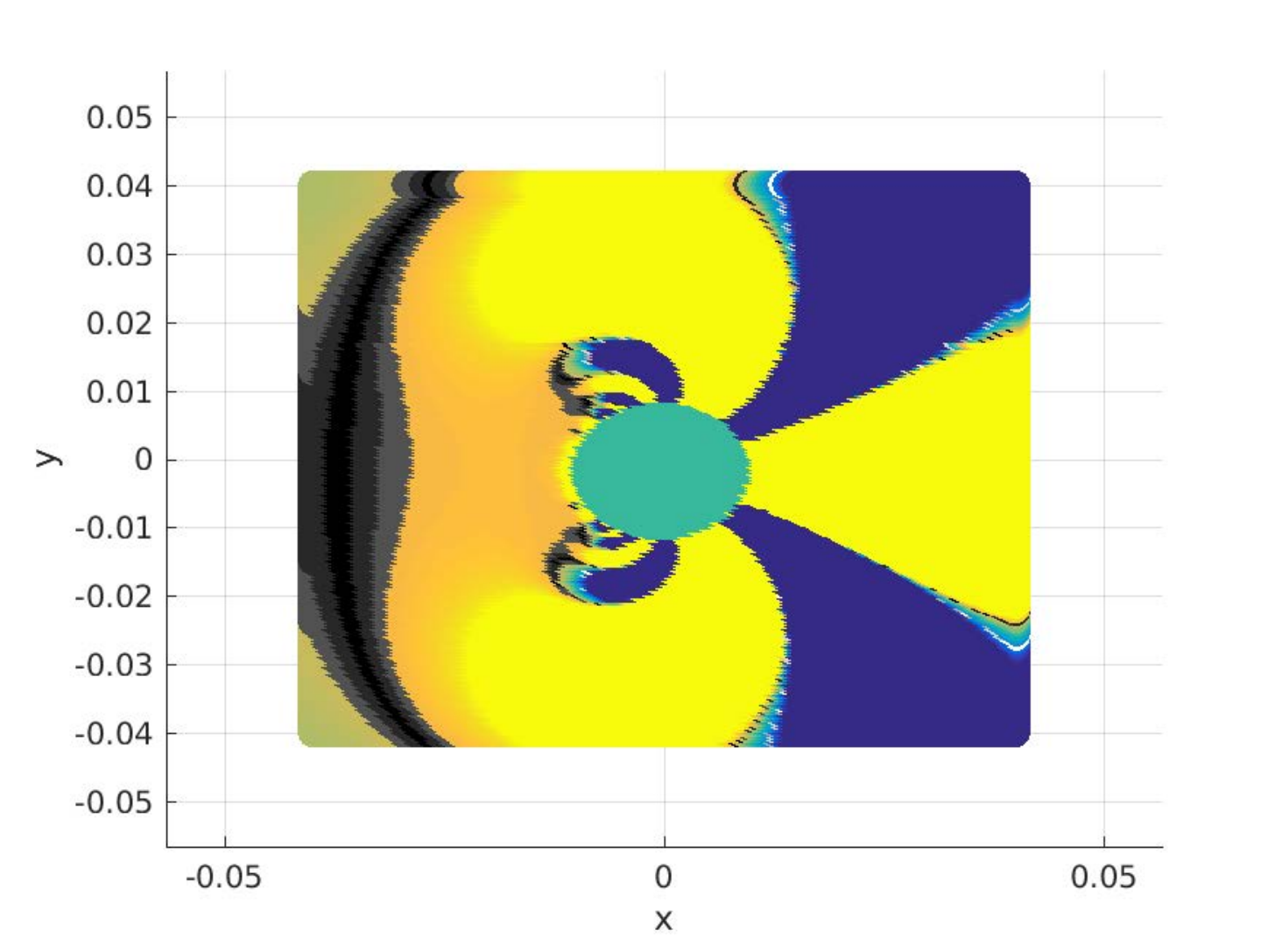}
		\label{fig:k10h30F10_t2_ss84}
		\caption{$\frac{84}{50}\pi$}
	\end{subfigure}
	\caption{Cross-sectional ($z=0$) time snapshots of the propagating generated acoustic field for different values of $kct$.}
	\label{fig:k10h30F10_ss}
\end{figure}

Figure \ref{fig:k10h30F10_ss} shows nine a cross-sectional views of the generated field along $z=0$ in a near-field region characterized by $(x,y)\in[-0.1,0.1]^2$.
The figure describes in order left to right from top left to bottom right plot, nine cross-sectional ($z=0$) time-snapshots ($kct=\{\frac{79}{50}\pi,\frac{80}{50}\pi,\frac{81}{50}\pi,\frac{81.5}{50}\pi,\frac{81.7}{50}\pi,\frac{82}{50}\pi,\frac{82.1}{50}\pi, \frac{83}{50}\pi, \frac{84}{50}\pi\}$) of the time-harmonic field generated by the synthesized  source in it's near field region. The color scheme in the plots is (truncated to 1 light yellow and -1 dark blue) with the antenna region (coloured cyan) not included in the numerical simulations and with the black stripe representing field amplitudes of $\approx0.6$. 

Following the plots in order from top left to bottom right plot it can be observed how the source works to approximate the incoming plane wave $u_1 = e^{i x k}e^{-ikct}$ in $D_1$. Indeed, the $(a), (b)$ plots show how the source works on creating a plane wave in region $D_1$. The $(c), (d), (e), (f)$ plots are zoomed in closer to the antenna in a region $(x,y)\in[-0.02,0.02]^2$ so that the approximation of the incoming plane wave is better observed. It can be seen in these plots how the portion of the fields with values $\approx 0.6$ (black stripe) enters region $D_1$ at time $kct= \frac{81.5}{50}\pi$ in a nearly rectilinear shape (i.e. corresponding to plane wave character) and propagates towards the source while keeping the same rectilinear throughout a neighbourhood of region $D_1$. In the last two plots we see how the generated field looses form near the source and propagates away from it in the far field (i.e., corresponding to a causal source).
 
 The time domain animation
 \label{HL_ANR_left_t2} 
 \href{https://drive.google.com/open?id=0B7nf-pdU3Z4DcDI2dmN4U292VVE}{animation 6}
 presents the cross-sectional view along $z=0$ of the time-harmonic evolution of the field generated by the synthesized source and respectively the incoming plane wave $u_1= e^{i x k}e^{-ikct}$ in a near field region given by $(x,y)\in[-0.1,0.1]^2$. The multimedia file shows two animations: the top one describing the time propagation of the generated field and the bottom one describing the time propagation of the plane wave $u_1= e^{i x k}e^{-ikct}$. The color scheme in the movies is (truncated to 1 light yellow and -1 dark blue) with the antenna region removed from the simulations and with the black stripe representing amplitude values around $\approx 0.6$, and the white stripe representing amplitude values around $\approx -0.6$ respectively. As above note that there will be two black stripes and respectively two white stripes per period for the approximated plane wave. The animations show the good accuracy of the approximation of the outgoing plane wave $u_1= e^{-i x k}e^{-ikct}$ in region $D_1$.

\begin{figure}[!htbp] \centering
	\begin{subfigure}{\figsizeC}
		\includegraphics[width=\textwidth]{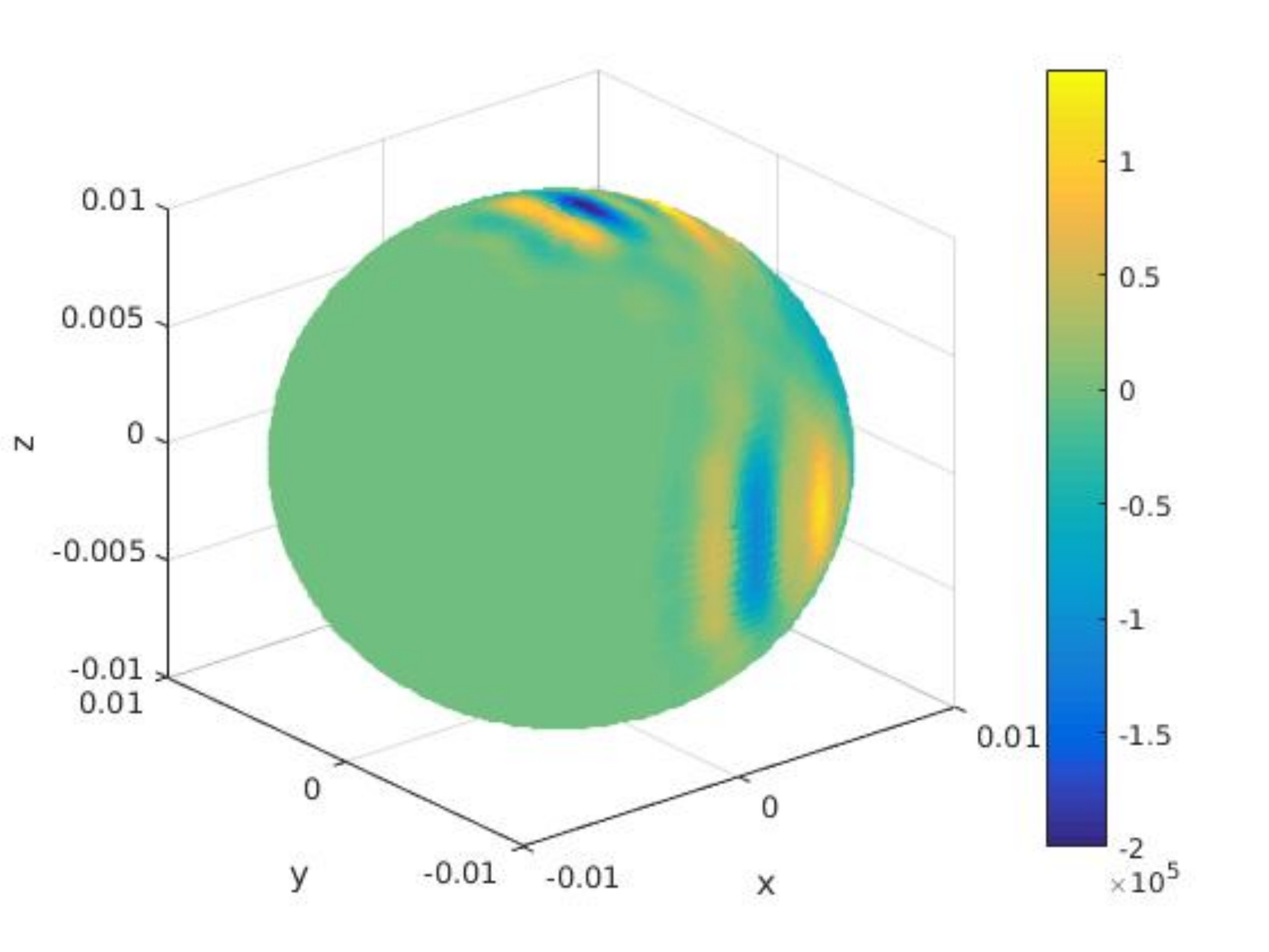}
		\label{fig:k10h30F10_ad_s1}
		\caption{side}
	\end{subfigure}
	\begin{subfigure}{\figsizeC}
		\includegraphics[width=\textwidth]{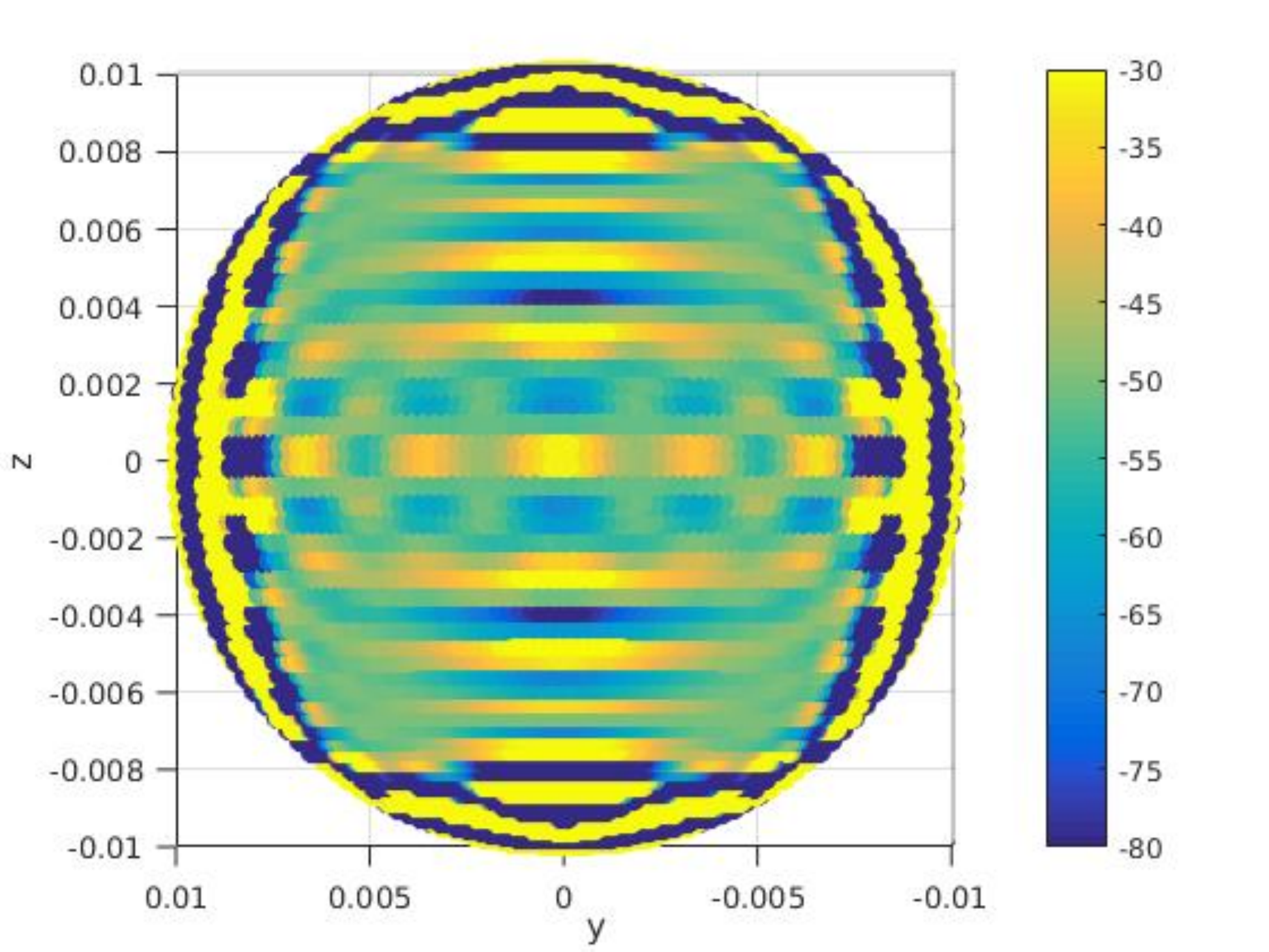}
		\label{fig:k10h30F10_ad_a2}
		\caption{front}
	\end{subfigure}
	\begin{subfigure}{\figsizeC}
		\includegraphics[width=\textwidth]{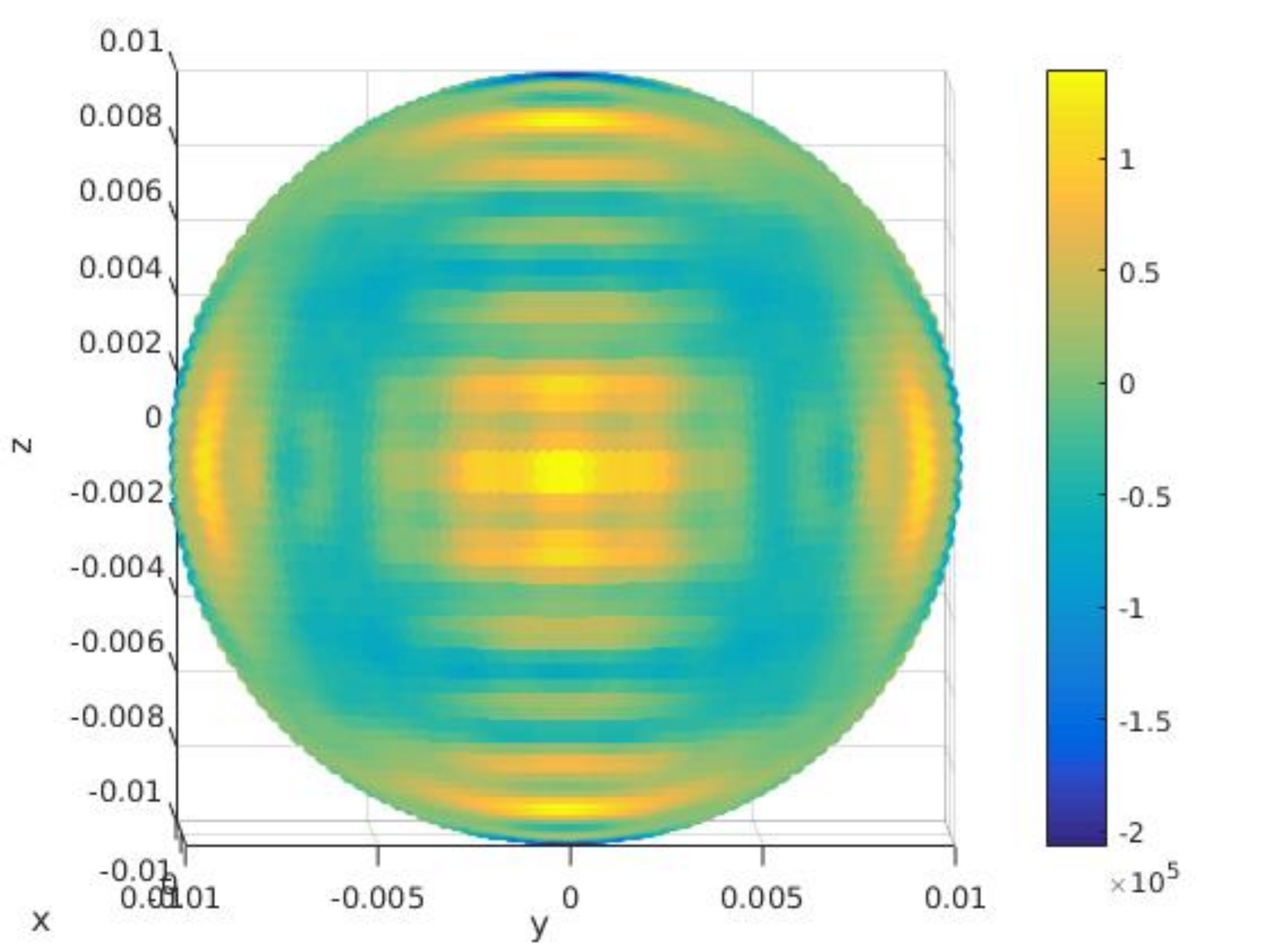}
		\label{fig:k10h30F10_ad_b1}
		\caption{back}
	\end{subfigure}
\caption{Density $w_\alpha$ on the antenna, with various color maps}
	\label{fig:k10h30F10_ad}
\end{figure}
Figure \ref{fig:k10h30F10_ad} shows the density $w_\alpha$ (see \eqref{lp}) on the boundary of the synthesized source. In the left plot in the figure we present the density values on the surface of the source viewed from a 3D side perspective and for better visualization we show two more plots in the figure; the center plot shows the density values on the part of the surface facing region $D_1$; the right plot of the figure presents the density values on the opposite pole of the source.

\section{Conclusions and Future Work} 
\label{III}
In the time-harmonic regime we described a unified framework where the possibility to characterize boundary source inputs for the control of acoustic fields in homogeneous infinite media (Section \ref{fs}) or finite depth homogeneous ocean environments (Section \ref{ho}) was theoretically established. We then presented a 3D optimization scheme (Section \ref{numerical}) based on the Method of Moments and Tikhonov regularization with Morozov discrepancy principle for the numerical characterization of boundary inputs (normal velocity or pressures) necessary on the boundary of the source to obtain the desired control effects. We then numerically discussed the performance of our scheme for the problem in homogeneous infinite  environments in the three distinct geometrical situations described at \eqref{geom-setting}. We did not show any explicit numerical simulations for the case of finite depth homogeneous ocean environments but, as shown in Section \ref{ho}, in this case the associated Green's function is computed explicitly in \cite{24} and is a continuous perturbation of the free space Helmholtz fundamental solution, hence the optimization scheme presented in the paper can be adapted to this case as well and this will be a part of our future work.

We believe that the results presented in this report may be relevant to: the problem of acoustic rendering or covert communications (Section \ref{O}); the problem of generating personal audio spots (Section \ref{I}); the question of acoustic protection where, paired with a time control loop for the detection of interrogating signals, a planar (or conformal) array of source elements similar to the one described Section \ref{II} is synthesized as a very weak radiator (thus unobservable to a far field measurement device) with a controllable near-field such that it nulls (by destructive interference in a region to the left of the array) an interrogating field coming from the left thus protecting a region behind the array, i.e., to its right side, (or within its convex hull if a conformal array is used). 

We mention also that the actual physical source boundary $\partial D_a$ does not need to be very smooth in general (Lipschitz suffice) and it must exists between the fictitious domain $D_{a'}$ and the boundary of $D_1$. In this regard, our investigations suggest also that, assuming more harmonics in the expansion of the density $w_\alpha$ we could achieve the same degree of control with region $D_1$ located further away from the source (giving thus more freedom in the choice of an actual physical boundary $\partial D_a$) or when one considers more then two regions of control. In this regard, in the spirit of \cite{22}, a detailed study of the sensitivity of the optimization scheme with respect to parameters such as, relative position of the control regions $D_1$ and $D_2$ and their distance from the source $D_a$, power budget and oscillatory character of the source input as well as acoustic intensity of the source will be presented in forthcoming reports. 

We observe that Figure \ref{fig:synthdensity}, Figure \ref{fig:double_ad} and Figure \ref{fig:k10h30F10_ad}  indicate the fact that, for all the cases considered in the paper, the synthesized source requires a very complex input on its boundary, i.e., with sub-areas characterized by small values and fast oscillations (e.g. , the part facing $D_1$ ) and other sub-areas characterized by very large values and slower oscillations (e.g., the pole opposite to $D_1$). In the context of linear approximation where everything can be scaled down appropriately, the results presented above in Figures \ref{fig:double_ue}, \ref{fig:double_ad}, \ref{fig:k10h30F10_ue} and \ref{fig:k10h30F10_ad} corroborated with the superposition principle suggest the possibility of approximating arbitrary given patterns (with small amplitudes) in each of the regions of interest. In this context, an investigation of similar optimization methods minimizing the discrepancy functional introduced at \eqref{eq:tikhonovobjective} but where instead of the penalty term $\alpha||w||_{L^2}^2$ one penalizes for the $TV$ or $L^1$ norm of the normal velocity (or pressure) (\eqref{2} or \eqref{2''} for free space and \eqref{4} or \eqref{4''} for homogeneous ocean environments), will be considered in a forthcoming report.

Last but not least, we believe that by using the superposition principle our strategy can be extended to the time-domain where different prescribed broadband signals could in principle be synthesised in disjoint near field regions. This, and, if needed, a direct time-domain analysis will be another important part of our future work.

	\section*{Acknowledgment}
	This work was supported by the Office of Naval Research under the award N00014-15-1-2462.

\end{document}